\documentclass{amsart}

\usepackage{amsmath}
\usepackage{xspace}
\usepackage{xypic}
\usepackage[psamsfonts]{amssymb}
\usepackage[latin1]{inputenc}
\usepackage{graphicx,color}
\usepackage{calc}
\usepackage[hypertex]{hyperref}
\usepackage{enumerate}

\theoremstyle{remark}
\newtheorem*{rema}{\bf Remark}

\theoremstyle{plain}
\newtheorem{theorem}{\bf Theorem}
\newtheorem*{theorem*}{\bf Theorem}
\newtheorem{lemma}{\bf Lemma}
\newtheorem{definition}{\bf Definition}
\newtheorem{proposition}{\bf Proposition}
\newtheorem{corollary}{\bf Corollary}

\newtheorem*{is}{\bf Result of Inou-Shishikura}

\newcounter{saveenum}

\def \eps{\varepsilon}
\def \epsilon{\varepsilon}
\def \e{\varepsilon}
\def\a{\alpha}

\def\d{\,{\rm d}}
\def\C{{\mathbb C}}
\def\D{{\mathbb D}}
\def\sF{{\mathcal IS}}
\def\H{{\mathbb H}}
\def\N{{\mathbb N}}
\def\P{{{\mathbb P}^1}}
\def\Q{{\mathbb Q}}
\def\R{{\mathbb R}}

\def\Z{{\mathbb Z}}

\def\att{\text{att}}
\def\rep{\text{rep}}

\def\mod{{\rm mod}}
\def\arg{{\rm arg}}
\def\re{{\rm Re}}
\def\im{{\rm Im}}

\def\area{{\rm area}}
\def\dens{{\rm dens}}

\def\PC{{\cal{PC}}}
\def\irrat{{\rm Irrat}}
\def\Pet{{\cal P}}
\def\Exp{{\rm Exp}}

\def\cal{\mathcal}

\def\ds{\displaystyle}

\def\Tq{\bigm|}
\def\To_#1{\mathrel{\mathop{\longrightarrow}\limits_{#1}}}
\def\eqdef{:=}

\begin{document}

\title{Quadratic Julia Sets with Positive Area.}
\author{Xavier Buff}
\author{Arnaud Chéritat}
\begin{abstract}We prove the existence of quadratic polynomials having
a Julia set with positive Lebesgue measure. We find such examples
with a Cremer fixed point, with a Siegel disk, or with infinitely
many satellite renormalizations.
\end{abstract}

\maketitle

\setcounter{tocdepth}{2}
\tableofcontents

\section*{Introduction}

Assume $P:\C\to \C$ is a polynomial of degree $2$. Its Julia set
$J(P)$ is a compact subset of $\C$ with empty interior. Fatou
suggested that one should apply to $J(P)$ the methods of
Borel-Lebesgue for the measure of sets.

It is known that the area (Lebesgue measure) of $J(P)$ is zero in
several cases including:
\begin{itemize}
\item if $P$ is hyperbolic;\footnote{Conjecturally, this is true for
a dense and open set of quadratic polynomials. If there were an open
set of non-hyperbolic quadratic polynomials, those would have a
Julia set of positive area (see \cite{mss}).}

\item if $P$ has a parabolic cycle (\cite{dh});

\item if $P$ is not infinitely renormalizable (\cite{l} or
\cite{s});

\item if $P$ has a (linearizable) indifferent cycle with
multiplier $e^{2i \pi \a}$ such that $\ds \a={\rm
a}_0+\cfrac{1}{{\rm a}_1+\cfrac{1}{{\rm a}_2+\ddots}}$ with $\log
{\rm a}_n = {\cal O}(\sqrt n)$ (\cite{pz}).\footnote{This is true
for almost every $\a\in \R/\Z$.}
\end{itemize}

Recently, we completed a program initiated by Douady with major
advances by the second author in \cite{c}: there exist quadratic
polynomials with a Cremer fixed point and a Julia set of positive
area. In this article, we present a slightly different approach (the
general ideas are essentially the same).

\begin{theorem}\label{theo_areacremer}
There exist quadratic polynomials which have a Cremer fixed point
and a Julia set of positive area.
\end{theorem}

We also have the following two results.

\begin{theorem}\label{theo_arealin}
There exist quadratic polynomials which have a Siegel disk and a
Julia set of positive area.
\end{theorem}

\begin{theorem}\label{theo_areainfren}
There exist infinitely satellite renormalizable quadratic
polynomials with a Julia set of positive area.
\end{theorem}

We will give a detailed proof of theorems \ref{theo_areacremer} and
\ref{theo_arealin}. We will only sketch the proof of theorem
\ref{theo_areainfren}.

The proofs are based on
\begin{itemize}
\item McMullen's results \cite{mcm} regarding the measurable density of the
filled-in Julia set near the boundary of a Siegel disk with bounded
type rotation number;

\item Chéritat's techniques of parabolic explosion \cite{c} and
Yoccoz's renormalization techniques \cite{y} to control the shape of
Siegel disks;

\item Inou and Shishikura's results \cite{is}
to control the post-critical sets of perturbations of polynomials
having an indifferent fixed point.
\end{itemize}

\section*{Acknowledgements}

We would like to thank Adrien Douady, John H. Hubbard, Hiroyuki
Inou, Curtis T. McMullen, Mitsuhiro Shishikura, Misha Yampolsky and
Jean-Christophe Yoccoz whose contributions were decisive in proving
these results.

\section{The Cremer case\label{sec_cremer}}

Let us introduce some notations.

\begin{definition}
For $\a\in \C$, we denote by $P_\a$ the quadratic polynomial
\[P_\a:z\mapsto e^{2i\pi \a}z+z^2.\]
We denote by $K_\a$ the filled-in Julia set of $P_\a$ and by $J_\a$
its Julia set.
\end{definition}

\subsection{Strategy of the proof}

The main gear is the following

\begin{proposition}\label{prop_keycremer0}
There exists a non empty set ${\cal S}$ of bounded type irrationals
such that: for all $\a\in {\cal S}$ and all $\eps>0$, there exists
$\a'\in {\cal S}$ with
\begin{itemize}
\item $|\a'-\a|<\eps$,

\item $P_{\a'}$ has a cycle in $D(0,\eps)\setminus \{0\}$ and

\item $\area(K_{\a'})\geq (1-\eps)\area (K_\a)$.
\end{itemize}
\end{proposition}

The proof of proposition \ref{prop_keycremer0} will occupy sections
\ref{sec_simplifyprop1} to \ref{sec_theproof}.

\begin{rema}
Since $\a\in {\cal S}$ has bounded type, $K_\a$ contains a Siegel
disk \cite{si} and thus, has positive area.
\end{rema}

\begin{figure}[htbp]
\centerline{ \scalebox{.45}{\includegraphics{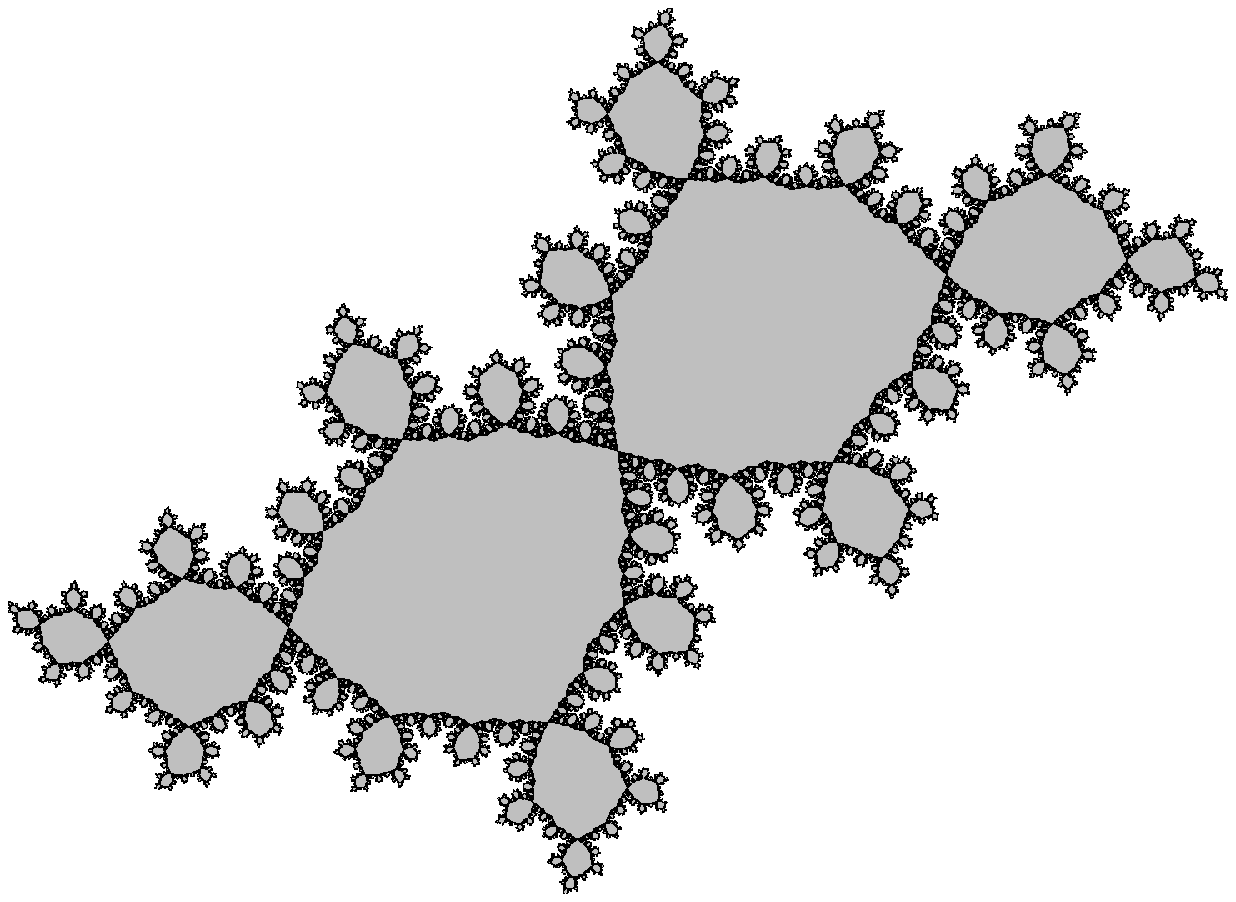}}
\scalebox{.45}{\includegraphics{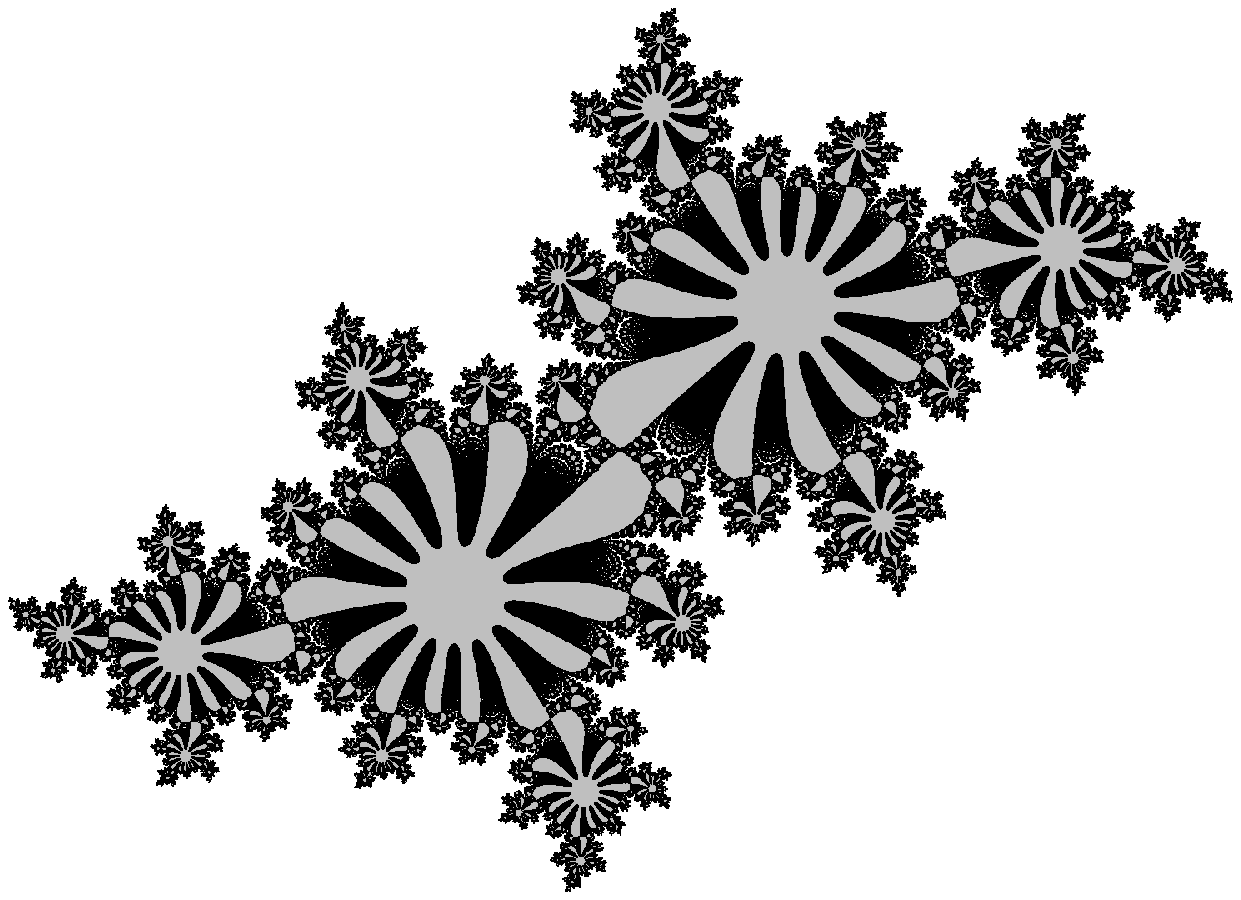}} } \caption{Two
filled-in Julia sets $K_{\a}$ and $K_{\a'}$, with $\a'$ a
well-chosen perturbation of $\a$ as in proposition
 \ref{prop_keycremer0}. This proposition
asserts that if $\a$ and $\a'$ are chosen carefully enough the loss
of measure from $K_{\a}$ to $K_{\a'}$ is small.}
\end{figure}

\begin{figure}[htbp]
\centerline{%
\setlength{\fboxsep}{0cm}%
\framebox{\scalebox{1.8}{\includegraphics{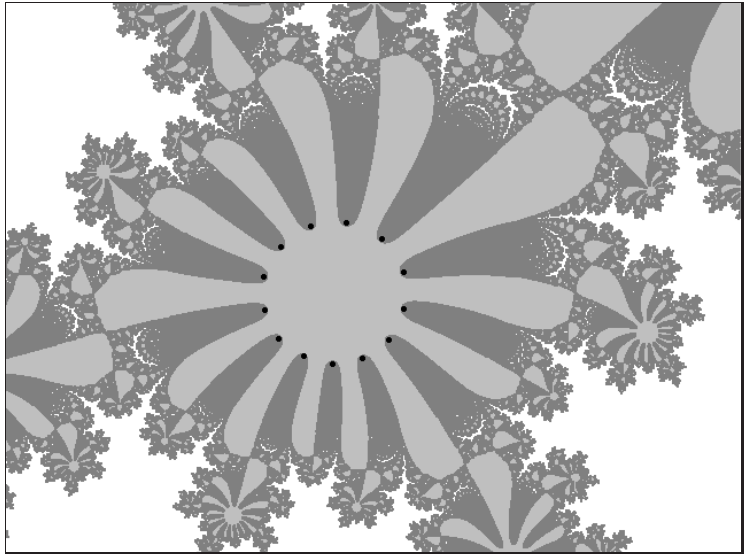}}}}%
\caption{A zoom on $K_{\a'}$ near its linearizable fixed point. The
small cycle is highlighted.}
\end{figure}

\begin{rema}
We do not know what is the largest set ${\cal S}$ for which
proposition \ref{prop_keycremer0} holds. It might be the set of all
bounded type irrationals.
\end{rema}

\begin{proposition}[Douady]\label{prop_semicontmeasure}
The function $\a\in \C\mapsto \area(K_\a)\in \left[0,+\infty\right[$
is upper semi-continuous.
\end{proposition}

\begin{proof} Assume $\a_n\to \a$. By \cite{d}, for any neighborhood
$V$ of $K_\a$, we have $K_{\a_n}\subset V$ for $n$ large enough.
According to the theory of Lebsegue measure, $\area(K_\a)$ is the
infimum of the area the open sets containing $K_\a$. Thus,
\[\area(K_\a)\geq \limsup_{n\to +\infty} \area(K_{\a_n}).\]
\end{proof}

\noindent{\it Proof of theorem \ref{theo_areacremer} assuming
proposition \ref{prop_keycremer0}.}~We choose a sequence of real
numbers $\e_n$ in $(0,1)$ such that $\prod(1-\e_n)>0$. We construct
inductively a sequence $\theta_n\in {\cal S}$ such that for all
$n\geq 1$
\begin{itemize}
\item $P_{\theta_n}$ has a cycle in $D(0,1/n)\setminus \{0\}$,

\item $\area(K_{\theta_{n}}) \geq (1-\eps_n)\area(K_{\theta_{n-1}})$.
\end{itemize}
Every polynomial $P_\theta$ with $\theta$ sufficiently close to
$\theta_n$ has a cycle in $D(0,1/n)\setminus \{0\}$. By choosing
$\theta_{n}$ sufficiently close to $\theta_{n-1}$ at each step, we
guarantee that
\begin{itemize}
\item the sequence $(\theta_n)$ is a
Cauchy sequence that converges to a limit $\theta$,

\item for all $n\geq 1$, $P_\theta$ has a cycle in $D(0,1/n)\setminus
\{0\}$.
\end{itemize}
So, the polynomial $P_\theta$ has small cycles and thus is a Cremer
polynomial. In that case, $J_\theta=K_\theta$. By proposition
\ref{prop_semicontmeasure}:
\[\area(J_\theta)=\area(K_\theta) \geq \limsup_{n\to +\infty}
\area(K_{\theta_n}) \geq \area(K_{\theta_0})\cdot\prod_{n\geq 1}
(1-\eps_n) >0.\]
\qed\par\medskip

\subsection{A stronger version of proposition
\ref{prop_keycremer0}\label{sec_simplifyprop1}}

For a finite or infinite sequence of integers, we will use the
following continued fraction notation:
\[[{\rm a}_0,{\rm a}_1,{\rm a}_2,\ldots]\eqdef {\rm a}_0+\cfrac{1}{{\rm a}_1+
\cfrac{1}{{\rm a}_2+\ddots}}.\] For $\a\in \R$, we will denote by
$\lfloor \a\rfloor$ the integral part of $\a$.

\begin{definition}\label{def_sn}
If $N\geq 1$ is an integer, we set
\[{\cal S}_N \eqdef \bigl\{\a=[{\rm a}_0,{\rm a}_1,{\rm a}_2,\ldots]\in \R\setminus
\Q\Tq({\rm a}_k)\text{ is bounded  and }{\rm a}_k\geq N\text{ for
all }k\geq 1\bigr\}.\]
\end{definition}

Note that ${\cal S}_{N+1}\subset {\cal S}_N\subset \cdots\subset
{\cal S}_1$ and ${\cal S}_1$ is the set of bounded type irrationals.
If $\a\in {\cal S}_1$, the polynomial $P_\a$ has a Siegel disk
bounded by a quasicircle containing the critical point (see
\cite{d1}, \cite{he}, \cite{sw}). In particular, the post-critical
set of $P_\a$ is contained in the boundary of the Siegel disk.

Proposition \ref{prop_keycremer0} is an immediate consequence of the
following proposition.

\begin{proposition}\label{prop_keycremer}
If $N$ is sufficiently large then the following holds.
\footnote{The choice of $N$ will be specified in equation \ref{eq_chooseN}}\\
Assume $\a\in {\cal S}_N$, choose a sequence $(A_n)$ such that
\[\sqrt[q_n]{A_n}\To_{n\to +\infty}+\infty\quad \text{and}\quad
\sqrt[q_n]{\log A_n}\To_{n\to +\infty}1.\footnote{For example, one
can choose
  $A_n\eqdef q_n^{q_n}$. However, we think that the
  proposition holds for more general sequences $(\alpha_n)$ for
  instance as soon as $\sqrt[q_n]{A_n}\to +\infty$.}\]
Set
\[\a_n\eqdef [{\rm a}_0,{\rm a}_1,\ldots, {\rm a}_n, A_n,
N,N,N,\ldots].\] Then, for all $\eps>0$, if $n$ is sufficiently
large,
\begin{itemize}
\item $P_{\a_n}$ has a cycle in $D(0,\eps)\setminus \{0\}$ and

\item $\area(K_{\a_n})\geq (1-\eps)\area (K_\a)$.
\end{itemize}
\end{proposition}

The rest of section \ref{sec_cremer} is devoted to the proof of
proposition \ref{prop_keycremer}. In the sequel, unless otherwise
specified,
\begin{itemize}
\item $\a$ is an irrational number of
bounded type,

\item $p_k/q_k$ are the approximants to $\a$ given by the continued
fraction algorithm and

\item $(\a_n)$ is a sequence converging to $\a$, defined as in
proposition \ref{prop_keycremer}.
\end{itemize}
Note that for $k\leq n$, the approximants $p_k/q_k$ are the same for
$\a$ and for $\a_n$. The polynomial $P_\a$ (resp. $P_{\a_n}$) has a
Siegel disk $\Delta$ (resp. $\Delta_n$). We let $r$ (resp. $r_n$) be
the conformal radius of $\Delta$ (resp. $\Delta_n$) at $0$ and we
let $\phi:D(0,r)\to \Delta$ (resp. $\phi_n:D(0,r_n)\to \Delta_n$) be
the conformal isomorphism which maps $0$ to $0$ with derivative $1$.

\subsection{The control of the cycle}

We first recall results of \cite{c} (see also \cite{bc0} Props. 1
and 2), which we reformulate as follows.

The first proposition asserts that as $\theta$ varies in the disk
$D(p/q,1/q^3)$, the polynomial $P_\theta$ has a cycle of period $q$
which depends holomorphically on $\sqrt[q]{\theta-p/q}$ and
coalesces at $z=0$ when $\theta=p/q$.

\begin{proposition}\label{prop_chi1}
For each rational number $p/q$ (with $p$ and $q$ coprime), there
exists a holomorphic function
\[\chi : D(0,1/q^{3/q}) \to \C\] with the following
properties.
\begin{enumerate}
\item $\chi(0)=0$.

\item $\chi'(0)\neq 0$.

\item If $\delta \in D(0,1/q^{3/q})\setminus \{0\}$, then $\chi(\delta)\neq
0$.

\item If $\delta \in D(0,1/q^{3/q})\setminus \{0\}$ and if we set $\zeta\eqdef e^{2i\pi p/q}$ and
$\ds \theta \eqdef \frac{p}{q} + \delta^{q_k}$, then,
\(\Big<\chi(\delta), \chi(\zeta \delta), \ldots,
\chi(\zeta^{q-1}\delta) \Big> \) forms a cycle of period $q$ of
$P_\theta$. In particular,
\[\forall \delta\in D(0,1/q^{3/q}),\qquad
\chi(\zeta \delta) = P_\theta \bigl(\chi(\delta)\bigr).
\]
\end{enumerate}
\end{proposition}

A function $\chi:D(0,1/q^{3/q})\to \C$ as in proposition
\ref{prop_chi1} is called an {\em explosion function} at $p/q$. Such
a function is not unique. However, if $\chi_1$ and $\chi_2$ are two
explosions functions at $p/q$, they are related by
$\chi_1(\delta)=\chi_2(e^{2i\pi kp/q}\delta)$ for some integer $k\in
\Z$.

The second proposition studies how the explosion functions behave as
$p/q$ ranges in the set of approximants of an irrational number $\a$
such that $P_\a$ has a Siegel disk.

\begin{proposition}\label{prop_chi2}
Assume $\a\in \R\setminus \Q$ is an irrational number such that
$P_\a$ has a Siegel disk $\Delta$. Let $p_k/q_k$ be the approximants
to $\a$. Let $r$ be the conformal radius of $\Delta$ at $0$ and let
$\phi:D(0,r)\to \Delta$ be the isomorphism which sends $0$ to $0$
with derivative $1$. For $k\geq 1$, let $\chi_k$ be an explosion
function at $p_k/q_k$ and set $\lambda_k\eqdef \chi_k'(0)$. Then,
\begin{enumerate}
\item $\ds |\lambda_k|\To_{k\to +\infty} r$ and

\item the sequence of maps
$\psi_k:\delta\mapsto\chi_k(\delta/\lambda_k)$ converges uniformly
on every compact subset of $D(0,r)$ to $\phi:D(0,r)\to \Delta$.
\end{enumerate}
\end{proposition}

\begin{corollary}
Let $(\a_n)$ be the sequence defined in proposition
\ref{prop_keycremer}. Then, for all $\eps>0$, if $n$ is sufficiently
large, $P_{\a_n}$ has a cycle in $D(0,\eps)\setminus \{0\}$.
\end{corollary}

\begin{proof} Let $\chi_n$ be an explosion at $p_n/q_n$ and let $C_n$
be the set of $q_n$-th roots of
\[\a_n-\frac{p_n}{q_n} = \frac{(-1)^n}{q_n (q_n A'_n + q_{n-1})}\quad{\rm with}\quad
A'_n \eqdef [A_n,N,N,N,\ldots].\] Since $\ds
\sqrt[q_n]{A'_n}\To_{n\to +\infty}+\infty$, for $n$ large enough,
the set $C_n$ is contained in an arbitrarily small neighborhood of
$0$ and $\chi_n(C_n)$ is a cycle of $P_{\a_n}$ contained in an
arbitrarily small neighborhood of $0$.
\end{proof}

\subsection{Perturbed Siegel disks\label{sec_pertsigeldisk}}

\begin{definition}
If $U$ and $X$ are measurable subsets of $\C$, with
$0<\area(U)<+\infty$, we use the notation
\[\dens_U(X)\eqdef \frac{\area(U\cap X)}{\area (U)}.\]
\end{definition}

In the whole section, $\a$ is a Bruno number, $p_n/q_n$ are its
approximants, and $\chi_n:D_n\eqdef D(0,1/q_n^{3/q_n})\to \C$ are
explosion functions at $p_n/q_n$.

\begin{figure}[htbp]
\centerline{\scalebox{.7}{\includegraphics{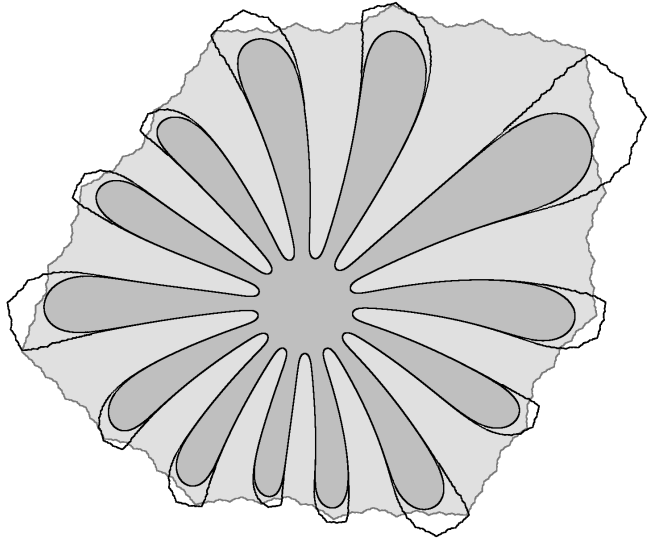}}}
\caption{Illustration of proposition \ref{theo_disquesdigites} for
$\a= \theta= [0,1,1,\ldots]$, $n=7$ and $A_n=10^{10}$. We see the
Siegel disk $\Delta$ of $P_\a$ (light grey), the Siegel disk
$\Delta'_n$ of the restriction of $P_{\a_n}$ to $\Delta$ (dark grey)
and the boundary of the Siegel disk of $P_{\a_n}$.\label{fig_amibe}}
\end{figure}

\begin{proposition}[see figure
\ref{fig_amibe}]\label{theo_disquesdigites} Assume $\a\eqdef [{\rm
a}_0,{\rm a}_1,\ldots]$ and $\theta\eqdef[0,{\rm t}_1,\ldots]$ are
Brjuno numbers and let $p_n/q_n$ be the approximants to $\alpha$.
Assume
\[\a_n \eqdef [{\rm a}_0,{\rm a}_1,\ldots,{\rm a}_n, A_n,{\rm
t}_1,{\rm t}_2,\ldots]\] with $(A_n)$ a sequence of positive
integers such that
\begin{equation}\label{eq_logAnsubexp}
\limsup_{n\to +\infty} \sqrt[q_n]{\log A_n}\leq 1.\footnote{We think
that the
  condition $\ds \limsup \sqrt[q_n]{\log A_n}\leq 1$ is not necessary.}
\end{equation}
Let $\Delta$ be the Siegel disk of $P_\a$ and $\Delta'_n$ the Siegel
disk of the restriction of $P_{\a_n}$ to
$\Delta$.\footnote{$\Delta'_n$ is the largest connected open subset
of $\Delta$ containing $0$, on which $P_{\a_n}$ is conjugate to a
rotation. It is contained in the Siegel disk of $P_{\a_n}$} For all
non empty open set $U\subset \Delta$,
\[\liminf_{n\to +\infty}\dens_U(\Delta'_n)\geq \frac{1}{2}.\]
\end{proposition}

\begin{proof} Set
\[\e_n \eqdef \a_n-\frac{p_n}{q_n} =
\frac{(-1)^n}{q_n^2(A_n+\theta)+q_nq_{n-1}}.\] Note that
\[\sqrt[q_n]{|\e_n|}\underset{n\to +\infty}\sim
\frac{1}{\sqrt[q_n]{A_n}}.\]
For $\rho<1$, define
\[X_n(\rho)\eqdef\left\{z \in \C ~;~
\frac{z^{q_n}}{z^{q_n}-\e_n} \in
D(0,s_n)\right\}\quad\text{with}\quad s_n\eqdef
\frac{\rho^{q_n}}{\rho^{q_n}+|\e_n|}.\]

\begin{figure}[htbp]
\begin{picture}(150,150)(0,0)%
\put(0,0){\scalebox{.4}{\includegraphics{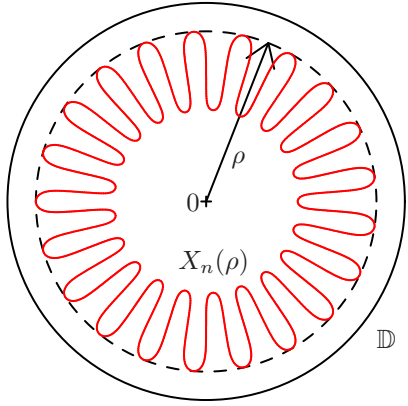}}}%
\put(68,72){$0$}%
\put(85,90){$\rho$}%
\put(140,20){$\D$}%
\put(65,50){$X_n(\rho)$}
\end{picture}
\caption{The boundary of a set $X_n(\rho)$.}
\end{figure}
This domain is star-like with respect to  $0$ and avoids the
$q_n$-th roots of $\e_n$.\footnote{It is the preimage by the map
$z\mapsto z^{q_n}$ of a disk which is not centered at $0$, contains
$0$ but not $\e_n$.} It is contained but not relatively compact in
$D(0,\rho)$. For all non empty open set $U$ contained in
$D(0,\rho)$,
\[\liminf_{n\to +\infty} \dens_U\bigl(X_n(\rho)\bigr)\geq \frac{1}{2}.\]
Since the limit values of the sequence $(\chi_n:D_n\to \C)$ are
isomorphisms $\chi:\D\to \Delta$, proposition
\ref{theo_disquesdigites} is a corollary of proposition
\ref{prop_SiegelXn} below.
\end{proof}

\begin{proposition}\label{prop_SiegelXn}
Under the same assumptions as in proposition
\ref{theo_disquesdigites}, for all $\rho<1$, if $n$ is large enough,
the Siegel disk $\Delta'_n$ contains $\chi_n\bigl(X_n(\rho)\bigr)$.
\end{proposition}

\proof We will proceed by contradiction. Assume there exist $\rho<1$
and an increasing sequence of integers $n_k$ such that
$\chi_{n_k}\bigl(X_{n_k}(\rho)\bigr)$ is not contained in
$\Delta'_{n_k}$. Extracting a subsequence, we may assume
\[A_{n_k}^{1/q_{n_k}}\to A\in [1,+\infty].\]
To simplify notations, we will drop the index $k$.

\begin{itemize}
\item
Assume $A=1$. Then, any compact $K\subset \Delta$ is contained in
$\Delta'_{n}$ for $n$ large enough (for a proof, see for example in
\cite{abc}, proposition 2, the remark following proposition 2 and
theorem 3). Note that $X_n(\rho)\subset D(0,\rho)$ and the limit
values of the sequence $(\chi_n:D_n\to \C)$ are isomorphisms
$\chi:\D\to \Delta$. It follows that for $n$ large enough,
\[
\chi_{n}\bigl(X_{n}(\rho)\bigr)\subset
\chi_{n}\bigl(D(0,\rho)\bigr)\subset \chi\bigl(D(0,\sqrt\rho)\bigr)
\subset \Delta'_{n}.\] This contradicts our assumption.

\item
Assume $A>1$. Without loss of generality, increasing $\rho$ if
necessary, we may assume that $\rho>1/A$. We will show that for
$\rho<\rho'<1$, if $n$ is large enough, the orbit under iteration of
$P_{\a_{n}}$ of any point $z\in \chi_{n}\bigl(X_{n}(\rho)\bigr)$
remains in $\chi_{n}\bigl(D(0,\rho')\bigr)\subset \Delta$. This will
show that $\chi_{n}\bigl(X_{n}(\rho)\bigr)\subset \Delta'_n$,
completing the proof of proposition \ref{prop_SiegelXn}.

Since the limit values of the sequence $\chi_n:D_n\to \C$ are
isomorphisms $\chi:\D\to \Delta$, there is a sequence $r'_n$ tending
to $1$ such that $\chi_n$ is univalent on $D'_n\eqdef D(0,r'_n)$ and
the domain of the map
\[f_n\eqdef \bigl(\chi_n|_{D'_n}\bigr)^{-1}\circ P_{\a_n}\circ \chi_n|_{D'_n}\]
eventually contains any compact subset of $\D$. So, proposition
\ref{prop_SiegelXn} is a corollary of proposition
\ref{prop_SiegelXn}' below.
\end{itemize}
\qed\par\medskip

\noindent{\bf Proposition \ref{prop_SiegelXn}'.} {\em Assume \[0\leq
\frac{1}{A}<\rho<\rho'<1.\] If $n$ is large enough, the orbit under
iteration of $f_n$ of any point $z\in X_n(\rho)$ remains in
$D(0,\rho')$.}
\medskip

The rest of section \ref{sec_pertsigeldisk}, is devoted to the proof
of proposition \ref{prop_SiegelXn}'.

\subsubsection{A vector field}

It is not enough to compare the dynamics of $f_n$ with the dynamics
of a rotation. Instead, we will compare it with the (real) dynamics
of the polynomial vector field
\[\xi_n=\xi_n(z)\frac{\partial}{\partial z}
\eqdef 2i\pi q_n z(\e_n-z^{q_n})\frac{\partial}{\partial z}.\] As we
shall see later, the time-1 map of $\xi_n$ very well approximates
$f_n^{\circ q_n}$ (the coefficient $2\pi q_n$ has been chosen so
that their derivatives coincide at $0$). For simplicity we will
assume that $n$ is even in which case $\eps_n>0$.

Note that the polynomial vector field $\xi_n$ is tangent to the
boundary of $X_n(\rho)$, which is therefore invariant by the (real)
dynamics of $\xi_n$.
\begin{figure}[htbp]
\centerline{\scalebox{.4}{\includegraphics{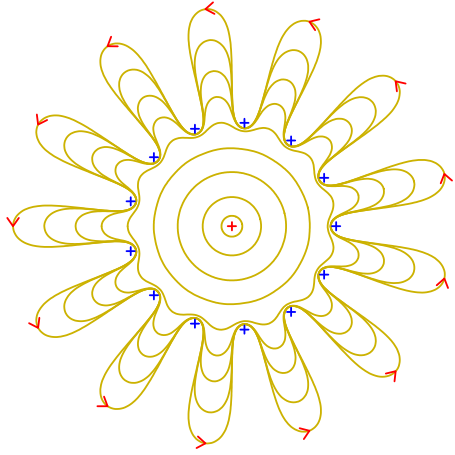}}} \caption{Some
real trajectories for the vector field $\xi_n$; zeroes of the vector
field are shown.}
\end{figure}

In order to compare the dynamics of $f_n$ to that of $\xi_n$, we
will work in a coordinate that straightens the vector field $\xi_n$.
Let us first consider the open set
\[\Omega_n\eqdef \left\{ z\in \C~;~
\frac{z^{q_n}}{z^{q_n}-\eps_n}\in \D\right\}\] which is invariant by
the real flow of the vector field $\xi_n$.

\begin{figure}[htbp]
\centerline{\scalebox{.4}{\includegraphics{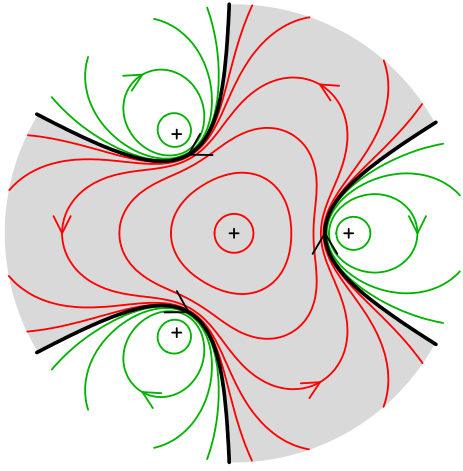}}}
\caption{An example of open set $\Omega_n$ for $q_n=3$. It is
bounded by the black curves. Some trajectories of the vector field
$\xi_n$ (red in $\Omega_n$ and green outside).}
\end{figure}

The map
\[z\mapsto \frac{z^{q_n}}{z^{q_n}-\e_n}:\Omega_n\to \D\]
is a ramified covering of degree $q_n$, ramified at $0$. Thus, there
is an isomorphism $\psi_n:\Omega_n\to \D$ such that
\[\bigl(\psi_n(z)\bigr)^{q_n}=\frac{z^{q_n}}{z^{q_n}-\e_n}.\]
We note $\phi_n:\D\to \Omega_n$ its inverse and $\pi_n:\H\to
\Omega_n\setminus\{0\}$ ($\H$ is the upper half-plane) the universal
covering given by
\[\pi_n(Z) \eqdef \phi_n\bigl(e^{2i\pi q_n \e_n Z}\bigr).\]
Then,
\[\pi_n^* \xi_n = \frac{\partial}{\partial z}.\]

For $r<1$, we have $X_n(r)\subset \Omega_n$ and the preimage of
$X_n(r)$ is the half-plane
\[\H_n(r) \eqdef \bigl\{Z\in \C~;~ \im(Z) >
\tau_n(r)\bigr\}\quad\text{with}\quad \tau_n(r)\eqdef \frac{1}{2\pi
q_n^2 \e_n} \log \left(1+\frac{\e_n}{r^{q_n}}\right).\] The map
$\pi_n:\H_n(r)\to X_n(r)\setminus \{0\}$ is a universal covering.

\begin{rema}
Note that $\tau_n(r)$ increases exponentially fast with respect to
$q_n$. More precisely,
\[\sqrt[q_n]{\tau_n(r)} \To_{n\to +\infty} \frac{1}{r}.\]
\end{rema}

\subsubsection{Working in the coordinate straightening the vector field}

\begin{definition}
We say that a sequence $(B_n)$ is {\em sub-exponential with respect
to} $q_n$ if
\[\limsup_{n\to +\infty} \sqrt[q_n]{|B_n|}\leq 1.\]
\end{definition}

\begin{proposition}\label{prop_liftFnGn}
Assume $r<1$. If $n$ is large enough, there exist holomorphic maps
$F_n:\H_n(r)\to \H$ and $G_n:\H_n(r)\to \H$ such that
\begin{itemize}
\item $\pi_n$ semi-conjugates $F_n$ to $f_n^{\circ q_n}$ and $G_n$ to
$f_n^{\circ q_{n-1}}$:
\[\pi_n \circ F_n = f_n^{q_n} \circ \pi_n
\quad\text{and}\quad \pi_n \circ G_n = f_n^{q_{n-1}} \circ \pi_n,\]

\item $F_n$ and $G_n$ are periodic of period $1/(q_n\eps_n)$ and

\item as $\im(Z)\to +\infty$, we have
\[F_n(Z) =Z+ 1 + o(1) \quad\text{and}\quad
G_n(Z)= Z-(A_n+\theta)+ o(1).\]
\end{itemize}
In addition, the sequences
\[\sup_{Z\in \H_n(r)} \bigl|F_n(Z)-Z-1\bigr|\quad\text{and}\quad
\sup_{Z\in \H_n(r)} \bigl|G_n(Z)-Z+A_n+\theta\bigr|\] are
sub-exponential with respect to $q_n$.
\end{proposition}

\begin{proof} We will use the following theorem of Jellouli (see
\cite{j1} or \cite{j2} Theorem 1) to show that the domains of
$f_n^{\circ q_n}$ and $f_n^{\circ q_{n-1}}$ eventually contain any
compact subset of $\D$.

\begin{theorem*}[Jellouli]
Assume $P_\a$ has a Siegel disk $\Delta$ and let $\chi:\D\to \Delta$
be a linearizing isomorphism. For $r<1$, set $\Delta(r)\eqdef
\chi\bigl(D(0,r)\bigr)$. Assume $\a_n\in \R$ and $b_n\in \N$ are
such that $b_n\cdot |\alpha_n-\alpha|=o(1)$.\footnote{In
  fact, Jellouli's theorem
  is stated for
the sequence $\alpha_n=p_n/q_n$ and $b_n = o(q_n q_{n+1})$ but the
adaptation to $b_n\cdot|\alpha_n-\alpha|=o(1)$ is straightforward.}
For all $r'_1<r'_2<1$, if $n$ is sufficiently large,
\[\Delta(r'_1)\subset \bigl\{z\in \Delta(r'_2)~;~
\forall j\leq b_n,~P_{\a_n}^{\circ j}(z)\in \Delta(r'_2)\bigr\}.\]
\end{theorem*}

\begin{corollary}\label{coro_jellouli}
For all $r_1<r_2<1$, if $n$ is sufficiently large, then for all
$z\in D(0,r_1)$ and for all $j\leq q_n$, we have $f_n^{\circ
j}(z)\in D(0,r_2)$.
\end{corollary}

\begin{proof} Choose $r'_1$ and $r'_2$ such that $r_1<r'_1<r'_2<r_2$.
Let $\chi:\D\to \Delta$ be a linearizing isomorphism of $P_\a$. Set
\[\Delta(r'_1)\eqdef \chi\bigl(D(0,r'_1)\bigr)\quad\text{and}\quad
\Delta(r'_2)\eqdef \chi\bigl(D(0,r'_2)\bigr).\] Since limit values
of the sequence $\chi_n:D'_n\to \C$ are linearizing isomorphisms
$\chi:\D\to \Delta$, for $n$ sufficiently large,
\[\chi_n\bigl(D(0,r_1)\bigr)\subset \Delta(r'_1)\subset
\Delta(r'_2)\subset \chi_n\bigl(D(0,r_2)\bigr).\] It is therefore
enough to show that for $n$ large enough,
\[\Delta(r'_1)\subset \bigl\{z\in \Delta(r'_2)~;~
\forall j\leq q_n,~P_{\a_n}^{\circ j}(z)\in \Delta(r'_2)\bigr\}.\]
This is Jellouli's theorem with $b_n=q_n$ since
\[q_n|\alpha_n-\alpha| \underset{n\to +\infty}\sim q_n
\left|\frac{p_n}{q_n}-\alpha\right| \underset{n\to +\infty}= o(1).\]
\end{proof}

In particular, for $r<1$, if $n$ is large enough, then $f_n^{\circ
q_n}$ and $f_n^{\circ q_{n-1}}$ are defined on $X_n(r)$. In order to
lift them via $\pi_n$ as required, it is enough to show that if $n$
is large enough, then
\[\forall z\in X_n(r)\setminus \{0\},
\quad f_n^{q_n}(z)\in \Omega_n\setminus \{0\}\quad\text{and}\quad
f_n^{q_{n-1}}(z)\in \Omega_n\setminus\{0\}.\] The periodicity of
$F_n$ and $G_n$ then follows from
\[\pi_n\left(Z+\frac{1}{q_n\eps_n}\right) = \pi_n(Z)\]
and the behavior as $\im(Z)\to +\infty$ follows by computing the
derivatives of $f_n^{\circ q_n}$ and $f_n^{\circ q_{n-1}}$ at $0$.

Lemma \ref{lem_approxqn} below asserts that $f_n^{\circ q_n}$ is
very close to the identity and bounds the difference.

\begin{lemma}\label{lem_approxqn}
There exist a holomorphic function $g_n$, defined on the same set as
$f_n^{\circ q_n}$, such that
\[f_n^{\circ q_n}(z)  = z + \xi_n(z)\cdot g_n(z).\]
For all $r<1$, the sequence $\ds \sup_{D(0,r)} |g_n|$ is
sub-exponential with respect to $q_n$.
\end{lemma}

\begin{proof} The map $f_n^{\circ q_n}$ fixes $0$ and the $q_n$-th
roots of $\eps_n$. This shows that $f_n^{\circ q_n}$ can be written
as prescribed. To prove the estimate on the modulus of $g_n$, note
that $f_n^{\circ q_{n}}$ takes its values in $\D$ and thus,
$\bigl|\xi_n(z)\cdot g_n(z)\bigr|\leq 2$. Choose a sequence $r_n\in
\left]0,1\right[$ tending to $1$ so that $g_n$ is defined on
$D(0,r_n)$. By the maximum modulus principle, if $n$ is large enough
so that $r_n>\max(r,1/A)$, we have
\[\sup_{|z|\leq r} \bigl|g_n(z)\bigr|\leq
\sup_{|z|\leq r_n} \bigl|g_n(z)\bigr| \leq B_n\eqdef \sup_{|z|=r_n}
\frac{2}{\bigl|\xi_n(z)\bigr|}.\] As $n\to +\infty$,
\[\inf_{|z|=r_n}\bigl|\xi_n(z)\bigr|\sim 2\pi q_n r_n^{1+q_n}
\quad\text{and thus}\quad \sqrt[q_n]{B_n}\sim r_n \to 1.\]
\end{proof}

Recall that we assume $n$ even, in which case
\[\e_n>0\quad\text{and}\quad q_{n-1}\cdot\frac{p_n}{q_n} = -\frac{1}{q_n}~\mod ~(1).\]
Lemma \ref{lem_approxqnmoins1} below asserts that $f_n^{\circ
q_{n-1}}$ is very close to the rotation of angle $-1/q_n$ and bounds
the difference.

\begin{lemma}\label{lem_approxqnmoins1}
There exists a holomorphic function $h_n$, defined on the same set
as $f_n^{\circ q_{n-1}}$, such that
\[e^{2i \pi/q_n} f_n^{\circ q_{n-1}}(z) =  z+ \xi_n(z)\cdot h_n(z).\]
For all $r<1$, the sequence $\ds \sup_{D(0,r)} |h_n|$ is
sub-exponential with respect to $q_n$.
\end{lemma}

\begin{proof} The map $f_n$ coincides with the rotation of angle
$p_n/q_n$ on the set of $q_n$-th roots of $\eps_n$ and
$q_{n-1}\cdot(p_n/q_n) = -1/q_n ~\mod(1)$. Thus, $e^{2i\pi
/q_n}f_n^{\circ q_{n-1}}(z)$ fixes $0$ and the $q_n$-th roots of
$\eps_n$. This shows that $e^{2i\pi /q_n}f_n^{\circ q_{n-1}}$ can be
written as prescribed. The same method as in
lemma~\ref{lem_approxqn} yields the bound on $h_n$.
\end{proof}

\noindent{\it Proof of proposition \ref{prop_liftFnGn}, continued}.
Now, given $r<1$, set
\[R_n\eqdef \min\left(\frac{1}{q_n\eps_n},\tau_n(r)\right).\]
Note that
\[\sqrt[q_n]{R_n}\To_{n\to +\infty} \min
\left(A,\frac{1}{r}\right).\] Hence, $R_n$ increases exponentially
fast with respect to $q_n$.

For all $n$ and all $Z\in \H_n(r)$, the map $\pi_n$ is univalent on
$D(Z,R_n)$ and takes its values in $\Omega_n\setminus \{0\}$. By
Koebe $1/4$-theorem, its image contains a disk centered at $z\eqdef
\pi_n(Z)$ with radius
\[\pi_n'(Z)\cdot \frac{R_n}{4} = \xi_n(z)\cdot \frac{R_n}{4}.\]
In particular, if the sequence $(B_n)$ is sub-exponential  with
respect to $q_n$ and if $n$ is large enough so that $B_n\leq R_n/4$,
we have
\[\forall z\in X_n(r),\quad D\bigl(z,\xi_n(z)\cdot B_n\bigr)\subset
\Omega_n\setminus \{0\}.\] Therefore, it follows from lemmas
\ref{lem_approxqn} and \ref{lem_approxqnmoins1} that for all $r<1$,
if $n$ is large enough, then
\[\forall z\in X_n(r)\setminus \{0\},
\quad f_n^{q_n}(z)\in \Omega_n\setminus \{0\}\quad\text{and}\quad
f_n^{q_{n-1}}(z)\in \Omega_n\setminus\{0\}.\]

Lemmas \ref{lem_approxqn} and \ref{lem_approxqnmoins1} and Koebe
distortion theorem applied to $\pi_n:D\bigl(Z,R_n\bigr)\to\C$ imply
that the sequences
\[\sup_{Z\in \H_n(r)} \bigl|F_n(Z)-Z-1\bigr|\quad\text{and}\quad
\sup_{Z\in \H_n(r)} \bigl|G_n(Z)-Z+A_n+\theta\bigr|\] are
sub-exponential  with respect to $q_n$.

This completes the proof of proposition \ref{prop_liftFnGn}.
\end{proof}

We will need the following improved estimate for $F_n$.

\begin{proposition}\label{prop_Fnfin}
Assume $r<1$. There exists a sequence $(B_n)$, sub-exponential with
respect to $q_n$, such that for all $Z\in \H_n(r)$,
\[\bigl|F_n(Z)-Z-1\bigr|\leq
B_n\cdot \left(|\eps_n|+\bigl|\eps_n-\pi_n(Z)^{q_n}\bigr|\right).\]
\end{proposition}

\begin{proof} Lemma \ref{lem_approxqnfin} below gives a similar
estimate for $f_n^{\circ q_n}$ on $X_n(r)$. This estimate transfers
to the required one by Koebe distortion theorem as in the previous
proof.
\end{proof}

\begin{lemma}\label{lem_approxqnfin}
There exist a complex number $\eta_n$ and a holomorphic function
$k_n$, defined on the same set as $f_n^{\circ q_n}$, such that
\[f_n^{\circ q_n}(z)  = z + \xi_n(z)\cdot
\bigl(1+\eta_n+(\eps_n-z^{q_n})k_n(z)\bigr).\] For all $r<1$, there
exists a sequence $(B_n)$, sub-exponential with respect to $q_n$,
such that
\[|\eta_n| \leq B_n\cdot |\eps_n| \quad\text{and}\quad
 \forall z\in D(0,r)\quad
\bigl|k_n(z)\bigr|\leq B_n.\]
\end{lemma}

\begin{proof} By lemma \ref{lem_approxqn}, we know that
\[f_n^{\circ q_n}(z)= z + \xi_n(z)\cdot h_n(z)\]
with, $\ds B_n\eqdef \sup_{D(0,r)} |h_n|$ a sub-exponential sequence
with respect to $q_n$. The map $f_n^{\circ q_n}$ has the same
multiplier at each $q_n$-th roots of $\eps_n$. It follows that
\[h_n(z) = 1+\eta_n+(\epsilon_n-z^{q_n})k_n(z)\]
as prescribed. Since $\sqrt[q_n]{\eps_n} \to 1/A <r<1$, the bound on
$h_n$, taken at any of the $q_n$-th roots of $\eps_n$, shows that
for $n$ large enough,
\[|1+\eta_n| \leq B_n\]
and thus
\[\forall z\in D(0,r)\quad
\bigl|(\eps_n-z^{q_n})k_n(z)\bigr| \leq 2B_n.\] As in
lemma~\ref{lem_approxqn}, we have for some sequence $r_n \to 1$ and
for $n$ large enough:
\[\sup_{|z|\leq r} |k_n(z)| \leq
B'_n := \frac{2 B_n}{r_n^{q_n}-\epsilon_n}
\]
and $(B'_n)$ is sub-exponential with respect to $q_n$. Looking at
$z=0$ gives:
\[ \frac{e^{2i\pi q_n \eps_n} -1}{2i\pi q_n \eps_n} = 1+ \eta_n +\eps_n k_n(0).\]
As $n\to +\infty$, the left hand of this equality expands to $1+i\pi
q_n \eps_n + o(q_n \eps_n)$. Therefore
\[|\eta_n| \leq \eps_n\bigl(|k_n(0)|+\pi q_n + o(q_n)\bigr).\]
Since $|k_n(0)| \leq B'_n$, we get the desired bound on $\eta_n$.
\end{proof}

\begin{corollary}\label{coro_Fnclosetranslation}
Assume $r<1$. Then,
\[\sup_{Z\in \H_n(r)} \bigl|F_n(Z)-Z-1\bigr|\underset{n \to
  +\infty}\longrightarrow 0
\quad\text{and}\quad \sup_{Z\in \H_n(r)}
\bigl|F_n'(Z)-1\bigr|\underset{n \to
  +\infty}\longrightarrow 0.\]
\end{corollary}

\begin{proof} The first is an immediate consequence of proposition
\ref{prop_Fnfin}. For the second, use the first on $\H_n(r')$ with
$r<r'<1$.
\end{proof}

\subsubsection{Iterating the commuting pair $(F_n,G_n)$}

\begin{proposition}\label{prop_SiegelFnGn}
Assume $1/A<r_1<r_2<1$. If $n$ is sufficiently large, the following
holds. Given any point $Z\in \H_n(r_1)$, there exists a sequence of
integers $(j_\ell)_{\ell\geq 0}$ such that for any integer $\ell\geq
0$ and any integer $j\in [0,j_\ell]$, the point
\[F_n^{\circ j}\circ G_n\circ F_n^{\circ j_{\ell-1}}\circ G_n\circ
\cdots \circ F_n^{\circ j_1}\circ G_n\circ F_n^{\circ j_0} (Z)\] is
well defined and belongs to $\H_n(r_2)$.
\end{proposition}

\begin{proof} We will need to control iterates of $F_n$ for a large
number of iterates. We will use the following lemma.

\begin{lemma}\label{lem_courbeminore}
Assume $F: \H \to \C$ satisfies
\[\bigl|F(Z)-Z-1\bigr|< u\bigl(\re(Z)\bigr)\]
with $u:\R\to \left]0,1/10\right[$ a function such that $\log u$ is
$1/2$-Lipschitz. Let $\Gamma$ be the graph of an antiderivative of
$-2u$. Then, every $Z\in \H$ which is above $\Gamma$ has an image
above $\Gamma$.
\end{lemma}

\begin{proof} Let $U$ be the antiderivative whose graph is $\Gamma$. Let
$Z=X+iY\in\H$. The point $Z'=X'+iY'=F(Z)$ satisfies $X' \in
[X+\frac{9}{10},X+\frac{11}{10}]$. Note that
\[\forall x\in \left[X,X+\frac{11}{10}\right],\quad
\log u(x) \geq \log u(X) -\frac{11}{20}.\] Therefore, from $X$ to
$X'$, $U$ decreases of at least
\[(X'-X)2e^{-11/20}u(X) \geq \frac{18}{10}
e^{-11/20}u(X) > u(X) > Y-Y'.\]

\begin{figure}[htbp]
\begin{picture}(225,80)(0,0)%
\put(0,0){\scalebox{0.5}{\includegraphics{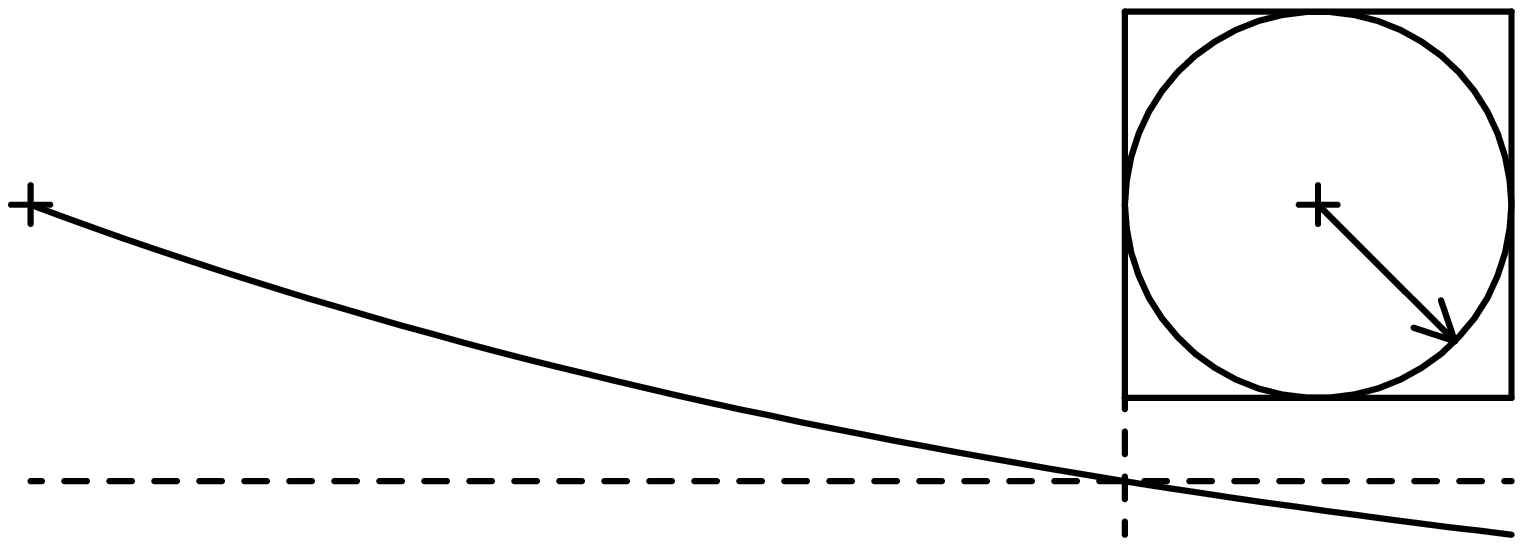}}}%
\put(100,35){$\Gamma$}%
\put(5,58){$Z=X+iY$}%
\put(183,58){\scalebox{0.8}{$Z+1$}}%
\put(183,30){$u(X)$}
\end{picture}
\end{figure}
\end{proof}

\begin{lemma}\label{lemma_return}
Assume $1/A<r<r'<1$. If $n$ is sufficiently large, then for all
$Z\in \H_n(r)$ there exists an integer $j(Z)$ such that
\begin{itemize}
\item for all $j\leq j(Z)$, we have $F_n^{\circ j} \circ G_n(Z)\in \H_n(r')$
and

\item $\re\bigl(F_n^{\circ j(Z)}\circ G_n(Z)\bigr)>\re(Z)$.
\end{itemize}
\end{lemma}

\begin{proof} Let us first recall that there exists a
sequence $(B_n)$, sub-exponential with respect to $q_n$, such that
for $n$ large enough, for all $Z\in \H_n(r)$,
\[\bigl|G_n(Z)-Z+A_n+\theta\bigr|\leq B_n.\] In
particular, if $n$ is sufficiently large,
\[\re\bigl(G_n(Z)\bigr)\geq \re(Z)-A_n-\theta-B_n
\quad\text{and}\quad \im\bigl(G_n(Z)\bigr)\geq \tau_n(r)-B_n.\] We
will apply lemma~\ref{lem_courbeminore} to control the orbit of
$G_n(Z)$ under iteration of $F_n$. More precisely, we will prove the
existence of a function $u_n$ such that:
\begin{itemize}
\item[a)] $\bigl|F_n(Z)-Z-1\bigr| \leq u_n\bigl(\re(Z)\bigr)$,

\item[b)] for $n$ large enough $u_n\in \left]0,1/10\right[$,

\item[c)] for $n$ large enough, $\log u_n$ is
$1/2$-Lipschitz and

\item [d)] the sequence $\ds C_n\eqdef \int_{\re\bigl(G_n(Z)\bigr)}^{\re(Z)} 2u_n(X)\d
X$ is sub-exponential with respect to $q_n$.
\end{itemize}

If $n$ is taken sufficiently large so as to have
\[\tau_n(r)\geq \tau_n(r')+B_n+C_n+\frac{1}{10},\]
it then follows from lemma \ref{lem_courbeminore} that there is an
integer $j(Z)$ such that \begin{itemize} \item for all $j\leq j(Z)$,
we have $F_n^{\circ j} \circ G_n(Z)\in \H_n(r')$ and

\item
$\re\bigl(F_n^{\circ j(Z)}\circ G_n(Z)\bigr)>\re(Z)$.
\end{itemize}

\begin{figure}[htb]
\begin{picture}(215,145)(0,0)%
\put(0,0){\scalebox{.5}{\includegraphics{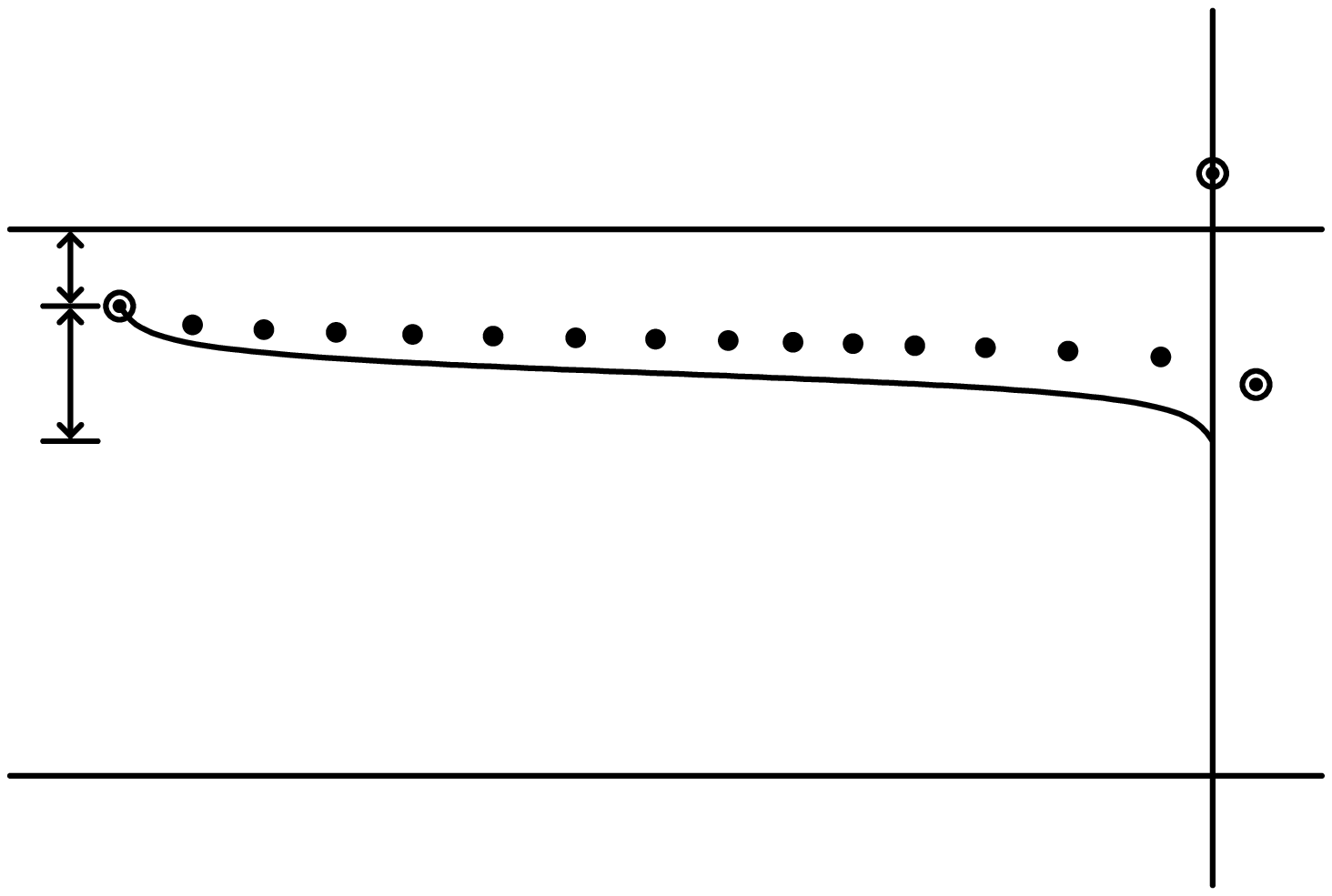} }}%
\put(-7,100){$\scriptstyle \leq B_n$}%
\put(0,82){$\scriptstyle C_n$}%
\put(20,100){$\scriptstyle G_n(Z)$}%
\put(197,120){$\scriptstyle Z$}%
\put(5,25){$\scriptstyle \H_n(r')$}%
\put(5,112){$\scriptstyle \H_n(r)$}%
\put(205,85){$\scriptstyle F_n^{\circ j(Z)}\circ G_n(Z)$}
\end{picture}
\end{figure}

\medskip
\noindent{a)} By proposition \ref{prop_Fnfin}, there is a sequence
$(B'_n)$, sub-exponential with respect to $q_n$, such that for all
$Z\in \H_n(r')$,
\[\bigl|F_n(Z)-Z-1\bigr|\leq B'_n\bigl(\eps_n +
\bigl|\eps_n-\pi_n(Z)^{q_n}\bigr|\bigr).\] Set $T_n\eqdef 1/(2\pi
q_n^2\eps_n) \to +\infty$. We have
\[\bigl(\pi_n(Z)\bigr)^{q_n} =  \frac{\eps_n }{1-e^{-i Z/T_n}}.\]
Using
\[ B'_n\bigl(\eps_n + \bigl|\eps_n-\pi_n(Z)^{q_n}\bigr|\bigr)
   \leq B'_n \bigl(2\eps_n + \bigl|\pi_n(Z)^{q_n}\bigr|\bigr)
\]
we see that for all $Z\in \H_n(r')$,
\begin{eqnarray*}
\bigl|F_n(Z)-Z-1\bigr|&\leq& B'_n\eps_n\left(2+\frac{1}
{\bigl|1-e^{-i
Z/T_n}\bigr|}\right) \\
&\leq& B'_n\eps_n\left(2+\frac{1} {\bigl|s_n e^{i
\re(Z)/T_n}-1\bigr|}\right)
\end{eqnarray*}
with \[ s_n = 1+\frac{\eps_n}{(r')^{q_n}}.\] Since $1/A<r'$, we have
$\eps_n/(r')^{q_n}\to 0$ and thus $s_n\to 1$. Thus, for $n$ large
enough
\[\frac{1}{3}\leq \frac{1}{|s_n e^{i
\re(Z)/T_n}-1|},\] and for all $Z\in \H_n(r')$,
\[\bigl|F_n(Z)-Z-1\bigr| \leq u_n\bigl(\re(Z)\bigr)
\quad\text {with}\quad u_n(X) \eqdef \frac{7 B'_n\eps_n}{|s_n e^{i
    X/T_n}-1|}.\]

\medskip  \noindent{b)}
Let us show that for $n$ large enough $u_n\in \left]0,1/10\right[$.
Note that
\[\forall X\in \R,\quad \bigl|u_n(X)\bigr|\leq \frac{7B'_n
\eps_n}{s_n-1} = 7B'_n(r')^{q_n}\To_{n\to +\infty} 0.\] Thus $u_n$
tends uniformly to $0$ as $n\to +\infty$.

\medskip  \noindent{c)}
Let us now check that for $n$ large enough, $\log u_n$ is
$1/2$-Lipschitz. One computes
\[\frac{u_n'(X)}{u_n(X)} =
   -\frac{s_n}{T_n}\cdot \frac{\sin(X/T_n)}{s_n^2+1-2s_n\cos(X/T_n)}.
\]
This function reaches its extrema when $(s_n^2+1)\cos(X/T_n) -
2s_n=0$. It follows that:
\[\ds \left|\frac{u_n'(X)}{u_n(X)}\right|\leq
\frac{s_n}{T_n(s_n^2-1)} \underset{n\to +\infty}\sim \pi
q_n^2(r')^{q_n}.\] Thus, $\ds \frac{\partial \log u_n(X)}{\partial
  X}$ converges uniformly to $0$ as $n\to +\infty$, and for $n$ large
enough, $\log u_n$ is $1/2$-Lipschitz.

\medskip  \noindent{d)}
Let us finally show that the sequence
\[C_n\eqdef \int_{\re\bigl(G_n(Z)\bigr)}^{\re(Z)} 2u_n(X) \d X\]
is sub-exponential with respect to $q_n$. If $n$ is large enough,
\[\re\bigl(G_n(Z)\bigr)\geq \re(Z) - A_n-\theta_n - B_n\geq -4\pi
T_n.\] Thus,
\[C_n\leq B''_n\eqdef \int_{\re(Z)-4\pi T_n}^{\re(Z)} 2u_n(X)\d X =
4\int_{-\pi T_n}^{\pi T_n} u_n(X)\d X.\] The change of variable
$\theta = X/T_n$, which yields
\[B''_n = \frac{14B'_n}{\pi q_n^2}\int_{-\pi}^{\pi}
\frac{\d\theta}{\sqrt{s_n^2+1 -2 s_n\cos \theta}}.\] It follows that
\[B''_n\underset{n\to +\infty}\sim \frac{28B'_n}{\pi q_n^2}
\log \frac{1}{s_n-1} \underset{n\to +\infty}\sim \frac{28 B'_n}{\pi
q_n}\log (r' A_n).\] By assumption (condition \eqref{eq_logAnsubexp}
in the statement of proposition \ref{theo_disquesdigites}), the
sequence $\log A_n$ is sub-exponential with respect to $q_n$. As a
consequence, $(B''_n)$, and thus $(C_n)$, is sub-exponential with
respect to $q_n$.
\end{proof}

\medskip
\noindent{\em Proof of proposition \ref{prop_SiegelFnGn},
continued}. Remember that we are given $r_1$ and $r_2$ with
$1/A<r_1<r_2<1$ and we want to prove that for $n$ sufficiently
large, any point of $\H_n(r_1)$ has an infinite orbit remaining in
$\H_n(r_2)$ along a well chosen composition of $F_n$ and $G_n$. It
is enough to show that this is true for any sequence of points
\[Z_n=X_n+iY_n\in \H_n(r_1).\]
We will use Douady-Ghys-Yoccoz's renormalization techniques and
follow the presentation in \cite{abc} section 3.2.

\medskip
\noindent{\bf Step 1. Construction of a Riemann surface: ${\cal
V}_n$.} Choose $n$ sufficiently large so that $F_n$ is defined in
the upper half-plane $\bigl\{Z\in \C~;~\im(Z)\geq
\tau_n(r_2)-1/10\}$ with
\[\bigl|F_n(Z)-Z-1\bigr|\leq \frac{1}{10} \quad\text{and}\quad
\bigl|F'_n(Z)-1\bigr|\leq \frac{1}{10}. \footnote{This is possible
by corollary \ref{coro_Fnclosetranslation} applied with $r>r_2$.
Indeed, for $n$ large enough, $\tau_n(r_2)>\tau_n(r)+1/10$.}\] Set
\[P_n\eqdef X_n+i \left(\tau_n(r_2)-\frac{1}{10}\right).\]
Let
\[L_n\eqdef \bigl\{X_n+it~;~ t> \im(P_n)\bigr\}\]
be the vertical half-line starting at $P_n$ and passing through
$Z_n$. The union
\[L_n \cup \bigl[P_n,F_n(P_n)\bigr] \cup  F_n(L_n)\cup
\{\infty\}\] forms a Jordan curve in the Riemann sphere bounding a
region $U_n$ such that for $Y>\im(P_n)$, the segment
$\bigl[iY,F_n(iY)\bigr]$ is contained in $\overline U_n$ (see
\cite{abc}). We set ${\cal U}_n\eqdef U_n\cup L_n$. If we glue the
sides $L_n$ and $F_n(L_n)$ of $\overline {\cal U}_n$ via $F_n$, we
obtain a topological surface $\overline {\cal V}_n$. We denote by
$\iota_n:\overline {\cal U}_n\to \overline{\cal V}_n$ the canonical
projection. The space $\overline{\cal V}_n$ is a topological surface
with boundary, whose boundary $\iota_n\bigl([P_n,F_n(P_n)]\bigr)$ is
denoted $\partial {\cal V}_n$. We set ${\cal V}_n=\overline{\cal
V}_n\setminus \partial {\cal V}_n$. Since the gluing map $F_n$ is
analytic, the surface ${\cal V}_n$ has a canonical analytic
structure induced by the one of ${\cal U}_n$. It is possible to show
that ${\cal V}_n$ is isomorphic to $\H/\Z\simeq \D^*$ (see \cite{y}
for details). Let $\phi_n:{\cal V}_n\to \D^*$ be an isomorphism.
Hence, we have the following composition:
\[\phi_n\circ \iota_n:{\cal U}_n\to \D^*.\]
We set
\[\zeta_n\eqdef  \phi_n\circ \iota_n(Z_n)\in \D.\]

\medskip
\noindent{\bf Step 2. The renormalized map $g_n$.} Choose $r_3\in
\left]r_1,r_2\right[$. Set
\[P'_n\eqdef X_n + i \left(\tau_n(r_3)+\frac{1}{10}\right).\]
Let ${\cal U}'_n$ be the set of points of ${\cal U}_n$ which are
above the segment $\bigl[P'_n,F_n(P'_n)\bigr]$ and let ${\cal V}'_n$
be the image of ${\cal U}'_n$ in ${\cal V}_n$. Choose $n$
sufficiently large so that lemma \ref{lemma_return} can be applied
with $r=r_3$ and $r'=r_2$. Then, for all $Z\in {\cal U}'_n\subset
\H_n(r_3)$, there exists an integer $j(Z)$ such that
\[W\eqdef
F_n^{\circ j(Z)}\circ G_n(Z)\in {\cal U}_n \quad\text{and}\quad
\forall j\in \bigl[0,j(Z)\bigr]\quad F_n^{\circ j}\circ G_n(Z)\in
\H_n(r_2).\]
 The map
$Z\mapsto W$ induces a univalent map $g_n:\phi_n({\cal V}'_n)\to
\D^*$. \footnote{The fact that $g_n:\phi_n({\cal V}'_n)\to \D^*$ is
continuous
  and univalent
  is not completely obvious; see \cite{y} for details.}
By the removable singularity theorem, this map extends
holomorphically to the origin by $g_n(0)=0$. Since
\[F_n(Z) = Z+ 1 + o(1) \quad\text{and}\quad G_n(Z) = Z-A_n-\theta +
  o(1)\]
as $\im(Z)\to +\infty$, it is possible to show that
\[g_n'(0)=e^{-2i\pi (A_n+\theta)}= e^{-2i\pi\theta}\]
(again, see \cite{y} for details).

\begin{figure}[htb]
\begin{picture}(252,189)(0,0)%
\put(0,0){\scalebox{.7}{\includegraphics{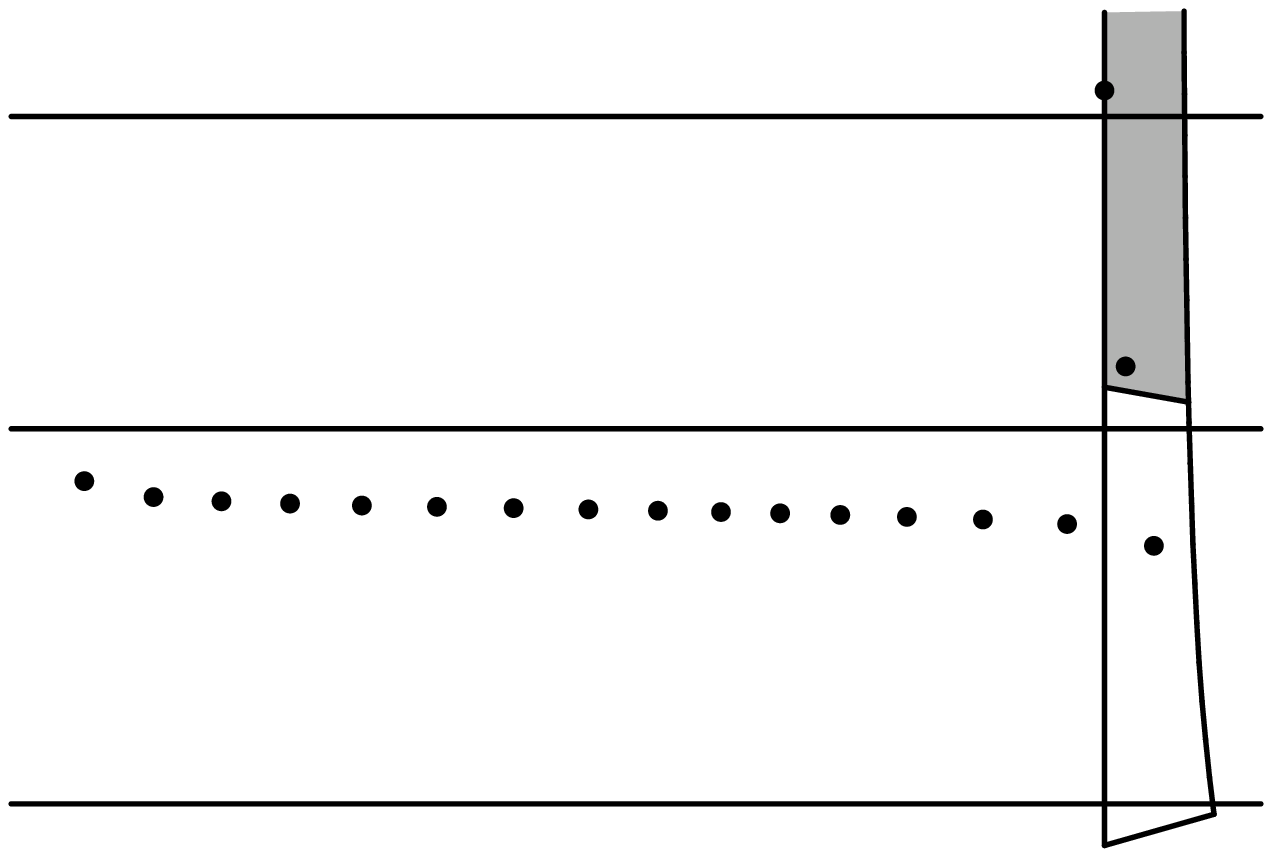} }}%
\put(227,101){$Z$}%
\put(-10,63){$G_n(Z)$}%
\put(229,50){$W$}%
\put(244,30){${\cal U}_n$}%
\put(244,130){${\cal U}'_n$}%
\put(0,14){$\H_n(r_2)$}%
\put(0,90){$\H_n(r_3)$}%
\put(0,153){$\H_n(r_1)$}%
\put(213,173){$L_n$}%
\put(233,173){$F_n(L_n)$}%
\put(213,-10){$P_n$}%
\put(208,153){$Z_n$}%
\put(245,0){$F_n(P_n)$}%
\put(208,95){$P'_n$}%
\put(240,91){$F_n(P'_n)$}%
\end{picture}
\end{figure}

\noindent{\bf Step 3. The orbit of $\zeta_n$.}

We will show that the orbit of $\zeta_n$ under iteration of $g_n$ is
infinite. For this, let $\rho_n$ be the radius of the largest disk
centered at $0$ and contained in $\phi_n({\cal V}'_n)$. We will show
that
\begin{itemize}
\item[a)] $\exists C>0$ such that $g_n$ has a Siegel disk which
  contains $D(0,C\rho_n)$

\item[b)] $|\zeta_n|=o(\rho_n)$.
\end{itemize}

\medskip
\noindent a) The restriction of $g_n$ to $D(0,\rho_n)$ is univalent.
It fixes $0$ with derivative $e^{-2i\pi
  \theta}$. Remember that $\theta$ is a Brjuno number.
It follows (see \cite{bru} or \cite{y} for example) that there is a
constant $C_\theta>0$ depending only on $\theta$ such that $g_n$ has
a Siegel disk containing $D(0,C_\theta\rho_n)$.

\medskip
\noindent b) Denote by $B_n$ the half-strip
\[B_n=\{Z\in \C~;~0<\re(Z)<1\text{ and }\im(Z)>\im(P_n)\}\]
and consider the map $H_n:\overline{B}_n\to \overline{{\cal U}}_n$
defined by
\[ H_n(Z) = (1-X)\cdot (X_n+iY)+X\cdot F_n(X_n+iY) \]
where $Z=X+i Y$, $(X,Y) \in [0,1]\times
\bigl[\im(P_n),+\infty\bigr[$. The map $H_n$ sends each segment
$[iY,iY+1]$ to the segment $\bigl[X_n+iY,F_n(X_n+iY)\bigr]$. An
elementary
    computation shows that $H_n$ is a
$5/4$-quasiconformal homeomorphism between $\overline{B}_n$ and
$\overline{\cal U}_n$.\footnote{For a proof that $H_n$ is
$5/4$-quasiconformal homeomorphism, see for example \cite{abc}
section 3.2 or \cite{sh} section 2.5.} Since $H_n(iY+1) =
F_n\bigl(H_n(iY)\bigr)$, the quasiconformal homeomorphism
$H_n:\overline B_n\to \overline {\cal U}_n$ induces a homeomorphism
between the half cylinder $\H/\Z$ and the Riemann surface ${\cal
V}_n$. This homeomorphism is clearly quasiconformal on the image of
$B_n$ in $\H/\Z$, i.e., outside a straight line. It is therefore
quasiconformal in the whole half cylinder ($\R$-analytic curves are
removable for quasiconformal homeomorphisms).

Let $R_n$ be the rectangle
\[R_n\eqdef \bigl\{Z\in \C~;~0\leq \re(Z)< 1~\text{and}~
\im(P'_n)<\im(Z)<\im(Z_n)\bigr\}.\] Note that $H_n(R_n)\subset {\cal
U}'_n$ and observe that
\[{\cal A}_n\eqdef \phi_n\circ \iota_n \circ H_n(R_n)\]
is an annulus contained in $\phi_n({\cal V}'_n)$ that surrounds $0$
and $\zeta_n$.
\begin{figure}[htb]
\begin{picture}(300,260)(0,0)%
\put(0,0){\scalebox{.5}{\includegraphics{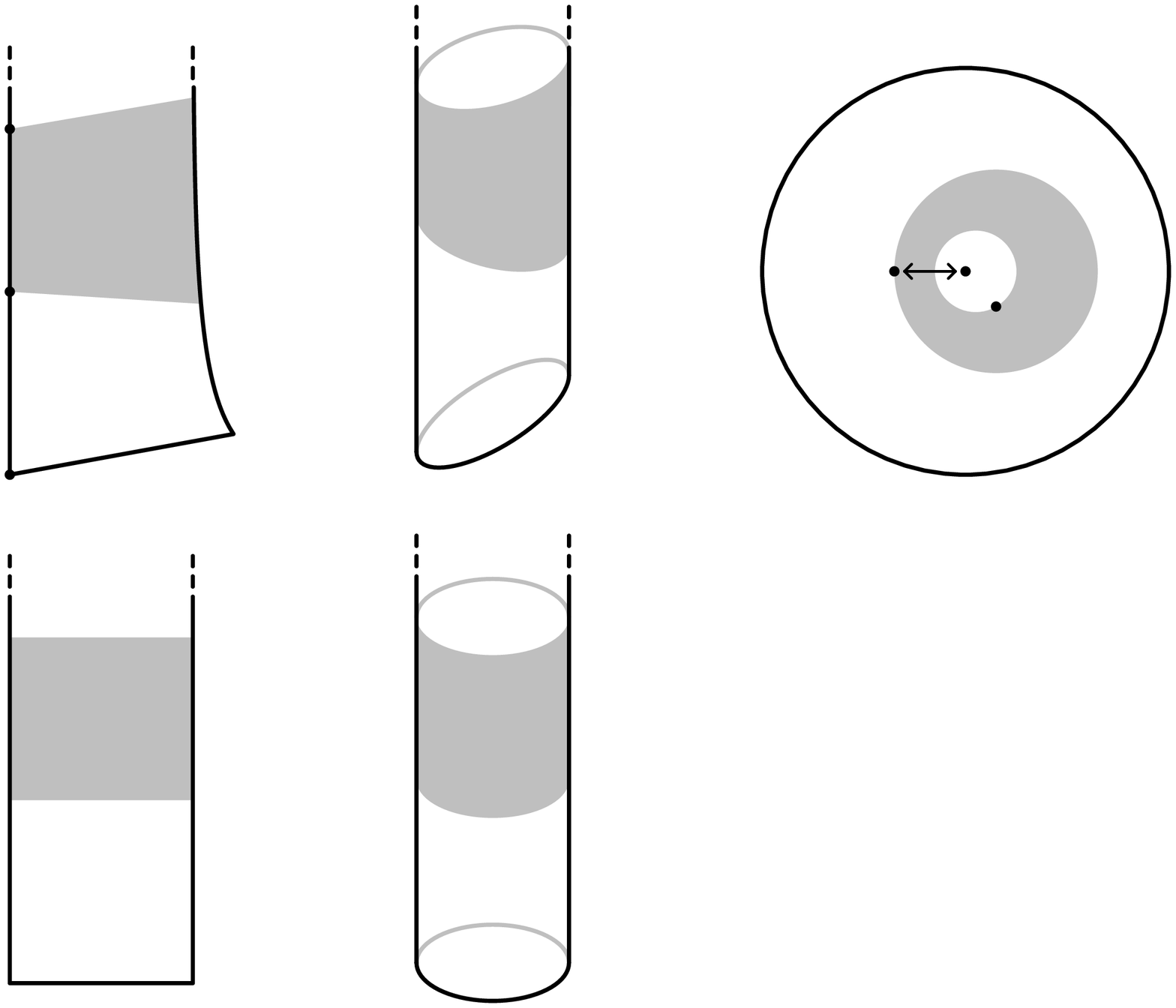} }}%
\put(20,70){$R_n$}%
\put(10,195){$H_n(R_n)$}%
\put(24,119){\rotatebox{90}{$\longrightarrow$}}%
\put(28,124){$H_n$}%
\put(-22,10){$\scriptstyle \im(P_n)$}%
\put(-22,55){$\scriptstyle \im(P'_n)$}%
\put(-22,95){$\scriptstyle \im(Z_n)$}%
\put(-8,135){$\scriptstyle P_n$}%
\put(-8,180){$\scriptstyle P'_n$}%
\put(-8,220){$\scriptstyle Z_n$}%
\put(147,10){$\H/\Z$}%
\put(147,155){${\cal V}_n$}%
\put(65,65){$\overset{Z\sim Z+1}\longrightarrow$}%
\put(65,200){$\overset{Z\sim F_n(Z)}{\underset{\iota_n}\longrightarrow}$}%
\put(160,190){$\underset{\phi_n}\longrightarrow$}%
\put(227.5,190){$\scriptstyle \rho_n$}%
\put(252,176){$\scriptstyle \zeta_n$}%
\put(242,187){$\scriptstyle 0$}%
\put(253,195){$\scriptstyle {\cal A}_n$}%
\put(65,150){${\cal U}_n$}%
\put(55,10){$B_n$}%
\end{picture}
\end{figure}

The image of $R_n$ in $\H/\Z$ is an annulus of modulus
\[M_n\eqdef \im(Z_n) -\im(P'_n) \geq \tau_n(r_1)-
\tau_n(r_3)-\frac{1}{10} \underset{n\to +\infty}\longrightarrow
+\infty.\] Note that $H_n$ induces a $5/4$-quasiconformal
homeomorphism between this annulus and ${\cal A}_n$. It follows that
\[{\rm modulus}({\cal A}_n)\geq \frac{4}{5} M_n \underset{n\to +\infty}\longrightarrow
+\infty.\] Since $A_n$ separates $0$ and $\zeta_n$ from $\infty$ and
a point of modulus $\rho_n$ in $\partial \phi_n ({\cal V}'_n)$, the
claim follows: as $n\to +\infty$, $|\zeta_n|=o(\rho_n)$.

\medskip
\noindent{\bf Step 4. Controlling the orbit of $Z_n$.}

We know that the orbit of $\zeta_n$ under iteration of $g_n$ is
infinite. Thus, we have a sequence
\[\zeta_n\in {\cal V}'_n \overset{g_n}\longrightarrow \zeta_n^1 \in {\cal
V}'_n \overset{g_n}\longrightarrow \zeta_n^2 \in {\cal V}'_n
\overset{g_n}\longrightarrow \cdots.\] Now, for each $\ell\geq 0$,
we have
\[\zeta_n^{\ell} = \phi_n\circ \iota_n(Z_n^{\ell})\quad \text{for
  some}
\quad  Z_n^\ell \in {\cal U}'_n.\] Moreover, by definition of $g_n$,
there exists an integer $j_\ell$ such that
\[Z_n^{\ell+1} = F_n^{\circ j_\ell}\circ
G_n(Z_n^\ell)\quad\text{and}\quad \forall j\in [0,j_\ell] \quad
F_n^{\circ j}\circ G_n(Z_n^\ell)\in \H_n(r_2).\] In other words,
$\ds \zeta_n^{\ell}\in {\cal V}'_n \overset{g_n}\longrightarrow
\zeta_n^{\ell+1} \in {\cal V}'_n$ corresponds to
\[Z_n^\ell \in {\cal U}'_n \overset{G_n}\longrightarrow \cdot \in
\H_n(r_2) \overset{F_n}\longrightarrow \cdot \in \H_n(r_2)
\overset{F_n}\longrightarrow \cdots \overset{F_n}\longrightarrow
Z_n^{\ell+1}\in {\cal U}'_n.\] Thus, for $n$ sufficiently large, any
point $Z_n\in \H_n(r_1)$ has an infinite orbit remaining in
$\H_n(r_2)$ along a well chosen composition of $F_n$ and $G_n$. This
completes the proof of proposition \ref{prop_SiegelFnGn}.
\end{proof}

\noindent{\em Proof of proposition \ref{prop_SiegelXn}', continued.}
Remember that $0<1/A<\rho<\rho'<1$. Choose $r_1=\rho<r_2<\rho'$. By
proposition \ref{prop_SiegelFnGn}, for $n$ sufficiently large, any
point $Z\in\H_n(\rho)$ has an infinite orbit remaining in
$\H_n(r_2)$ under a well chosen composition of $F_n$ and $G_n$. This
means that any point $z\in X_n(\rho)$ has an infinite orbit
remaining in $X_n(r_2)$ under a well chosen composition of
$f_n^{\circ q_n}$ and $f_n^{\circ q_{n-1}}$. By corollary
\ref{coro_jellouli}, if $n$ is sufficiently large, we know that any
point in $X_n(r_2)\subset D(0,r_2)$ has its first $q_n$ iterates in
$D(0,\rho')$. This shows that any point $z\in X_n(\rho)$ has an
infinite orbit remaining in $D(0,\rho')$ under iteration of $f_n$,
as required.

In other words,
\[\cdot
\in \H_n(r_2) \overset{G_n}\longrightarrow \cdot \in
\H_n(r_2)\quad\text{corresponds to}\quad \cdot \in X_n(r_2)
\overset{f_n^{\circ q_{n-1}}}\longrightarrow \cdot \in X_n(r_2)\]
and
\[\cdot
\in \H_n(r_2) \overset{F_n}\longrightarrow \cdot \in
\H_n(r_2)\quad\text{corresponds to}\quad \cdot \in X_n(r_2)
\overset{f_n^{\circ q_{n}}}\longrightarrow \cdot \in X_n(r_2).\]
Moreover, for $n$ sufficiently large,
\[\cdot \in X_n(r_2) \overset{f_n^{\circ q_{n-1}}}\longrightarrow \cdot
\in X_n(r_2)\quad\text{and}\quad \cdot \in X_n(r_2)
\overset{f_n^{\circ q_{n}}}\longrightarrow \cdot \in X_n(r_2)\]
decompose as
\[\cdot \in X_n(r_2)\subset D(0,r_2) \overset{f_n}\longrightarrow \cdot
\in D(0,\rho')\overset{f_n}\longrightarrow \cdots
\overset{f_n}\longrightarrow \cdot\in
D(0,\rho')\overset{f_n}\longrightarrow \cdot X_n(r_2).\] This
completes the proof of proposition \ref{prop_SiegelXn}'.
\qed\par\medskip

\subsection{The control of the post-critical
set\label{sec_inoushishikura}}

\begin{definition}
We denote by $\partial$ the Hausdorff semi-distance:
\[\partial(X,Y) = \sup_{x\in X} d(x,Y).\]
\end{definition}

\begin{definition}
We denote by $\PC(P_\a)$ the post-critical set of $P_\a$:
\[\PC(P_\a) \eqdef \bigcup_{k\geq 1} P_\a^{\circ
k}(\omega_\a)\quad\text{with}\quad \omega_\a\eqdef
-\frac{e^{2i\pi\a}}{2}.\]
\end{definition}

This section is devoted to the proof of the following proposition.
Remember that ${\cal S}_N$ is the set of irrational numbers of
bounded type whose continued fractions have entries greater than or
equal to $N$.

\begin{proposition}\label{theo_postcrit}
There exists $N$ such that as $\a'\in {\cal S}_N\to \a\in {\cal
S}_N$, we have
\[\partial\bigl(\PC(P_{\a'}),\overline \Delta_\a\bigr)\to 0,\]
with $\Delta_\a$ the Siegel disk of $P_\a$.
\end{proposition}

The corollary we will use later is the following.

\begin{corollary}\label{cor_postcrit}
Let $(\a_n)$ be the sequence defined in proposition
\ref{prop_keycremer}. If $n$ is large enough, the post-critical set
of $P_{\a_n}$ is contained in the $\delta$-neighborhood of the
Siegel disk of $P_\a$.
\end{corollary}

The proof of proposition \ref{theo_postcrit} will rely on some
(almost) classical results on Fatou coordinates and perturbed Fatou
coordinates. We refer the reader to appendix \ref{appendix_fatou}
and to \cite{sh} for more details. The proof will also rely on
results of Inou and Shishikura \cite{is} that we will now recall.

\subsubsection{The class of Inou and Shishikura}

Consider the cubic polynomial \[P(z)= z(1+z)^2.\] This polynomial
has a multiple fixed point at $0$, a critical point at $ -1/3$ which
is mapped to the critical value at $-4/27$, and a second critical
point at $-1$ which is mapped to $0$. We set \[R\eqdef
e^{4\pi}\quad\text{and}\quad v\eqdef-4/27.\] Let $U$ be the open set
defined by
\[U\eqdef
P^{-1}\bigl(D(0,|v|R)\bigr)\setminus
\bigl(\left]-\infty,-1\right]\cup B\bigr),\]  where $B$ is the
connected component of  $P^{-1}\bigl(D(0,|v|/R)\bigr)$ which
contains $-1$.

\begin{figure}[htbp]%
\centerline{\begin{picture}(0,0)%
\scalebox{.75}{\includegraphics{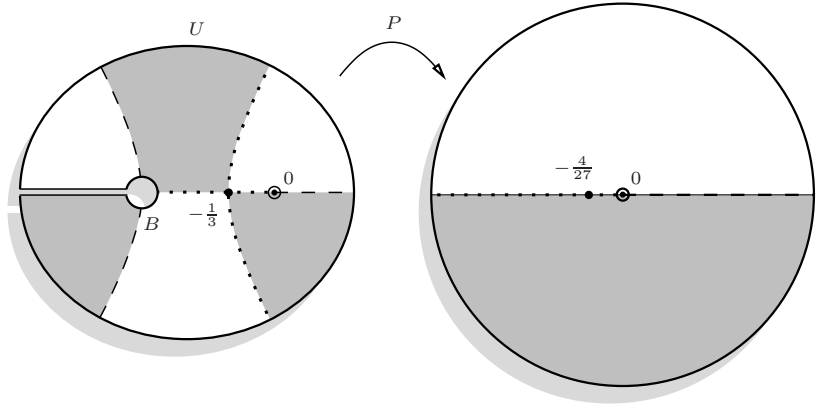}}%
\end{picture}%
\setlength{\unitlength}{1242sp} 
\begin{picture}(15875,7646)(1600,-7934)%
\put(4800,-4596){$\scriptstyle -\frac{1}{3}$}%
\put(3900,-4800){$\scriptstyle B$}%
\put(6701,-3886){$\scriptstyle 0$}%
\put(12101,-3661){$\scriptstyle -\frac{4}{27}$}%
\put(8746,-801){$\scriptstyle P$}%
\put(4816,-901){$\scriptstyle U$}%
\put(13636,-3886){$\scriptstyle 0$}%
\end{picture}%
}%
\caption{A schematic representation of the set $U$. We colored gray
the set of points in $U$ whose image by $P$ is contained in the
lower half-plane.}%
\end{figure}

Consider the following class of maps (Inou and Shishikura use the
notation ${\cal F}'_1$ in \cite{is}):
\[\sF_0\eqdef \left\{f=P\circ \varphi^{-1}:U_f\to \C\mbox{ with
}\begin{array}{l}\varphi:U\to U_f\mbox{ isomorphism such that }\\
\varphi(0)=0\mbox{ and }\varphi'(0)=1\end{array}\right\}.\]

\begin{rema}
The set $\sF_0$ is identified with the space of univalent maps in
$U$ fixing $0$ with derivative $1$, which is compact. A sequence of
univalent maps $(\varphi_n:U\to \C)$ satisfying $\varphi_n(0)=0$ and
$\varphi_n'(0)=1$ converges uniformly to $\varphi:U\to \C$ on every
compact subset of $U$, if and only if the sequence $(f_n=P\circ
\varphi_n^{-1})$ converges to $f=P\circ \varphi^{-1}$ on every
compact subset of $U_f=\varphi_f(U)$.
\end{rema}

A map $f\in \sF_0$ fixes $0$ with multiplier $1$. The map $f:U_f\to
D\bigl(0,|v|R\bigr)$ is surjective. It is not a proper map. Inou and
Shishikura call it a {\em partial covering}. The map $f$ has a
critical point $\omega_f \eqdef  \varphi_f(-1/3)$ which depends on
$f$ and a critical value $v\eqdef -4/27$ which does not depend on
$f$.

\subsubsection{Fatou coordinates\label{subsubfatou}}

Near $z=0$, elements $f\in \sF_0$ have an expansion of the form
\[f(z)=z+c_f z^2+ {\cal O}(z^3).\]
The following result of Inou and Shishikura is an immediate
consequence of the Koebe $1/4$-theorem.

\begin{is}[Main theorem 1 part a]
The set $\{c_f~;~f\in \sF_0\}$ is a compact
subset of $\C^*$.
\end{is}

In particular, for all $f\in \sF_0$, $c_f\neq 0$ and $f$ has a
multiple fixed point of multiplicity $2$ at $0$. If we make the
change of variables
\[w=\tau_f(z)\eqdef -\frac{1}{c_f z},\] we find $F(w)= w+1+o(1)$
near infinity. To lighten notation, we will write $f$ and $F$ for
pairs of functions related as above; $\omega_f\eqdef \phi_f(-1/3)$
and $\omega_F\eqdef \tau_f^{-1}(\omega_f)$ will denote their
critical points.

\begin{lemma}
There exists $R_0$ such that for all $f\in \sF_0$
\begin{itemize}
\item $F$ is defined and univalent in a neighborhood of
$\C\setminus D(0,R_0)$ and

\item for all $w\in \C\setminus D(0,R_0)$,
\[\bigl|F(w)-w-1\bigr|<\frac{1}{4}\quad\text{and}\quad
\bigl|F'(w)-1\bigr|<\frac{1}{4}.\]
\end{itemize}
\end{lemma}

\begin{proof}
This follows from the compactness of $\sF_0$.
\end{proof}

If $R_1>\sqrt 2R_0$, the regions
\[\Omega^\att\eqdef\bigl\{w\in \C~;~\re(w)>R_1-|\im(w)|\bigr\}
\]
and
\[\Omega^\rep\eqdef
\bigl\{w\in \C~;~\re(w)<-R_1+|\im(w)|\bigr\}\] are contained in
$\C\setminus D(0,R_0)$.

\begin{figure}[htbp]
\begin{picture}(300,150)
 \put(0,0){\scalebox{0.75}{\includegraphics{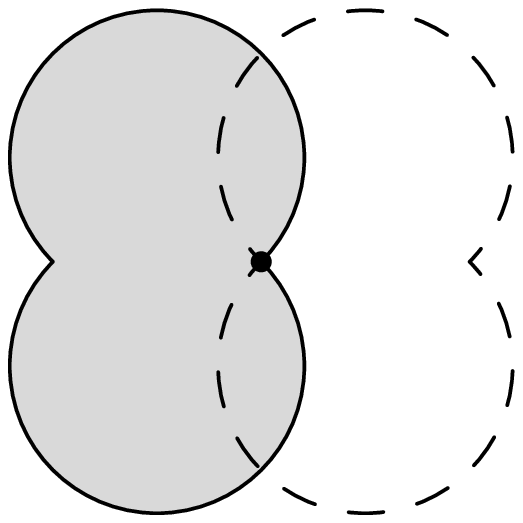}}}
 \put(150,0){\scalebox{0.75}{\includegraphics{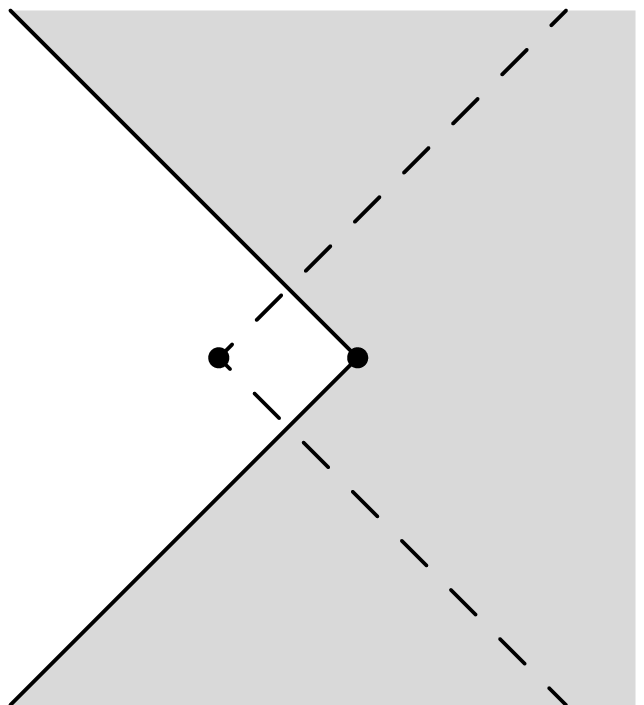}}}
 \put(145,73){$\overset{\tau_f}{\longleftarrow}$}
 \put(30,40){$\Omega_{\att,f}$}
 \put(95,40){$\Omega_{\rep,f}$}
 \put(260,110){$\Omega^\att$}
 \put(170,30){$\Omega^\rep$}
 \put(82,73){$0$}
 \put(245,73){$R_1$}
 \put(185,73){$-R_1$}
\end{picture}
\caption{Right: the sets $\Omega^\att$ and $\Omega^\rep$. Left: the
set $\Omega_{\att,f}$ and $\Omega_{\rep,f}$ for a map $f$ with
$c_f=1$. The sets $\Omega^\att$ and $\Omega_{\att,f}$ are shaded.
The boundaries of the sets $\Omega^\rep$ and $\Omega_{\rep,f}$ are
dashed. }
\end{figure}

Then, for all $f\in \sF_0$,
\[F(\Omega^\att)\subset \Omega^\att
\quad \text{and}\quad F\bigl(\Omega^\rep\bigr)\supset \Omega^\rep.\]
In addition, there are univalent maps $\Phi^\att_F:\Omega^\att\to
\C$ (attracting Fatou coordinate for $F$) and
$\Phi^\rep_F:\Omega^\rep\to \C$ (repelling Fatou coordinate for $F$)
such that
\[\Phi^\att_F\circ F(w) = \Phi^\att_F (w) +
1\quad\text{and}\quad \Phi^\rep_F\circ F(w)= \Phi^\rep_F (w) + 1\]
when both sides of the equations are defined. The maps $\Phi^\att_F$
and $\Phi^\rep_F$ are unique up to an additive constant.

\begin{is}[Main theorem 1 part a]
For all $f\in \sF_0$, the critical point $\omega_f$ is attracted to
$0$.
\end{is}

The following lemma easily follows, using the compactness of the
class $\sF_0$.

\begin{lemma} There exists $k$ such that for all $f\in \sF_0$ we have
$F^{\circ k}(\omega_F) \in \Omega^\att$.
\end{lemma}

\begin{proof} By contradiction, suppose that there is a sequence
$(f_n)\in \sF_0$ such that for $k\le n$ we have $F_n^{\circ
k}(\omega_{F_n})\notin \Omega^\att$. By compactness of $\sF_0$ we
may assume that the sequence $F_n$ converges to $F_\infty$. But
since $f_\infty\in \sF_0$, the orbit of the critical point
$\omega_{f_\infty}$ converges to $0$, so for some $k$ we have
$F_\infty^{\circ k}(\omega_{F_\infty}) \in \Omega^\att$. But
\[
  F_\infty^{\circ k}(\omega_{F_\infty})= \lim_{n\to \infty} F_n^{\circ k}(\omega_{F_n})
\]
and this is a contradiction.
\end{proof}

Since the maps $\Phi^\att_F$ and $\Phi^\rep_F$ are only defined up
to an additive constant, we can normalize $\Phi^\att_F$ so that
\[\Phi^\att_F\bigl(F^{\circ k}(\omega_F)\bigr) = k.\]
Then, we can normalize $\Phi^\rep_F$ so that
\[\Phi^\att_F(w)-\Phi^\rep_F(w)\to 0\quad\text{when}\quad \im(w)\to
+\infty \quad\text{with}\quad w\in \Omega^\att\cap
\Omega^\rep.
\]

Coming back to the $z$-coordinate, we define
\[\Omega_{\att,f}\eqdef \tau_f(\Omega^\att)\quad\text{and}\quad
\Omega_{\rep,f}\eqdef \tau_f(\Omega^\rep)\] and we set
\[\Phi_{\att,f}\eqdef \Phi^\att_F \circ
\tau_f^{-1}\quad\text{and}\quad \Phi_{\rep,f}\eqdef \Phi^\rep_F
\circ \tau_f^{-1}.\] The univalent maps
$\Phi_{\att,f}:\Omega_{\att,f}\to \C$ and
$\Phi_{\rep,f}:\Omega_{\rep,f}\to\C$ are called attracting and
repelling Fatou coordinates for $f$. Note that our normalization of
the attracting coordinates is given by
\[\Phi_{\att,f}\bigl(f^{\circ k}(\omega_f)\bigr) = k.\]

The following result of Inou and Shishikura asserts that the
attracting Fatou coordinate can be extended univalently up to the
critical point of $f$. It  easily follows from \cite{is} Proposition
5.6.

\begin{figure}[htbp]
\begin{picture}(330,150)
\put(0,20){\scalebox{0.6}{\includegraphics{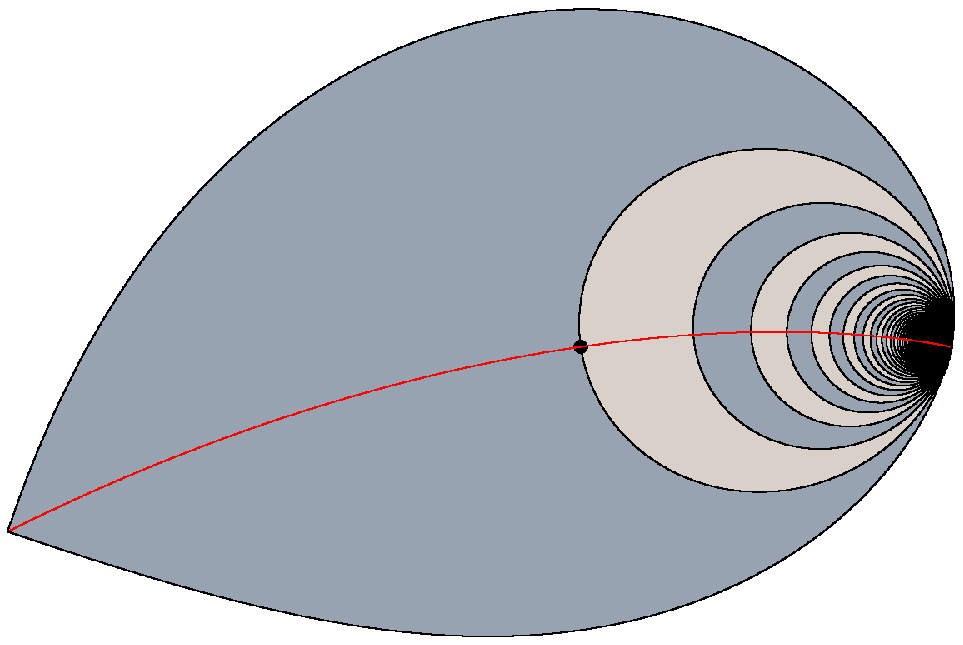}}}
\put(200,0){\scalebox{0.6}{\includegraphics{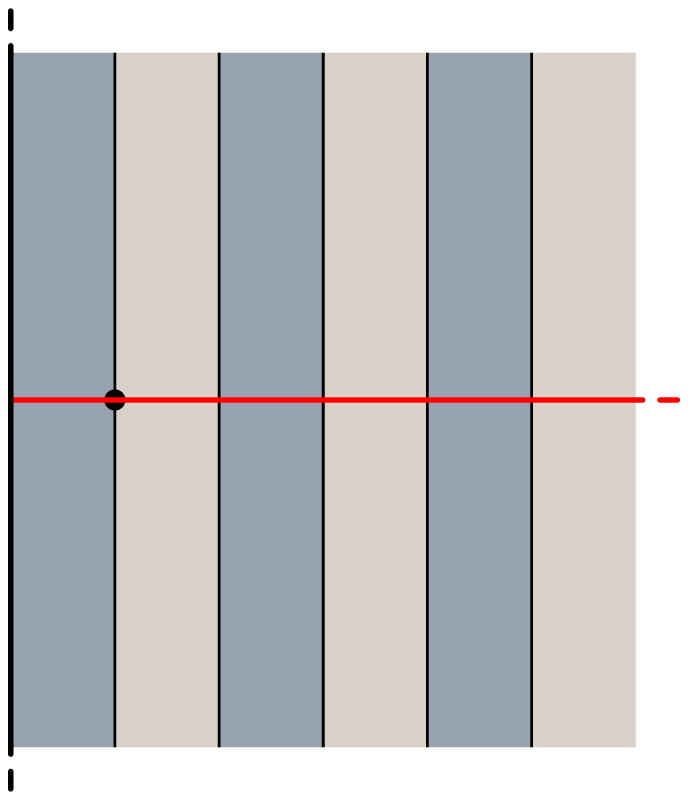}}}
\put(175,73){$\overset{\Phi_{\att,f}}{\longrightarrow}$}
\put(120,15){$\Pet_{\att,f}$} \put(106,68){$v$}
\put(-4,32){$\omega_f$} \put(169,65){$0$}
\end{picture}
\caption{\label{fig_petal}Left: the attracting petal $\Pet_{\att,f}$
of some map $f\in\sF_0$; the critical point is $\omega_f$, the
critical value $v$ and $0$ is a fixed point. Right: its image by
$\Phi_{\att,f}$; we divided the right half plane
$]0,+\infty[\times\R$ into vertical strips of width $1$ of
alternating color, highlighted the real axis in red, and put a black
dot at the point $z=1$. On the left, we pulled this coloring back by
$\Phi_{\att,f}$.}
\end{figure}

\begin{is}[see figure \ref{fig_petal}]
For all $f\in \sF_0$, there exists an attracting petal
$\Pet_{\att,f}$ and an extension of the Fatou coordinate, that we
will still denote $\Phi_{\att,f}:{\cal P}_{\att,f}\to \C$, such that
\begin{itemize}
\item $v\in \Pet_{\att,f}$,

\item $\Phi_{\att,f}(v)=1$ and

\item $\Phi_{\att,f}(\Pet_{\att,f}) = \bigl\{w~;~\re(w)>0\bigr\}$.
\end{itemize}
\end{is}

\begin{figure}[htbp]
\begin{picture}(300,170)%
\put(0,21.4){\scalebox{0.6}{\includegraphics{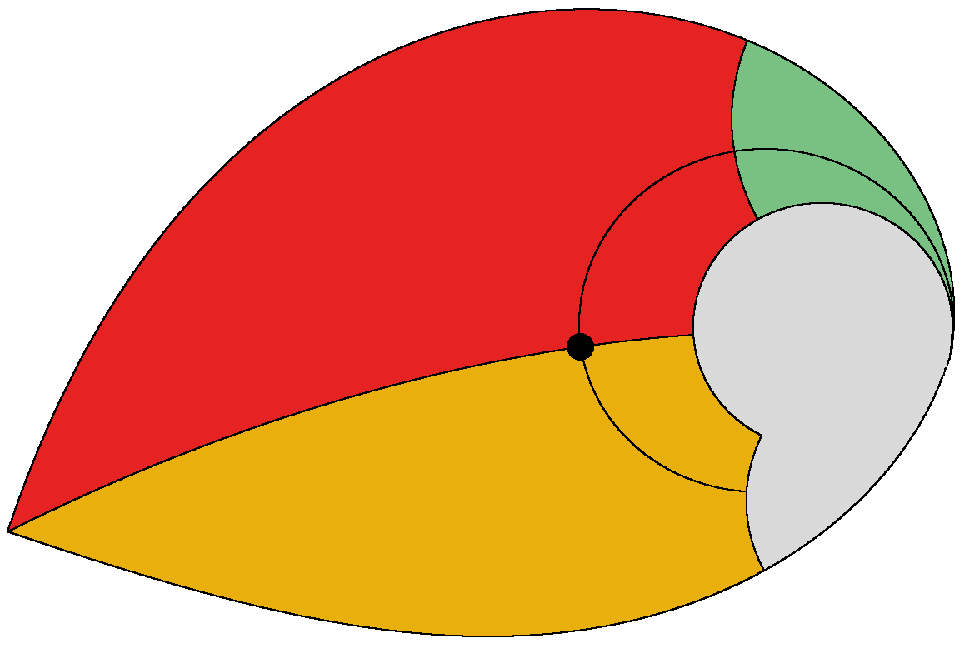}}}%
 \put(189,-4){\scalebox{0.47}{\includegraphics{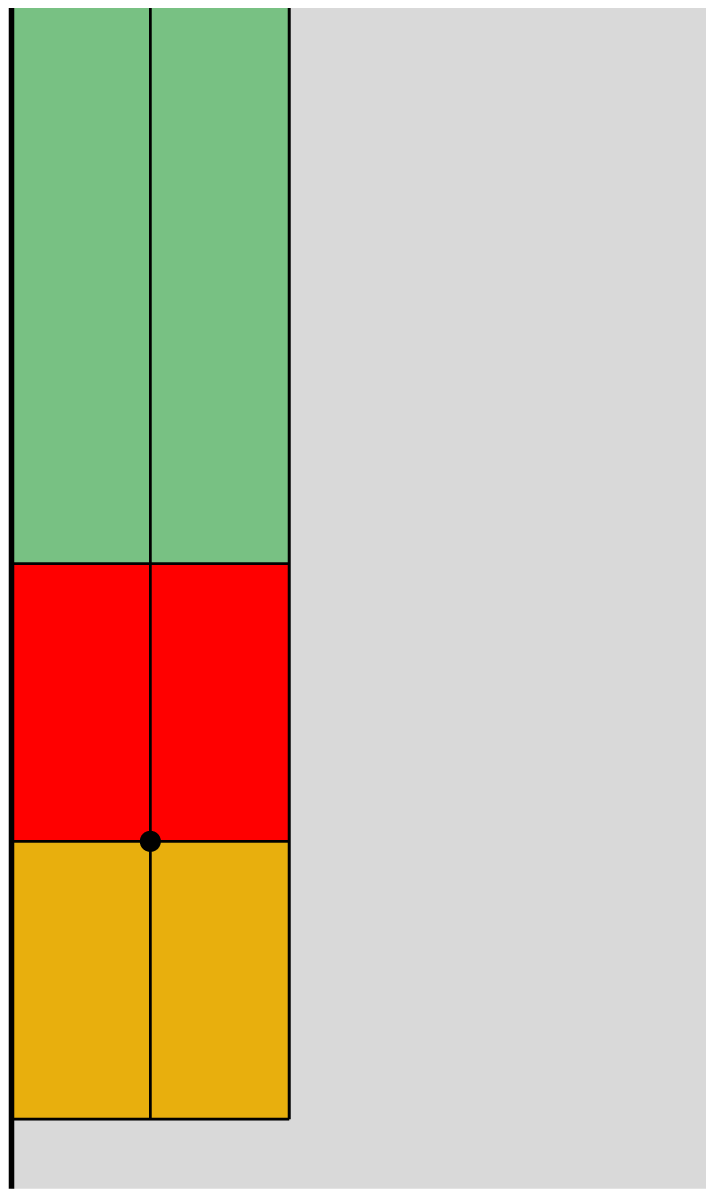}}}%
 \put(171,77){$\overset{\Phi_{\att,f}}{\longrightarrow}$}%
\end{picture}
\caption{\label{fig_defvfwf}On the right, we divided $]0,2[\times]-2,+\infty[$
  into $3$
regions of different colors. We subdivided each by a vertical line
through $z=1$. These $6$ pieces were then pulled back on the left by
$\Phi_{\att,f}$, for the same parabolic $f\in\sF_0$ as in
figure~\ref{fig_petal}.
}
\end{figure}

\begin{definition}[see figure \ref{fig_defvfwf}]
For $f\in \sF_0$, we set:
\[V_f\eqdef \left\{z\in {\cal
P}_{\att,f}~;~\im\bigl(\Phi_{\att,f}(z)\bigr)>0\mbox{ and }
0<\re\bigl(\Phi_{\att,f}(z)\bigr)<2 \right\}\] and
\[W_f\eqdef \left\{z\in \Pet_{\att,f}~;~
-2<\im\bigl(\Phi_{\att,f}(z)\bigr)<2\mbox{ and }
0<\re\bigl(\Phi_{\att,f}(z)\bigr)<2\right\}.\]
\end{definition}

We now come to the key result of Inou and Shishikura. The result
stated below easily follows from \cite{is} Prop. 5.5 and 5.7. Our
domain $V_f^{-k}\cup W_f^{-k}$ below corresponds in \cite{is} to
\[D_{-k}\cup D_{-k}^\sharp\cup D_{-k}''\cup D_{-k+1}\cup
D_{-k+1}^\sharp\cup D_{-k+1}'.\]

\begin{figure}[htbp]
\begin{picture}(344,144)%
 \put(0,0){\scalebox{0.25}{\includegraphics{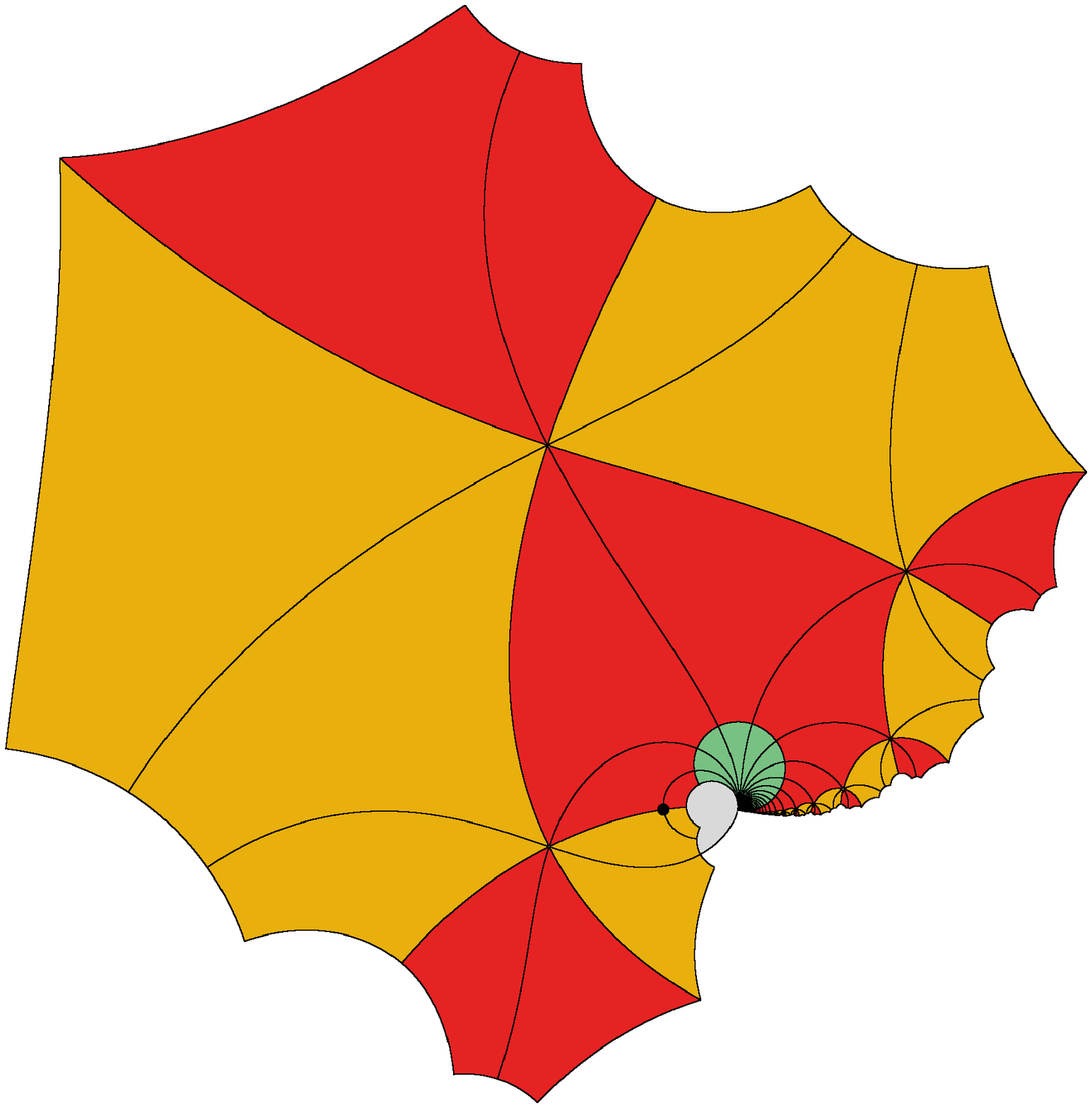}}}%
 \put(152,0){\scalebox{0.5}{\includegraphics{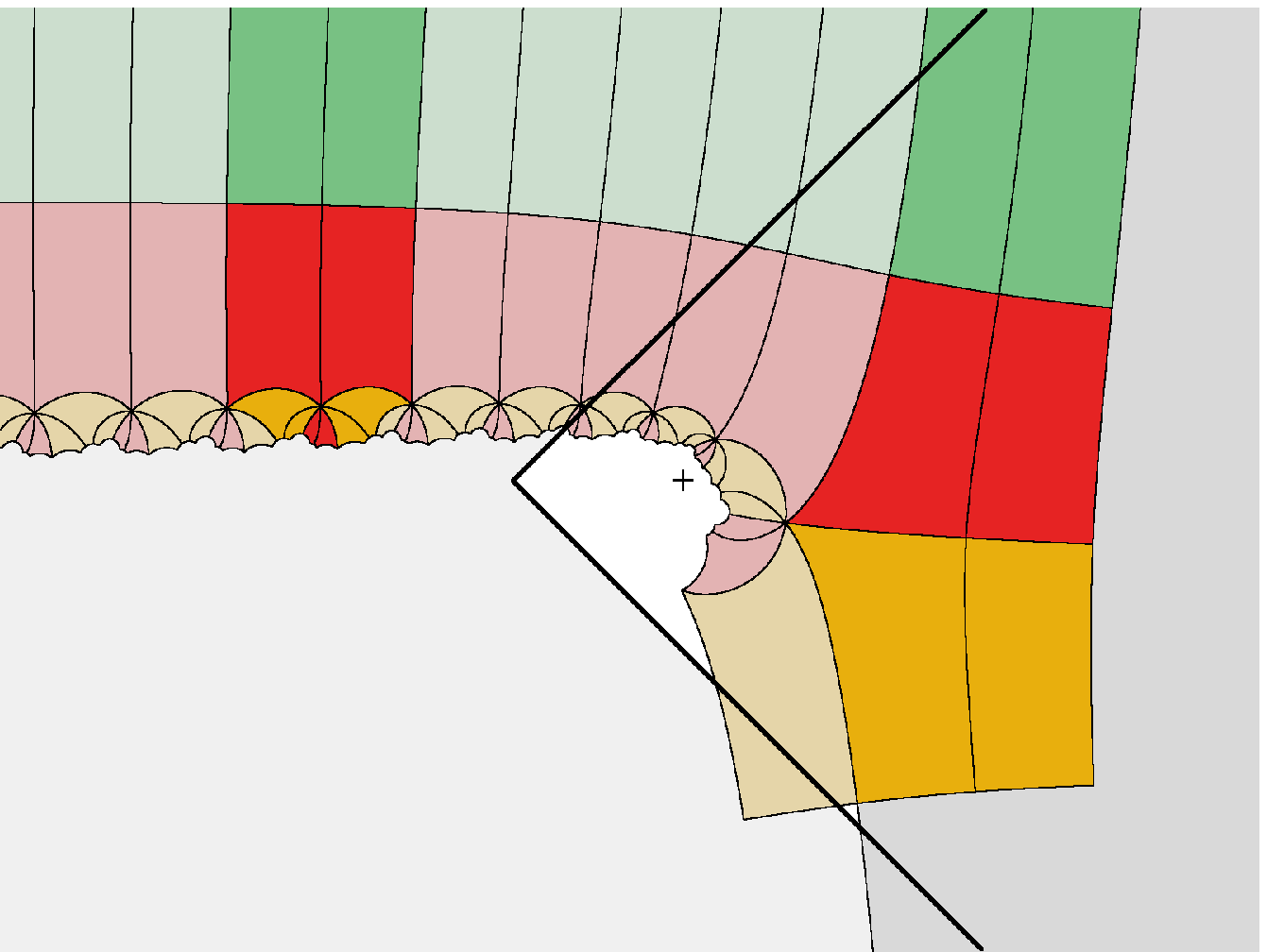}}}%
 \put(240,10){$\Omega^\rep$}
 \put(130,72){$\overset{\tau_0}{\longleftarrow}$}
\end{picture}
\caption{\label{fig_vks}%
Left: among the successive preimages of $V_f$ and $W_f$ by $f$,
those that compose the sets $V_f^{-k}$, $W_f^{-k}$ are shown. The
colors are preserved by $f$. Right: preimage of the left part by
$\tau_0$.  We highlighted $W_F\cup V_F$ and $W_F^{-7}\cup V_F^{-7}$.
}
\end{figure}

\begin{is}[see figure \ref{fig_vks}]
For all $f\in \sF_0$ and all $k\geq 0$,
\begin{itemize}
\item the unique connected component
$V_f^{-k}$ of $f^{-k}(V_f)$ which contains $0$ in its closure is
relatively compact in $U_f$ (the domain of $f$) and $f^{\circ
k}:V_f^{-k}\to V_f$ is an isomorphism and

\item the unique connected component $W_f^{-k}$ of
$f^{-k}(W_f)$ which intersects $V_f^{-k}$ is relatively compact in
$U_f$ and $f^{\circ k}:W_f^{-k}\to W_f$ is a covering of degree $2$
ramified above $v$.
\end{itemize}
In addition, if $k$ is large enough, then $V_f^{-k}\cup
W_f^{-k}\subset \Omega_{\rep,f}$.
\end{is}

The following lemma asserts that if $k$ is large enough, then for
all map $f\in \sF_0$, the set $V_f^{-k}\cup W_f^{-k}$ is contained
in a repelling petal of $f$, i.e. the preimage of a left half-plane
by $\Phi_{\rep,f}$.

\begin{figure}[htbp]
\begin{picture}(350,185)%
 \put(0,0){\scalebox{.75}{\includegraphics{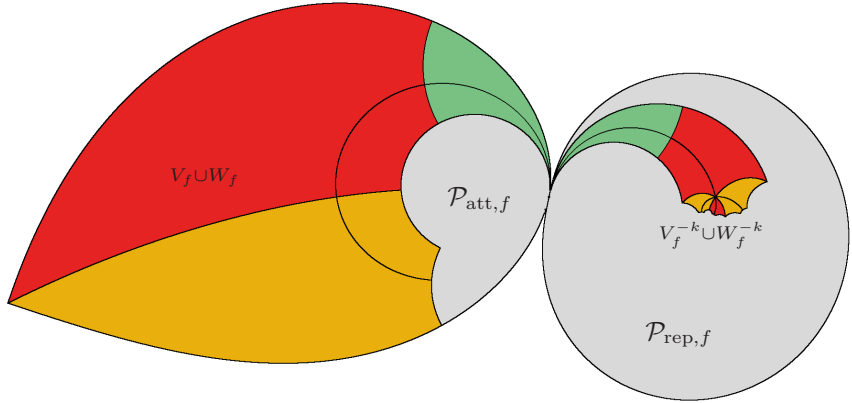}}}
 \put(75,105){$\scriptstyle V_f\cup W_f$}
 \put(260,83){$\scriptstyle V_f^{-k}\cup W_f^{-k}$}
 \put(255,45){$\Pet_{\rep,f}$}
 \put(180,96){$\Pet_{\att,f}$}
\end{picture}
\caption{\label{fig_vfwfinpetrep}%
If $k$ is large enough, $V_f^{-k}\cup W_f^{-k}$ is contained in the
repelling petal $\Pet_{\rep,f}$.}
\end{figure}

\begin{lemma}[see figure \ref{fig_vfwfinpetrep}]\label{lemma_r2k0}
There is an $R_2>0$ such that for all $f\in \sF_0$, the set
$\Phi_f(\Omega_{\rep,f})$ contains the half-plane $\{w\in
\C~;~\re~w<-R_2\}$. There is an integer $k_0>0$ such that for all
$k\geq k_0$, we have
\[V_f^{-k}\cup W_f^{-k}\subset  \bigl\{z\in
\Omega_{\rep,f}~;~\re\bigl(\Phi_{\rep,f}(z)\bigr)<-R_2\bigr\}.\]
\end{lemma}

\begin{rema}
Of course, $R_2$ can be replaced by any $R_3\geq R_2$, replacing if
necessary $k_0$ by $k_1\eqdef k_0+\lfloor R_3-R_2\rfloor +1$.
\end{rema}

\begin{proof}
For all $f\in \sF_0$, $\Phi_f(\Omega_{\rep,f})$ contains a left
half-plane. The existence of $R_2$ follows from the compactness of
$\sF_0$.

By Inou and Shishikura's result, we know that for all $f\in \sF_0$
there is an integer $k>0$ such that $W_f^{-k}$ is relatively compact
in $\Omega_{\rep,f}$. It follows from the compactness of $\sF_0$
that there is an integer $k_1>0$ and a constant $M$, such that for
all $f\in \sF_0$, $W_f^{-k_1}\subset \Omega_{\rep,f}$ and
\[\sup_{w\in
W_f^{-k_1}}\re\bigl(\Phi_{\rep,f}(w)\bigr)<M.\] Set $k_0\eqdef
k_1+M+\lfloor R_2\rfloor +3$. Then,
\[\sup_{w\in
W_f^{-k_0}}\re\bigl(\Phi_{\rep,f}(w)\bigr)<-R_2-2.\] We will show
that we then automatically have
\begin{equation}\label{eq_Wfinpetimpliesvfinpet}
V_f^{-k_0}\subset \Omega_{\rep,f}\quad\text{and}\quad \sup_{w\in
V_f^{-k_0}}\re\bigl(\Phi_{\rep,f}(w)\bigr)<-R_2.\end{equation} It
will follow immediately that \[\forall k\geq k_0\text{ and }\forall
w\in V_f^{-k}\cup W_f^{-k}, \qquad
\re\bigl(\Phi_{\rep,f}(w)\bigr)<-R_2,\] which will conclude the
proof of the lemma.

In order to get (\ref{eq_Wfinpetimpliesvfinpet}), we fix $f\in
\sF_0$ and consider $k\geq k_0$ large enough so that
$V_f^{-k}\subset \Omega_{\rep,f}$ (this is possible thanks to Inou
and Shishikura). Note that
\[\sup_{w\in W_f^{-k}}\re\bigl(\Phi_{\rep,f}(w)\bigr)<-R_2-2-k+k_0.\]
Denote by $g:\overline V_f\to \overline {V_f^{-k}}$ the inverse
branch of $f^{\circ k}:\overline {V_f^{-k}}\to \overline V_f$. Set
\[B\eqdef \bigl\{w\in \C~;~0<\re(w)<2~\text{and}~0<\im(w)\bigr\}.\]
Note that $B=\Phi_{\att,f}(V_f)$. Consider the map $\Psi:\overline
B\to \C$ defined by
\[\Psi\eqdef \Phi_{rep,f}\circ g\circ \Phi_{\att,f}^{-1}.\]
Since $\Psi$ commutes with translation by $1$, so that $\Psi(w)-w$
is $1$-periodic, the maximum modulus principle yields
\[\sup_{w\in B}\re\bigl(\Psi(w)-w\bigr) = \sup_{w\in
[0,2]}\re\bigl(\Psi(w)-w\bigr).\] Note that
\[g\circ \Phi_{\att,f}^{-1}\bigl([0,2]\bigr)\subset W_f^{-k}\]
and thus
\[\sup_{w\in
[0,2]}\re\bigl(\Psi(w)-w\bigr)<-R_2-2-k+k_0.\] Hence,
\[\sup_{w\in V_f^{-k}}\re\bigl(\Phi_{\rep,f}(w)\bigr) = \sup_{w\in
B}\re\bigl(\Psi(w)\bigr) < -R_2-k+k_0.\] It now follows that
\[\sup_{w\in V_f^{-k_0}}\re\bigl(\Phi_{\rep,f}(w)\bigr) < -R_2.\]
This completes the proof of (\ref{eq_Wfinpetimpliesvfinpet}) and of
lemma \ref{lemma_r2k0}.
\end{proof}

\subsubsection{Perturbed Fatou coordinates\label{subsubperturbedfatou}}

For $\a\in \R$, we denote by $\sF_\a$ the set of maps of the form
$z\mapsto f(e^{2i\pi \a}z)$ with $f\in \sF_0$. If $A$ is a subset of
$\R$, we denote by $\sF_A$ the set
\[\sF_A\eqdef \bigcup_{\alpha\in A} \sF_\a.\]
If $f\in \sF_{\left[0,1\right[}$, we denote by $\a_f\in
\left[0,1\right[$ the rotation number of $f$ at $0$, i.e. the real
number $\a_f\in \left[0,1\right[$ such that
\[f'(0)=e^{2i\pi \a_f}.\]

\begin{lemma}
There exist $\eps_0\in \left]0,1\right[$ and $r>0$ such that for all
$f\in \sF_{\left[0,\eps_0\right[}$, the map $f$ has two fixed points
in $D(0,r)$ (counting multiplicities), one at $z=0$ the other one
denoted by $\sigma_f$. The map
$\sigma:\sF_{\left[0,\eps_0\right[}\to D(0,r)$ defined by $f\mapsto
\sigma_f$ is continuous.
\end{lemma}

\begin{proof}
According to Inou and Shishikura, maps $f\in \sF_0$ have a double
fixed point at $0$. By compactness of $\sF_0$, there is an $r'>0$
such that maps $f\in \sF_0$ have only $2$ fixed points in $D(0,r')$.
Choose $r\in \left]0,r'\right[$. By Rouché's theorem and by
compactness of $\sF_0$, there is an $\eps_0>0$ such that maps $f\in
\sF_{\left[0,\eps_0\right[}$ have exactly two fixed points in
$D(0,r)$. The result follows easily.
\end{proof}

The following results are consequences of results in \cite{s}, the
compactness of the class $\sF_0$ and the results of the previous
paragraph.

\begin{figure}[htbp]
\begin{picture}(350,112)%
\put(-10,0){\scalebox{0.39}{\includegraphics{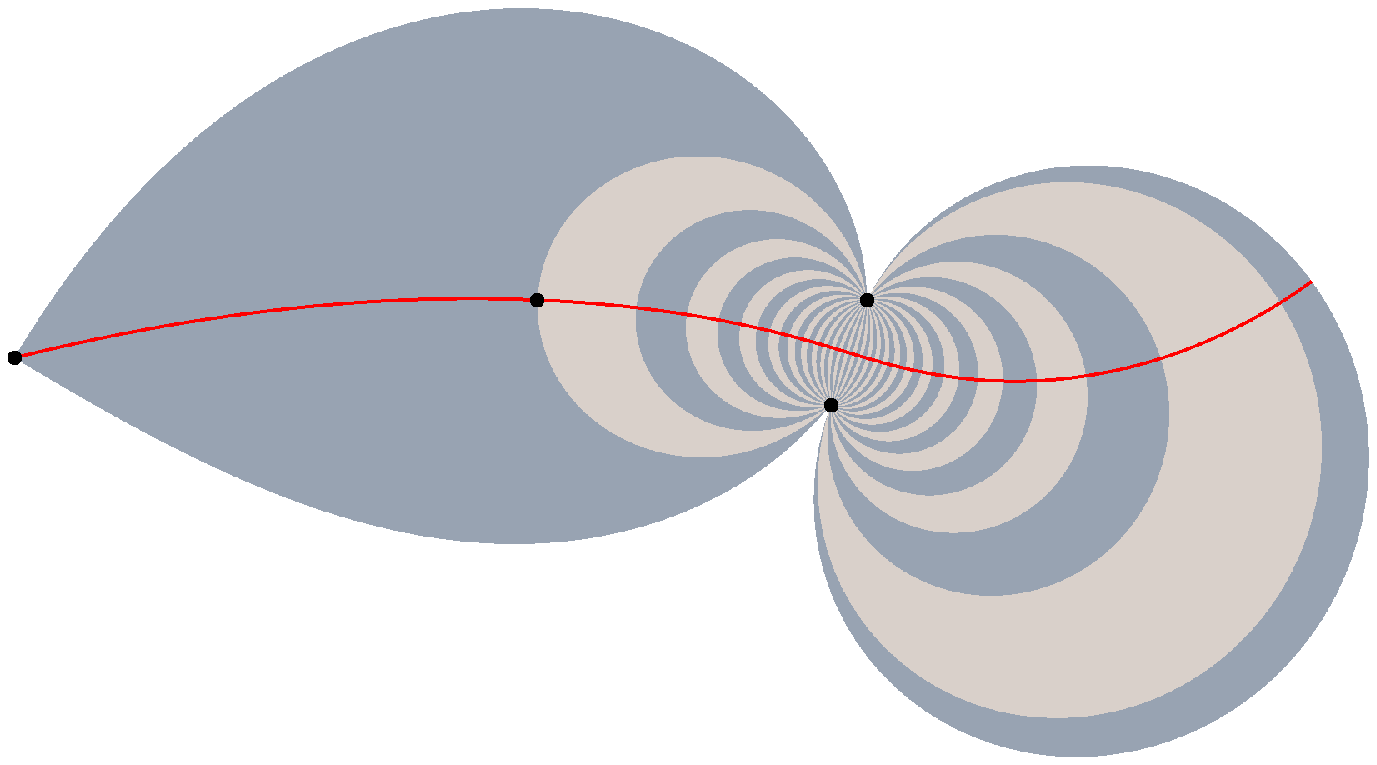}}}%
 \put(191,0){\scalebox{0.44}{\includegraphics{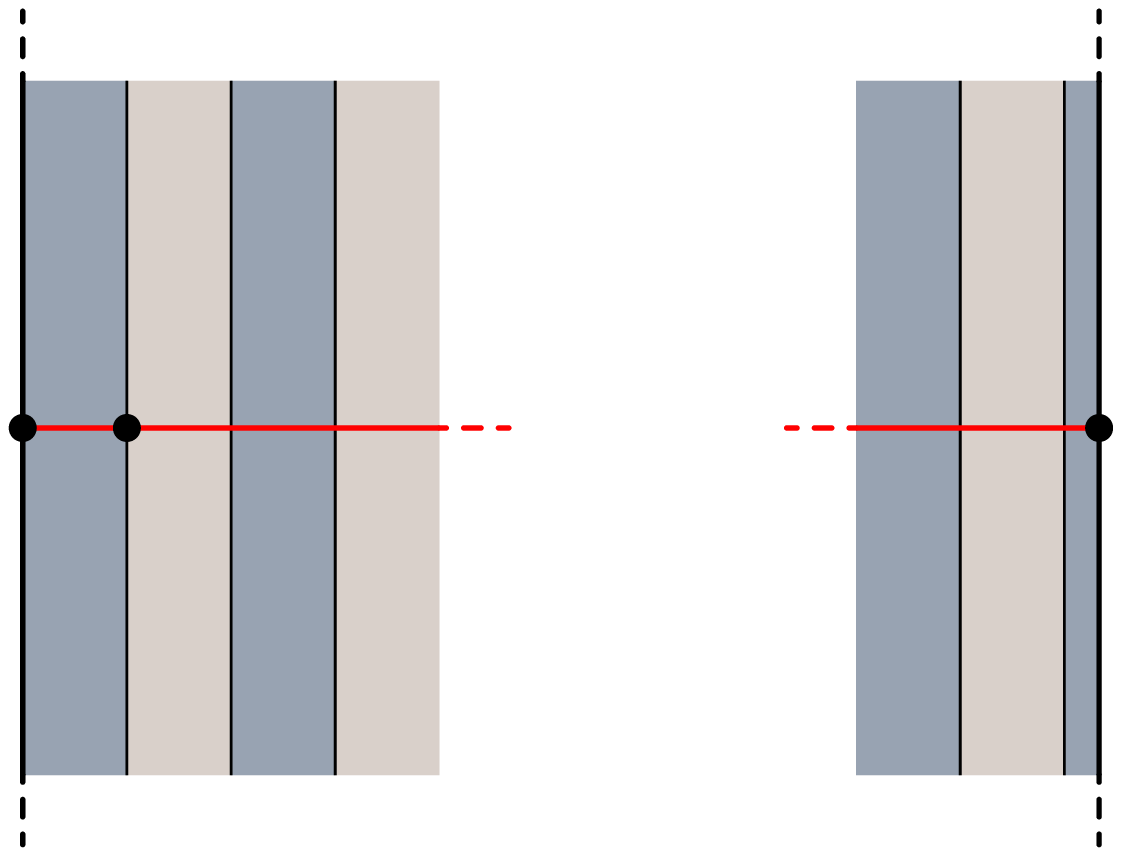}}}%
 \put(167,56){$\overset{\Phi_f}{\longrightarrow}$}%
 \put(50,75){$\Pet_f$}%
 \put(0,60){$\omega_f$}%
 \put(63,55){$v$}%
 \put(187,47){$0$}%
 \put(203,47){$1$}%
 \put(333,47){$\frac{1}{\a_f}-R_3$}%
 \put(12.62,54.6){\circle*{4}}
 \put(71.51,61.23){\circle*{4}}
 \put(104.27,49.14){\circle*{4}}
 \put(108.56,60.84){\circle*{4}}
\end{picture}
\caption{\label{fig_defpf}The perturbed petal $\Pet_f$ whose image
by the perturbed Fatou coordinate $\Phi_f$ is the strip
$\bigl\{0<\re(w)<1/\a_f-R_3\bigr\}$. }
\end{figure}

\begin{proposition}[see figure \ref{fig_defpf}]\label{prop_perturbedpetal}
There are constants $K>0$, $\eps_1>0$ and $R_3\geq R_2$ with
$1/\eps_1-R_3>1$, such that for all $f\in
\sF_{\left]0,\eps_1\right[}$ the following holds.
\begin{enumerate}
\item There is a Jordan domain $\Pet_f\subset U_f$ (a perturbed petal) containing $v$, bounded by two arcs
joining $0$ to $\sigma_f$ and there is a branch of argument defined
on $\Pet_f$ such that
\[\sup_{z\in \Pet_f} \arg(z) - \inf_{z\in \Pet_f} \arg(z) <K.\]

\item There is a univalent map $\Phi_f:{\cal P}_f\to \C$ (a perturbed Fatou coordinate)
such that
\begin{itemize}
\item $\Phi_f(v)=1$,

\item
$\ds \Phi_f(\Pet_f) = \bigl\{w\in \C~;~ 0<\re(w)<1/\a_f-R_3\bigr\}$,

\item
$\im\bigl(\Phi_f(z)\bigr)\to +\infty$ as $w\to 0$ and
$\im\bigl(\Phi_f(z)\bigr)\to -\infty$ as $w\to \sigma_f$ and

\item  $\Phi_f\circ f(z) = \Phi_f(z)+1$
when $z\in \Pet_f$ and $\re\bigl(\Phi_f(z) \bigr)<1/\a_f-R_3-1.$
\end{itemize}
\setcounter{saveenum}{\value{enumi}}
\end{enumerate}
For $f\in \sF_0$, we set
\[\Pet_{\rep,f}\eqdef \bigl\{z\in
\Omega_{\rep,f}~;~\re\bigl(\Phi_{\rep,f}(z)\bigr)<-R_3\bigr\}.\]
\begin{enumerate}
\setcounter{enumi}{\value{saveenum}}
\item If $(f_n)$ is a sequence of maps
in $\sF_{\left]0,\eps_1\right[}$ converging to a map $f_0\in \sF_0$,
then
\begin{itemize}
\item any compact $K\subset \Pet_{\att,f_0}$ is
contained in $\Pet_{f_n}$ for $n$ large enough and the sequence
$(\Phi_{f_n})$ converges to $\Phi_{\att,f_0}$ uniformly on $K$, and

\item any compact $K\subset \Pet_{\rep,f_0}$ is contained in
$\Pet_{f_n}$ for $n$ large enough and the sequence
$(\Phi_{f_n}-\frac{1}{\a_{f_n}})$ converges to $\Phi_{\rep,f_0}$
uniformly on $K$.
\end{itemize}
\end{enumerate}
\end{proposition}

\begin{proof}
Thanks to the compactness of the class $\sF_0$, it is enough to show
that if $(f_n)$ is a sequence of maps in $\sF_{\left]0,1\right[}$
converging to a map $f_0\in \sF_0$, there is a number $R_3\geq R_2$
such that properties (1), (2) and (3) hold.

So, assume $f_n$ is such a sequence, and for simplicity, write
$\a_n$, $\sigma_n$, \ldots instead of $\a_{f_n}$, $\sigma_{f_n}$,
\ldots

Let $\tau_n:\C\to\P\setminus \{0,\sigma_n\}$ be the universal
covering given by
\[\tau_n(w) \eqdef \frac{\sigma_n}{1-e^{-2i\pi \a_n w}}\]
so that
\[\tau_n(w)\underset{\im(w)\to +\infty}\longrightarrow
0\quad\text{and}\quad \tau_n(w)\underset{\im(w)\to
-\infty}\longrightarrow \sigma_n.\] Denote by $T_n:\C\to \C$ the
translation
\[T_n:w\mapsto w-\frac{1}{\a_n}.\]
Recall that $f_0(z) = z + c_0 z^2 + {\cal O}(z^3)$ with $c_0\neq 0$,
and
\[\tau_0(z)\eqdef  -\frac{1}{c_0 z}.\]

The following observations follow from \cite{sh}. We let $R_0$ and
$R_1$ be the constants introduced in paragraph \ref{subsubfatou}.
\begin{enumerate}
\item
The sequence $(\tau_n)$ converges to $\tau_0$ uniformly on every
compact subset of $\C^*$.

\item If $n$ is sufficiently large, there is a map $F_n:{\cal D}_n\to
\C$, defined and univalent in
\[{\cal D}_n\eqdef \C\setminus \bigcup_{k\in \Z}\overline D(k/\a_n,R_0)\]
which satisfies
\begin{itemize}
\item $f_n\circ \tau_n = \tau_n\circ F_n$,

\item $F_n(w)-w$ is $1/\a_n$-periodic (or equivalently, $F_n\circ
T_n = T_n\circ F_n$),

\item $F_n(w)-w\to 1$ as $\im(w)\to +\infty$.
\end{itemize}

\begin{rema}
This lift $F_n$ of $f_n$ may be defined by
\[F_n(w)\eqdef w+ \frac{1}{2i\pi \a_n} \log
\left(\frac{f_n(z)-\sigma_n}{f_n(z)}\cdot
\frac{z}{z-\sigma_n}\right)\quad\text{with}\quad z = \tau_n(w).\]
\end{rema}

\item As $n$ tends to $+\infty$, the sequence $(F_n)$ converges to
$F_0$ uniformly on every compact subset of $\C\setminus \overline
D(0,R_0)$.

\item The set
\[\Omega^n\eqdef \left\{w\in \C~;~\re(w)>R_1-\bigl|\im(w)\bigr|\text{ and }
\re(w)<\frac{1}{\a_n}-R_1+\bigl|\im(w)\bigr|\right\}\] is contained
in ${\cal D}_n$ (see figure \ref{fig_omeganfatou}.

\begin{figure}[htbp]
\begin{picture}(300,157)
 \put(0,0){\scalebox{.7}{\includegraphics{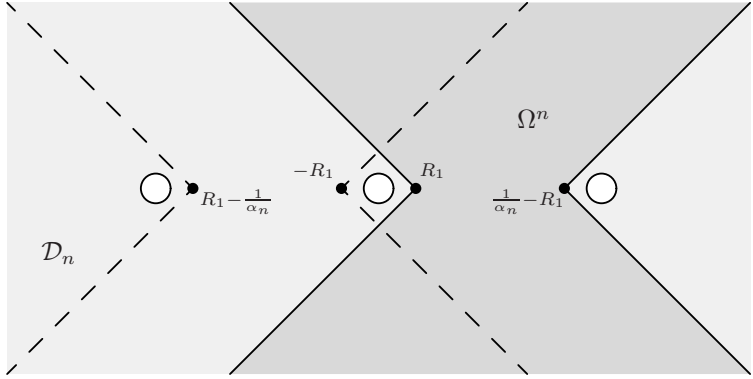}}}
 \put(200,100){$\Omega^n$}
 \put(190,71){$\scriptstyle \frac{1}{\a_n}-R_1$}
 \put(163,82){$\scriptstyle R_1$}
 \put(115,82){$\scriptstyle -R_1$}
 \put(80,71){$\scriptstyle R_1-\frac{1}{\a_n}$}
 \put(20,50){${\cal D}_n$}
\end{picture}
\caption{\label{fig_omeganfatou}The domain ${\cal D}_n$ (grey) is
the complement of a union of disks and the {\em hourglass}
$\Omega^n$ (drak grey) is contained in ${\cal D}_n$. }
\end{figure}

\item Remember that for all $w\in \C\setminus D(0,R_0)$,
\[\bigl|F_0(w)-w-1\bigr|<\frac{1}{4}\quad\text{and}\quad
\bigl|F_0'(w)-1\bigr|<\frac{1}{4}.\] It follows from the convergence
of $(F_n)$ to $F_0$ that if $n$ is sufficiently large, then for all
$w\in \Omega^n$,
\[\bigl|F_n(w)-w-1\bigr|<\frac{1}{4}\quad\text{and}\quad
\bigl|F_n'(w)-1\bigr|<\frac{1}{4}.\]

\item Increasing $n$ if necessary, we may assume that
$1/\a_n>2R_1+2$. Then, there is a univalent map $\Phi^n:\Omega^n\to
\C$, called a perturbed Fatou coordinate for $F_n$, such that
\[\Phi^n\circ F_n (w) = F_n(w)+1\]
when $w\in \Omega^n$ and $F_n(w)\in \Omega^n$. This map is unique up
to post-composition with a translation.

\item \label{point_convergenceliftedFatou}
Remember that there is a $k$ such that $f_0^{\circ k}(\omega_0)\in
\Omega_\att$, with $\omega_0$ the critical point of $f_0$. For $n$
large enough, $f_n^{\circ k}(\omega_n)$ is in $\tau_n(\Omega^n)$.
There is a point $w_n\in \Omega^n$ such that
\[\tau_n(w_n)=f_n^{\circ
k}(\omega_n)\quad\text{with}\quad w_n\underset{n\to
+\infty}\longrightarrow \tau_0^{-1}\bigl(f_0^{\circ
k}(\omega_0)\bigr).\] We can normalize $\Phi^n$ by $\Phi^n(w_n) =
k$. Then,
\[\Phi^n\underset{n\to
+\infty}\longrightarrow \Phi^\att_0\] uniformly on every compact
subset of $\Omega^\att$. Due to the normalization
$\Phi^\att_0(w)-\Phi^\rep_0(w)\to 0$ as $\im(w)\to +\infty$ with
$w\in \Omega^\att\cap \Omega^\rep$, we have \[T_n\circ \Phi^n\circ
T_n^{-1}\underset{n\to +\infty}\longrightarrow \Phi^\rep_0\]
uniformly on every compact subset of $\Omega^\rep$.
\end{enumerate}

Coming back to the $z$-coordinate is not immediate. Indeed, the map
$\tau_n$ is not injective on $\Omega^n$ and we cannot define a Fatou
coordinate for $f_n$ on $\tau_n(\Omega^n)$. We will instead restrict
to a subset $\Pet^n\subset \Omega^n$ whose image by $\Phi^n$ is a
vertical strip and on which $\tau_n$ is injective. The precise
statement is the following. The proof is given in appendix
\ref{appendix_fatou}. It is a consequence of results in \cite{sh},
but is not stated there.

\begin{lemma}[see figure \ref{fig_defunpn}]
If $K>0$ and $R \geq R_2$ are sufficiently large, then for $n$ large
enough:
\begin{itemize}
\item
$\Phi^n(\Omega^n)$ contains the vertical strip
\[U^n\eqdef \bigl\{w\in \C~;~R<\re(w)<1/\a_n-R\bigr\}\]
and

\item $\tau_n$ is injective on $\Pet^n\eqdef (\Phi^n)^{-1}(U^n)$.

\item there is a branch of argument defined on $\tau_n(\Pet^n)$ such that
\[\sup_{z\in \tau_n(\Pet^n)} \arg(z) - \inf_{z\in \tau_n(\Pet^n)} \arg(z) <K.\]
\end{itemize}
\end{lemma}

\begin{figure}[htbp]
\begin{picture}(330,210)
 \put(0,110){\scalebox{.5}{\includegraphics{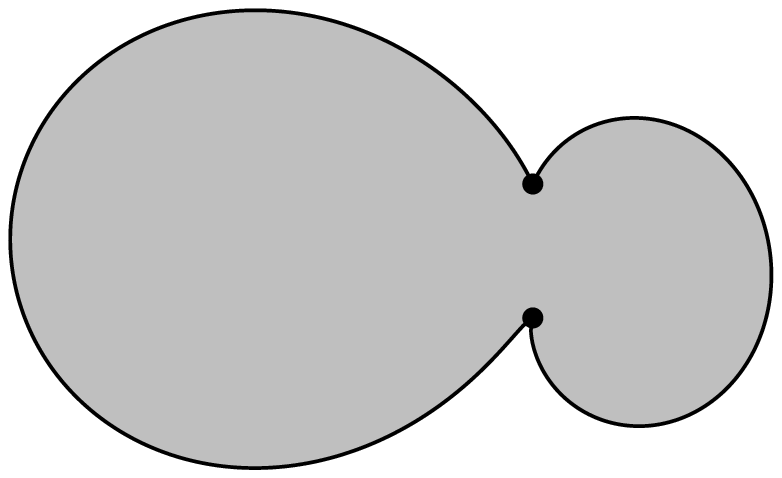}}}
 \put(130,110){\scalebox{.5}{\includegraphics{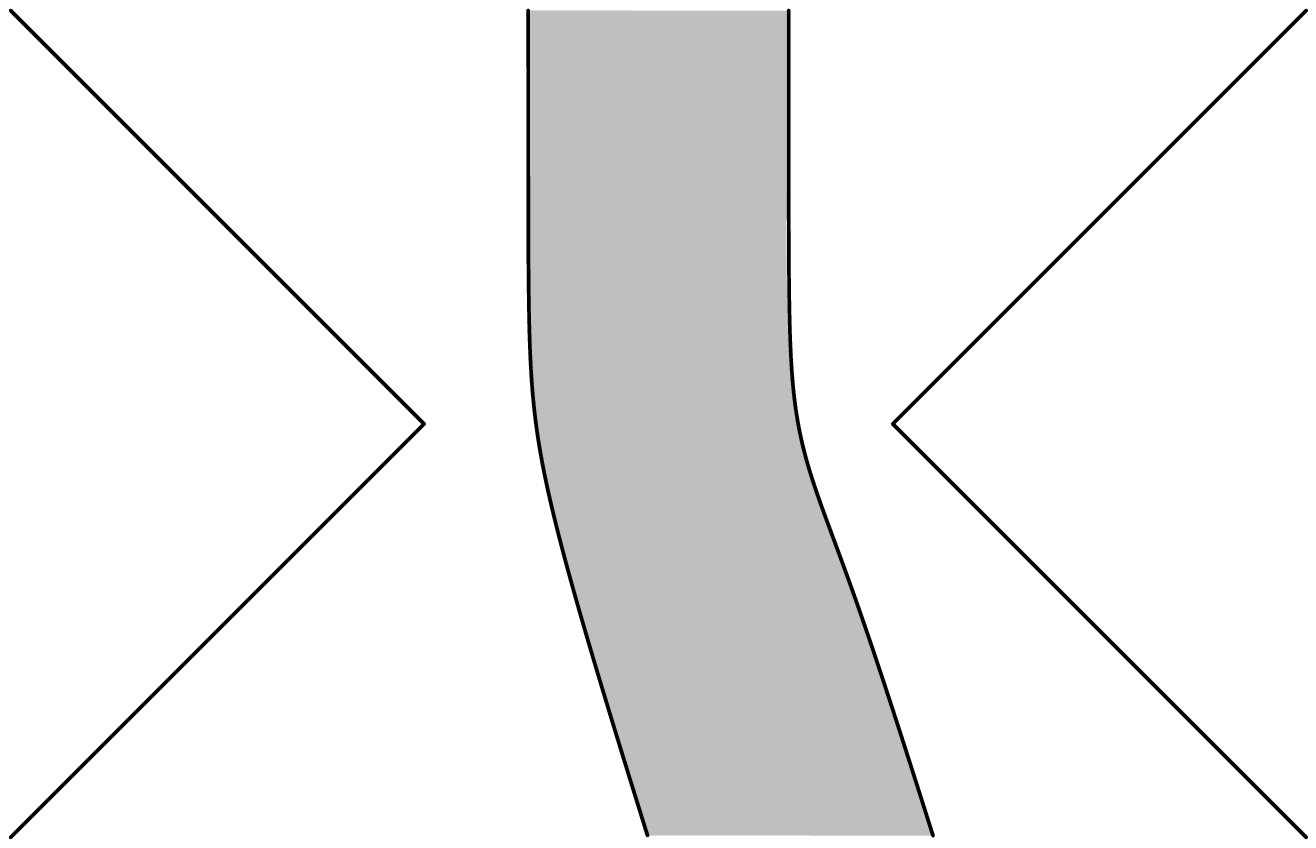}}}
 \put(50,0){\scalebox{.5}{\includegraphics{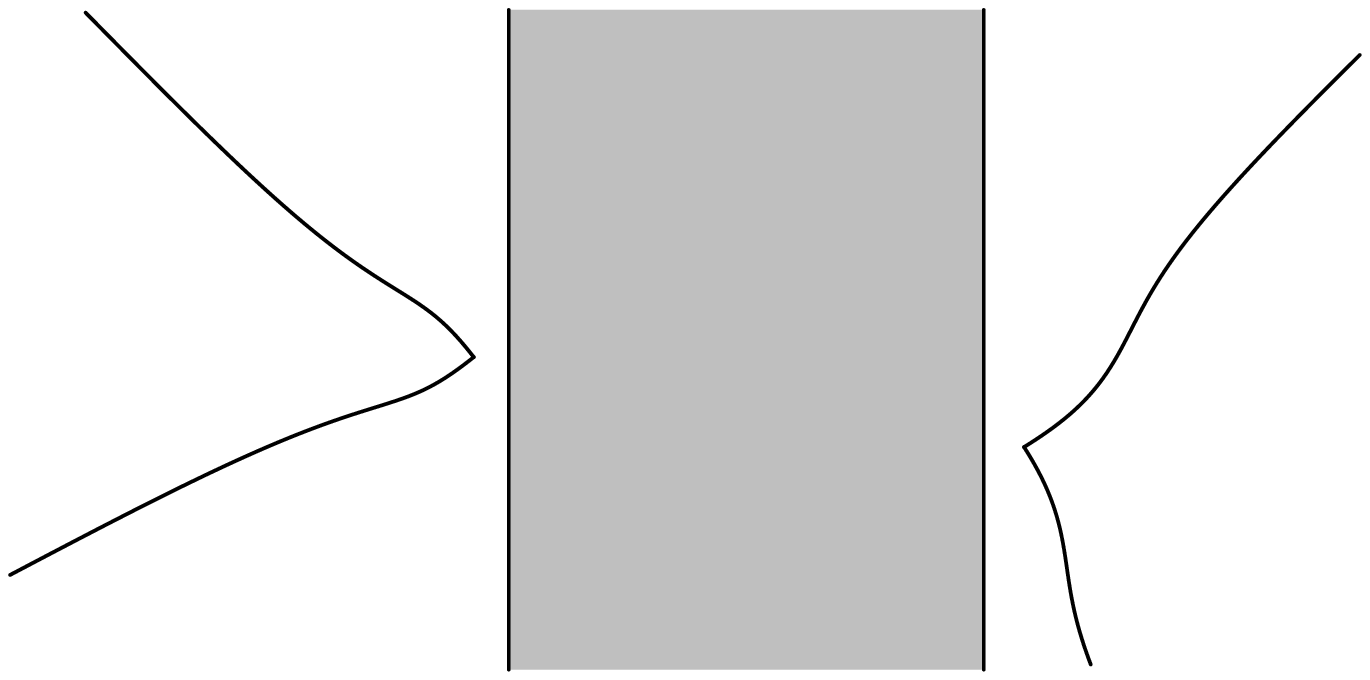}}}
 \put(140,155){$\tau_n$}
 \put(160,150){\vector(-1,0){30}}
 \put(180,105){$\Phi^n$}
 \put(220,113){\vector(-3,-1){40}}
 \put(30,155){$\tau_n(\Pet^n)$}
 \put(230,155){$\Pet^n$}
 \put(150,60){$U^n$}
 \put(87,158){$0$}
 \put(87,144){$\sigma^n$}
 \put(125,43){\circle*{4}}
 \put(193,43){\circle*{4}}
 \put(127,39){$\scriptstyle R$}
 \put(167,39){$\scriptstyle \frac{1}{\a_n}-R$}
\end{picture}
\caption{\label{fig_defunpn} The map $\tau_n$ is injective on
$\Pet^n\eqdef (\Phi^n)^{-1}(U^n)$.}
\end{figure}

Let $M>R$ be an integer. Note that
\[\bigl\{w\in
\C~;~\re(w)>M\bigr\}\subset \Phi_{\att,0}(\Omega_{\att,0})\] and
\[\bigl\{w\in
\C~;~\re(w)<-M\bigr\}\subset \Phi_{\rep,0}(\Omega_{\rep,0}).\] Set
\[\Pet'_0\eqdef \bigl\{z\in \Omega_{\att,0}~;~\re\bigl(\Phi_{\att,0}(z)\bigr)>M\bigr\}\cup
 \bigl\{z\in \Omega_{\rep,0}~;~\re\bigl(\Phi_{\rep,0}(z)\bigr)<-M\bigr\}\]
 and
\[\Pet'_n\eqdef \tau_n\bigl(\bigl\{w\in \Pet^n~;~M<\re\bigl(\Phi^n(w)\bigr)<1/\a_n-M\bigr\}\bigr).\]
For any $r>0$, if $n$ is sufficiently large so that $\sigma_n\in
D(0,r)$, then points with large (positive or negative) imaginary
part are mapped by $\tau_n$ in $D(0,r)$. It therefore follows from
point (7) above that $\overline{\Pet'_n} \to \overline {\Pet'_0}$ as
$n\to +\infty$.

Set
\[\Pet_0\eqdef \Pet_{\att,0}\cup
 \bigl\{z\in \Omega_{\rep,0}~;~\re\bigl(\Phi_{\rep,0}(z)\bigr)<-2M\bigr\}.\]
 Note that $\Pet_0$ is compactly contained in the domain of $f_0^{\circ
 M}$ and that $f_0^{\circ M}:\Pet_0\to \Pet'_0$ is an isomorphism.
In addition, for $n$ sufficiently large, $f_n^{\circ M}$ does not
have any critical value in $\Pet'_n$.

It follows from Rouché's theorem that for $n$ large enough, the
connected component $\Pet_n$ of $f_n^{-M}(\Pet'_n)$ which contains
$0$ in its boundary is relatively compact in the domain of $f_n$,
and $f_n^{\circ M}:\Pet_n\to \Pet'_n$ is an isomorphism. The
perturbed Fatou coordinate $\Phi_n:\Pet_n\to \C$ is defined by
\[\Phi_n(z) \eqdef \Phi^n(w)-M\quad\text{with}\quad w\in
\Pet^n\text{ and }\tau_n(w)=f_n^{\circ M}(z).\] In a simply
connected neighborhood of $\overline{\Pet'_0}$, the function
$f_0^{\circ M}(z)/z$ does not vanish (and extends by $1$ at $z=0$).
It follows that for $n$ large enough, there are branches of argument
of $f_n^{\circ M}(z)/z$ which are uniformly bounded on $\Pet_n$. It
is now easy to check that the proposition holds for $f_n$ with $n$
large enough.
\end{proof}

\subsubsection{Renormalization}

Recall that for maps $f\in \sF_0$ we defined sets $V_f\subset
\Pet_{\att,f}$ and $W_f\subset \Pet_{\att,f}$. We claimed (see lemma
\ref{lemma_r2k0}) that for $k\geq 0$ there are components $V_f^{-k}$
and $W_f^{-k}$ properly mapped by $f^{\circ k}$ respectively to
$V_f$ with degree $1$ and $W_f$ with degree $2$. In addition, there
is an integer $k_0>0$ such that \[\forall f\in \sF_0,\qquad
V_f^{-k_0}\cup W_f^{-k_0}\subset \Pet_{\rep,f}.\]

We will now generalize this to maps $f\in \sF_{\left]0,\eps\right[}$
with $\eps$ sufficiently small. If $f \in
\sF_{\left]0,\eps_1\right[}$, we set
\[V_f\eqdef \left\{z\in \Pet_f~;~\im\bigl(\Phi_f(z)\bigr)>0\text{ and }
0<\re\bigl(\Phi_f(z)\bigr)<2 \right\}\] and
\[W_f\eqdef \left\{z\in \Pet_f~;~
-2<\im\bigl(\Phi_f(z)\bigr)<2\text{ and }
0<\re\bigl(\Phi_f(z)\bigr)<2\right\}.\]

\begin{figure}[htbp]
\begin{picture}(336,192)%
 \put(0,-10){\scalebox{.7}{\includegraphics{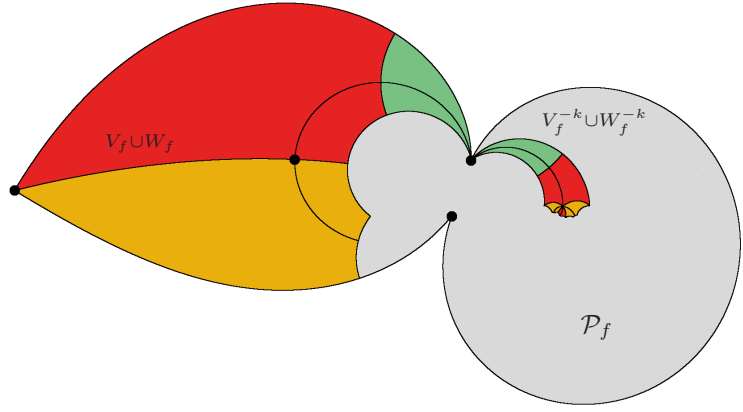}}}
 \put(75,105){$\scriptstyle V_f\cup W_f$}
 \put(240,113){$\scriptstyle V_f^{-k}\cup W_f^{-k}$}
 \put(255,35){$\Pet_f$}
\end{picture}
\caption{\label{fig_vfwfinpetf}%
If $k$ is large enough, $V_f^{-k}\cup W_f^{-k}$ is contained in the
perturbed petal $\Pet_f$.}
\end{figure}

\begin{proposition}[see figure \ref{fig_vfwfinpetf}]
There is a number $\eps_2>0$ and an integer $k_1\geq 1$ such that
for all $f\in \sF_{\left]0,\eps_2\right[}$ and  for all integer
$k\in [1,k_1]$,
\begin{enumerate}
\item the unique connected component $V_f^{-k}$ of $f^{-k}(V_f)$
which contains $0$ in its closure is relatively compact in $U_f$
(the domain of $f$) and $f^{\circ k}:V_f^{-k}\to V_f$ is an
isomorphism,

\item the unique connected component $W_f^{-k}$ of $f^{-k}(W_f)$
which intersects $V_f^{-k}$ is relatively compact in $U_f$ and
$f^{\circ k}:W_f^{-k}\to W_f$ is a covering of degree $2$ ramified
above $v$.

\item $V_f^{-k_1}\cup W_f^{-k_1}\subset\bigl\{z\in \Pet_f~;~
2<\re\bigl(\Phi_f(z)\bigr)<\frac{1}{\a_f}-R_3-5\bigr\}$.
\end{enumerate}
\end{proposition}

\begin{proof}
Set $k_1\eqdef k_0+7$. By compactness of $\sF_0$, there is an
$\eps_2>0$ such that for all $f\in \sF_{\left]0,\eps_2\right[}$,
properties (1) and (2) hold for all integers $k\in [1,k_1]$, and
further, $W_f^{-k_1}$ is contained in $\bigl\{z\in \Pet_f~;~
4<\re\bigl(\Phi_f(z)\bigr)<\frac{1}{\a_f}-R_3-7\bigr\}$.

To see that $V_f^{-k_1}$ is a subset of $\bigl\{z\in \Pet_f~;~
2<\re\bigl(\Phi_f(z)\bigr)<\frac{1}{\a_f}-R_3-5\bigr\}$, we proceed
as in the proof of lemma \ref{lemma_r2k0}.
\end{proof}

We now come to the definition of the renormalization of maps $f\in
\sF_{\left]0,\eps_2\right[}$.

\begin{is}[Main theorem 3]\label{theo_inoushishikura2}
If $f\in \sF_{\left]0,\eps_2\right[}$, the map
\[\Phi_f\circ f^{\circ k_1}\circ
\Phi_f^{-1}:\Phi_f\bigl(V_f^{-k_1}\cup W_f^{-k_1}\bigr)\to
\Phi_f\bigl(V_f\cup W_f\bigr)\] projects via $w\mapsto
-\frac{4}{27}e^{2i\pi w}$ to a map ${\cal R}(f)\in \sF_{-1/\a_f}$.
\end{is}

\begin{definition}
The map ${\cal R}(f)$ is called the renormalization of $f$.
\end{definition}

The construction we described also works for polynomials $P_\a$ with
$\a>0$ sufficiently close to $0$, i.e. the existence of perturbed
petals and perturbed Fatou coordinates, the existence of a
renormalization ${\cal R}(P_\a)$ which belongs to $\sF_{-1/\a}$ (the
only difference is that the critical value of $P_\a$ is not
normalized at $-4/27$). In the sequel, $\eps_2>0$ is chosen
sufficiently small so that for $\a\in \left]0,\eps_2\right[$, a map
$f$ which either is a polynomial $P_\a$, or belongs to $\sF_\a$, has
a renormalization ${\cal R}(f)\in \sF_{-1/\a}$.

\subsubsection{Renormalization tower}

Assume $1/N<\eps_2$. Denote by $\irrat_{\geq N}$ the set:
\[\irrat_{\geq N}\eqdef\bigl\{\a=[{\rm a}_0,{\rm a}_1,{\rm a}_2,\ldots]\in
\R\setminus \Q~;~ {\rm a}_k\geq N\text{ for all }k\geq 1\bigr\}.\]
Assume $\a=[{\rm a}_0,{\rm a}_1,{\rm a}_2,\ldots]\in \irrat_{\geq
N}$. For $j\geq 0$, set
\[\alpha_j\eqdef [0,{\rm a}_{j+1},{\rm a}_{j+2},\ldots].\]
Note that for all $j\geq 1$,
\[\a_{j+1}=\frac{1}{\a_j}-\left\lfloor
\frac{1}{\a_j}\right\rfloor.\] The requirement
$\alpha\in\irrat_{\geq N}$ translates into
\[\forall j,\qquad \alpha_j\in]0,1/N[.\]
Denote by $p_j/q_j$ the approximants to $\a_0$ given by the
continued fraction algorithm.

Now, if either $f_0=P_\a$ or $f_0\in \sF_\a$, we can define
inductively an infinite sequence of renormalizations, also called a
\emph{renormalization tower}, by
\[ f_{j+1}\eqdef s\circ
{\cal R}(f_{j})\circ s^{-1},\] the conjugacy by $s:z\mapsto \bar z$
being introduced so that \[f_j'(0)=e^{2i\pi\alpha_j}.\] It will be
convenient to define
\[\begin{array}{rccl}
\Exp:&\C&\to &\C^*\\
&w&\mapsto& -\frac{4}{27} s( e^{2i\pi w}). \end{array}\] For $j\geq
0$, we define
\[\phi_{j}\eqdef \Exp\circ \Phi_{f_j}:\Pet_{f_j} \to \C.\]
The map $\phi_j$ goes from the $j$-th level of the renormalization
tower to the next level.

We now want to relate the dynamics of maps at different levels of
the renormalization tower. For this purpose, we will use the
following lemma.

\begin{lemma}\label{lemma_widthlogpet}
There is a constant $K>0$ such that for all $f\in
\sF_{\left]0,\eps_2\right[}$, there is an inverse branch of $\Exp$
which is defined on $\Pet_f$ and takes its values in the strip
$\bigl\{w\in \C~;~0<\re(w)<K\bigr\}$.
\end{lemma}

\begin{proof}
This is an immediate consequence of proposition
\ref{prop_perturbedpetal} part (1).
\end{proof}

From now on, we assume that $N$ is sufficiently large so that
\begin{equation}\label{eq_chooseN}
\frac{1}{N}<\eps_2\quad \text{and}\quad \frac{1}{N}-R_3>K.
\end{equation}

Then, according to lemma \ref{lemma_widthlogpet}, for all $j\geq 1$,
there is an inverse branch $\psi_j$ of $\phi_{j-1}$ defined on the
perturbed petal $\Pet_{f_j}$ with values in $\Pet_{f_{j-1}}$ (there
are several possible choices, we choose any one).

\begin{figure}[htbp]
\begin{picture}(340,230)
 \put(-10,115){\scalebox{.4}{\includegraphics{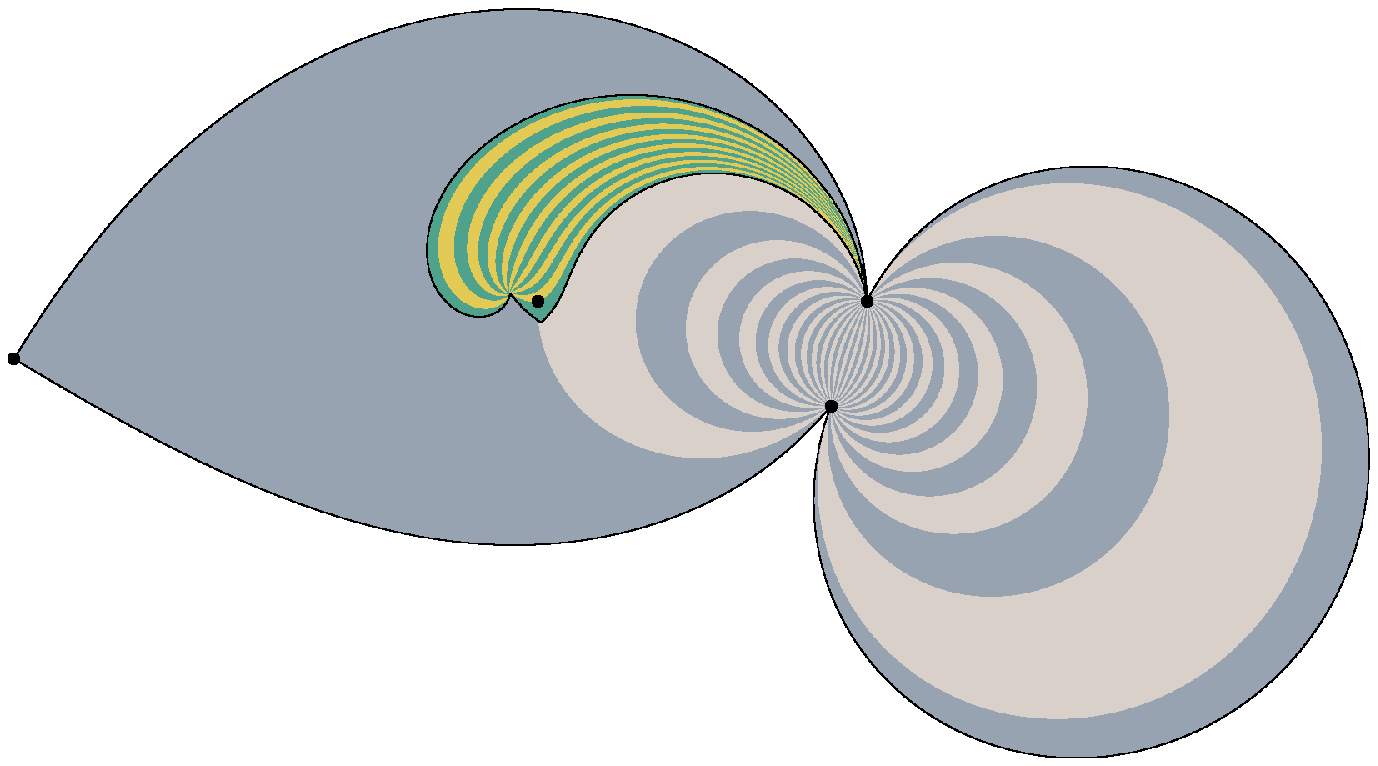}}}
 \put(200,120){\scalebox{.35}{\includegraphics{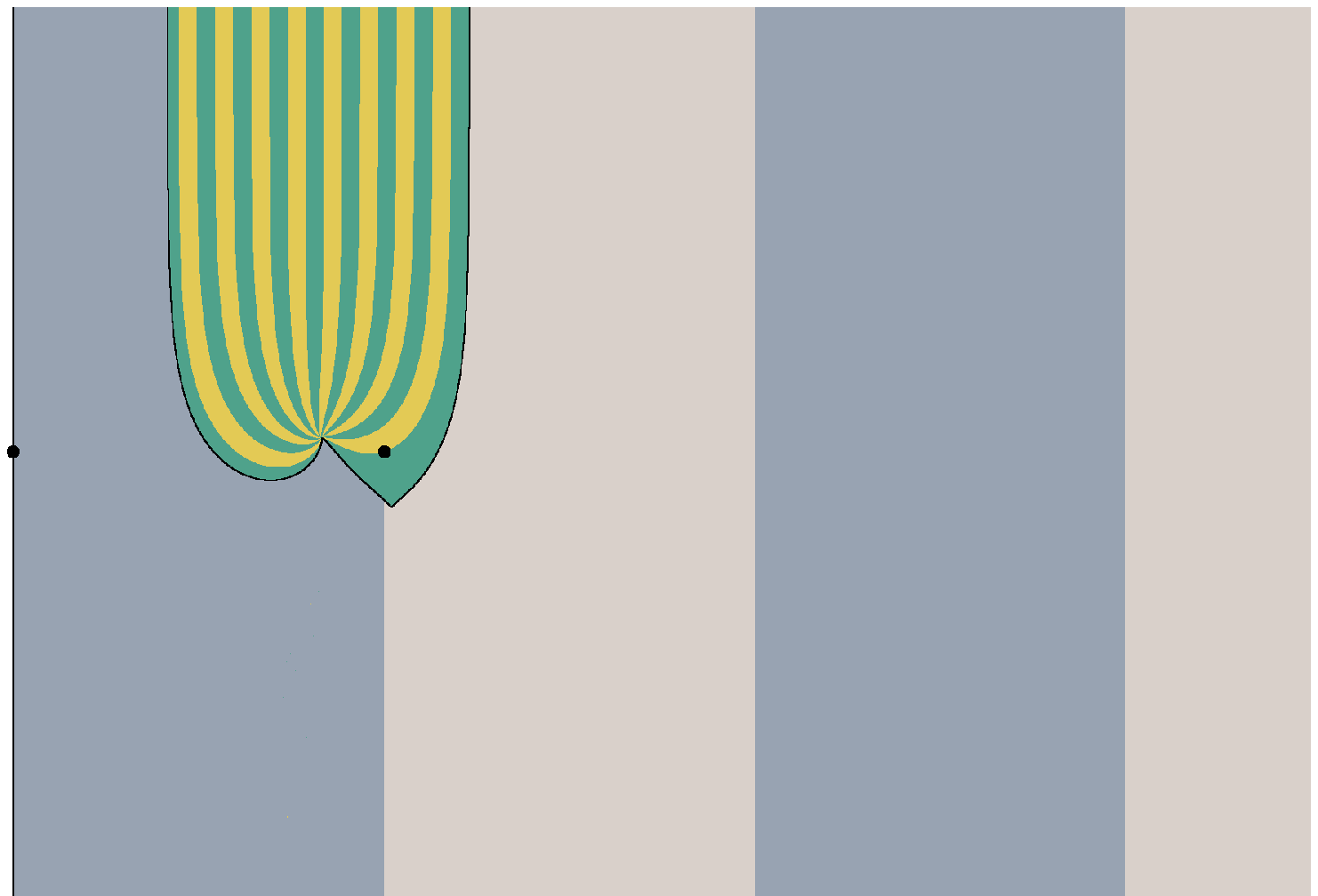}}}
 \put(100,0){\scalebox{.4}{\includegraphics{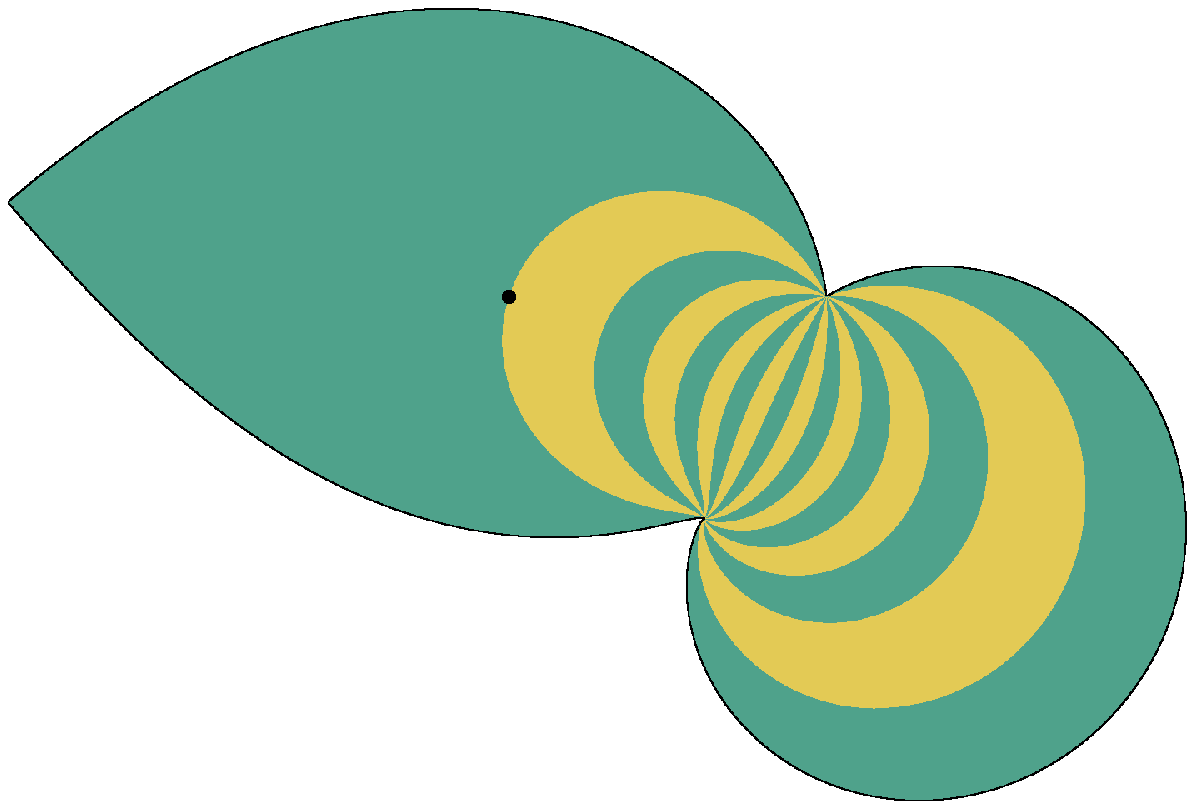}}}
 \put(130,30){$\Pet_{f_{j+1}}$}
 \put(10,150){$\Pet_{f_j}$}
 \put(181,170){\vector(1,0){20}}
 \put(183,175){$\Phi_{f_j}$}
 \put(250,166){\vector(-1,-2){42}}
 \put(225,100){$\Exp$}
 \put(70,92){\rotatebox{120}{\vector(1,0){104}}}
 \put(90,120){$\psi_{j+1}$}
\end{picture}
\caption{The branch $\psi_{j+1}$ maps $\Pet_{f_{j+1}}$ univalently
into $\Pet_{f_j}$. }
\end{figure}

The map
\[\Psi_{j}\eqdef \psi_{1}\circ \psi_{2}\circ \ldots\circ \psi_{j}\]
is then defined and univalent on  $\Pet_{f_{j}}$ with values in the
dynamical plane of the polynomial $f_0$.

Remember that
\[\Phi_{f_j}(\Pet_{f_j}) = \bigl\{w\in
\C~;~0<\re(w)<1/\a_j-R_3\bigr\}.\] Define $\Pet_j \subset
\Pet_{f_j}$ and $\Pet'_j \subset \Pet_{f_j}$ by
\[\Pet_j\eqdef \bigl\{z\in
\Pet_{f_j}~;~0<\re\bigl(\Phi_{f_j}(w)\bigr)<1/\a_j-R_3-1\bigr\}\]
and
\[\Pet'_j\eqdef \bigl\{z\in
\Pet_{f_j}~;~1<\re\bigl(\Phi_{f_j}(w)\bigr)<1/\a_j-R_3\bigr\}.\]
Note that $f_j$ maps $\Pet_{j}$ to $\Pet'_{j}$ isomorphically. Set
\[{\cal Q}_{j}\eqdef \Psi_{j}(\Pet_j)\quad\text{and}\quad
 {\cal Q}'_{j}\eqdef \Psi_{j}(\Pet'_j).\]

\begin{proposition}\label{prop_correspondanceorbit1}
The map $\Psi_{j}$ conjugates $f_j:\Pet_{j}\to \Pet'_{j}$ to
$f_0^{\circ  q_j}:{\cal Q}_{j}\to {\cal Q}'_{j}$.
\end{proposition}

In other words, we have the following commutative diagram:
\[\diagram
{\cal Q}_j \subset \Psi_j(\Pet_{f_j}) \rto^{f_0^{\circ q_j}} & {\cal
Q}'_j\subset \Psi_j(\Pet_{f_j})\\
\Pet_j\subset \Pet_{f_j} \rto_{f_j} \uto_{\Psi_j} & \Pet'_j\subset
\Pet_{f_j} \uto_{\Psi_j}.
\enddiagram\]

\begin{proof}
We must show that if $z_j\in \Pet_j$ and $z'_j\eqdef f_j(z_j)\in
\Pet'_j$, then the points $z_0\eqdef \Psi_j(z_j)$ and $z'_0\eqdef
\Phi_j(z'_j)$ are related by
\[z'_0 = f_0^{\circ q_j}(z_0).\]
Let us first show that there is an integer $k$ such that $z'_0 =
f_0^{\circ k}(z_0)$. Our proof is based on the following lemma.

\begin{lemma}\label{lemma_correspondanceorbit}
Assume $\ell \geq 0$, $w\in U_{f_{\ell+1}}$ and $w'\eqdef
f_{\ell+1}(w)$. Let $z\in \Pet_{f_\ell}$ and $z'\in \Pet_{f_\ell}$
be such that
\[\Exp\circ \Phi_{f_\ell}(z) = w\quad\text{and}\quad  \Exp\circ \Phi_{f_\ell}(z') =
w'.\] Then, there is an integer $k\geq 1$ such that
$z'=f_\ell^{\circ k}(z)$.
\end{lemma}

\begin{proof}
Let $z'_1\in \Pet_{f_\ell}$ be the unique point such that
\[\re\bigl(\Phi_{f_\ell}(z'_1)\bigr)\in
\left]0,1\right]\quad\text{and}\quad \Exp\circ \Phi_{f_\ell}(z'_1) =
w'.\] By definition of the renormalization $f_{\ell+1}$, there is a
point $z_1\in V_{f_\ell}^{-k_1}\cup W_{f_\ell}^{-k_1}$ such that
\[\Exp\circ \Phi_{f_\ell}(z_1) = w\quad\text{and}\quad f_\ell^{\circ
k_1}(z_1)=z'_1.\]
We then have
\[\Phi_{f_\ell}(z_1) = \Phi_{f_\ell}(z)+m_1\quad\text{and}\quad
\Phi_{f_\ell}(z')=\Phi_{f_\ell}(z'_1)+m'_1\] with $m_1\in \Z$ and
$m'_1\in \N$. If $m_1\geq 0$, we have
\[z_1=f_\ell^{\circ m_1}(z)\quad\text{and} z'=f_\ell^{\circ
m'_1}(z'_1).\] Since $k_1\geq 0$, we then have
\[z'= f^{\circ k}(z)\quad\text{with}\quad k\eqdef k_1+m_1+m'_1\geq 1.\]
If $m_1<0$, then $z=f_\ell^{\circ -m_1}(z'_1)$. However, for $m\leq
-m_1$, we have $f_\ell^{\circ m}(z'_1)\in \Pet_{f_\ell}$, and so,
$k_1\geq -m_1+1$. Thus, we can write
\[z'_1 = f_\ell^{\circ m_2}(z)\quad\text{with}\quad m_2\eqdef k_1+m_1\geq 1.\]
In that case,
\[z' = f^{\circ k}(z)\quad\text{with}\quad k\eqdef m_2+m'_1\geq 1.\]
\end{proof}

It follows by decreasing induction on $\ell$ from $j$ to $0$ that
for all $z_j\in \Pet_j$, there is an integer $k\geq 1$ such that
\[z'_0 = f_0^{\circ k}(z_0).\]
We will now show that we have a common integer $k$, valid for all
points $z_j\in \Pet_j$.

\begin{lemma}
There is an integer $k_0\geq 1$ such that for all point $z_j\in
\Pet_j$, we have
\[z'_0 = f_0^{\circ k_0}(z_0).\]
\end{lemma}

\begin{proof}
We will use the connectivity of $\Pet_j$. For $k\geq 1$, set
\[{\cal O}_k\eqdef \{z\in \Pet_j~;~f_0^{\circ k}\bigl(\Psi_j(z)\bigr)\text{
is defined}\bigr\}\] This is an open set. Set
\[X_k\eqdef \bigl\{z\in {\cal O}_k~;~f_0^{\circ
k}\bigl(\Psi_j(z)\bigr) = \Psi_j\bigl(f_j(z)\bigr)\}.\]

Note that for every component $O$ of ${\cal O}_k$, either $X_k\cap
{\cal O} = O$, or $X_k$ is discrete in $O$, in particular countable.
Indeed, $X_k$ is the set of zeroes of the holomorphic function
$f_0^{\circ k}\circ \Psi_j-\Psi_j\circ f_j:{\cal O}_k\to \C$.

Since
\[\Pet_j = \bigcup_{k\geq 1} X_k\]
there is a smallest integer $k_0\geq 1$ such that $X_{k_0}$ is not
countable. Then, there is a component $O$ of ${\cal O}_{k_0}$ such
that on $O$, we have $f_0^{\circ k_0}\circ \Psi_j = \Psi_j\circ
f_j$.

Since $O$ is a component of ${\cal O}_{k_0}$, we have
\[\partial O\cap \Pet_j \subset \C\setminus {\cal O}_{k_0}.\]
It follows that
\[\partial O\cap \Pet_j\subset X_1\cup \ldots X_{k_0-1}\]
since the remaining $X_k$'s are contained in ${\cal O}_{k_0}$. So,
$\partial O\cap \Pet_j$ is countable. This is only possible if
$\partial O\cap \Pet_j=\emptyset$ since in any neighborhood of a
point $z\in \C\setminus {\cal O}_{k_0}$, there are uncountably many
points in $\C\setminus {\cal O}_{k_0}$. As a consequence,
$O=\Pet_j$, which concludes the proof of the lemma.
\end{proof}

We must now show that $k_0=q_j$. Let $L_j\subset \Pet_j$ be the
curve defined by \[L_j\eqdef \bigl\{z\in
\Pet_j~;~\re\bigl(\Phi_{f_j}(z)\bigr)=1\bigr\}.\] Set $L'_j\eqdef
f_j(L_j)$, i.e. the curve
\[L'_j\eqdef \bigl\{z\in
\Pet_j~;~\re\bigl(\Phi_{f_j}(z)\bigr)=2\bigr\}.\]

Those curves both have an end point at $z=0$. They both have
tangents at $z=0$. Since the linear part of $f_j$ at $z=0$ is the
rotation of angle $\alpha_j$, the angle between $L_j$ and $L'_j$ at
$z=0$ is $\alpha_j$. It follows that the curves $\Psi_j(L_j)$ and
$\Psi_j(L'_j)$ have tangents at $z=0$ and the angle between those
curves is $\alpha_0\alpha_1\cdots \alpha_j$. So, the linear part of
$f_0^{\circ k_0}$ at $z=0$ is the rotation of angle
$\alpha_0\alpha_1\cdots \alpha_j$. It follows that $k_0=q_j$.
\end{proof}

Set
\begin{gather*}
D_j\eqdef V_{f_j}^{-k_1}\cup W_{f_j}^{-k_1},\quad
 D'_j\eqdef V_{f_j}\cup W_{f_j},
\\
C_j\eqdef \Psi_j(D_{j}) \quad \text{and} \quad  C'_j\eqdef
\Psi_j(D'_{j}).
\end{gather*}
Note that $f_j^{\circ k_1}$ maps $D_j$ to $D'_j$.

\begin{proposition}\label{prop_correspondanceorbit2}
The map $\Psi_{j}$ conjugates the map $f_j^{\circ k_1}:D_j\to D'_j$
to the map $f_0^{\circ(k_1q_j+q_{j-1})}:C_j\to C'_j$.
\end{proposition}

In other words, we have the following commutative diagram:
\[\diagram
C_j \subset \Psi_j(\Pet_{f_j}) \rrto^{f_0^{\circ(k_1q_j+q_{j-1})}} &
&
C'_j\subset \Psi_j(\Pet_{f_j})\\
D_j\subset \Pet_{f_j} \rrto_{f_j} \uto_{\Psi_j} & &D'_j\subset
\Pet_{f_j^{\circ k_1}} \uto_{\Psi_j}.
\enddiagram\]

\begin{proof} The proof is similar to the one of proposition
\ref{prop_correspondanceorbit1}.
\end{proof}

\subsubsection{Neighborhoods of the postcritical set}

We can now see that the post-critical set of maps $f\in \sF_\a$ with
$\a\in \irrat_{\geq N}$ is infinite.

\begin{proposition}[Inou-Shishikura]\label{prop_pcinfinite}
For all $\a\in \irrat_{\geq N}$ and all $f\in \sF_\a$, the
postcritical set of $f$ is infinite.
\end{proposition}

\begin{proof}
For $j\geq 1$, the map $f_j^{\circ k_1}:W_{f_j}^{-k_1}\to W_{f_j}$
is a ramified covering of degree $2$, ramified above $v$. Denote by
$w_j$ the critical point of this ramified covering. Set $w_0\eqdef
\Psi_j(w_j)$. According to proposition
\ref{prop_correspondanceorbit2}, we can iterate $f_0$ at least $k_1
q_j + q_{j-1}$ times at $w_0$, $w_0$ is a critical point of
$f_0^{\circ (k_1 q_j + q_{j-1})}$ and its critical value is
$\Psi_j(v)$. In particular, $\Psi_j(v)$ is a point of the
postcritical set of $f_0$.

Note that $v\in \Pet_j$. According to proposition
\ref{prop_correspondanceorbit1}, we can iterate $f_0$ at least $q_j$
times at $\Psi_j(v)$. This shows that we can iterate $f_0$ at least
$q_j$ times at $v$. Since $j\geq 1$ is arbitrary, the postcritical
set of $f_0$ is infinite.
\end{proof}

For every $\a\in \irrat_{\geq N}$, we are going to define a sequence
$(U_j)$ of open sets containing the post-critical set of $P_\a$. We
still use the notations of the previous paragraph. In particular,
for $j\geq 1$, the $j$-th renormalization of $f_0\eqdef P_\a$ has a
perturbed petal $\Pet_{f_j}$, a perturbed Fatou coordinate
\[\Phi_{f_j}:\Pet_{f_j}\to \bigl\{w\in \C~;~0<\re(w)<1/\a_j-R_3\bigr\}.\]
The set
\[D_j\eqdef V_{f_j}^{-k_1}\cup W_{f_j}^{-k_1}\subset \Pet_{f_j}\]
is mapped by $f_j^{\circ k_1}$ to
\[D'_j\eqdef \bigl\{z\in
\Pet_{f_j}~;~0<\re\bigl(\Phi_{f_j}(z)\bigr)<2\text{ and
}\im\bigl(\Phi_{f_j}(z)\bigr)>-2\bigr\}.\] There is a map $\Psi_j$,
univalent on $\Pet_{f_j}$, with values in the dynamical plane of
$P_\a$, conjugating $f_j^{\circ k_1}:D_j\to D'_j$ to $P_\a^{\circ
(k_1q_j+q_{j-1})}:C_j\to C'_j$ with
\[C_j\eqdef \Psi_j(D_j)\quad\text{and}\quad C'_j\eqdef
\Psi_j(D'_j).\]

\begin{definition}
For $\a\in \irrat_{\geq N}$ and $j\geq 1$ we set
\[U_j(\a)\eqdef \bigcup_{k=0}^{q_{j+1}+\ell q_j} P_\a^{\circ  k}(C_j)\]
where $\ell \eqdef k_1-\lfloor R_3\rfloor -4 \in\N$.
\end{definition}

\begin{figure}[htbp]
\centerline{\scalebox{.4}{\includegraphics{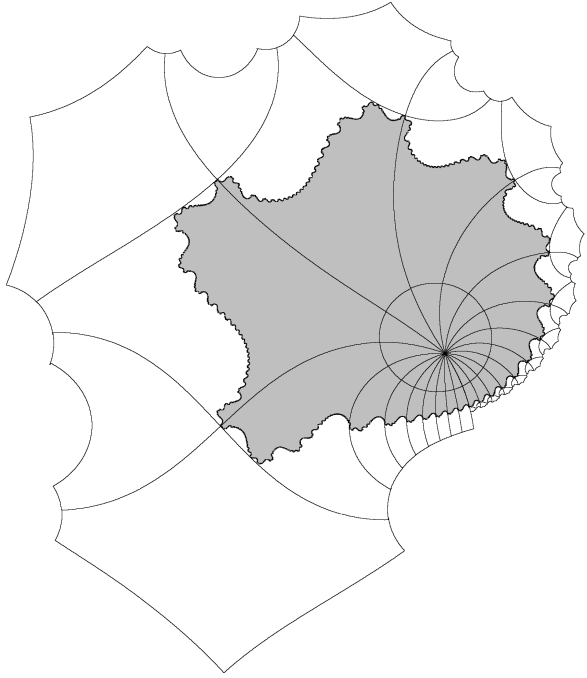}}}
\caption{If $f\in \sF_\a$ with $\a\in \irrat_{\geq N}$, the set
$U_1(f)$ contains the postcritical set $\PC(f)$. If $\a$ is of
bounded type, this post-critical set is dense in the boundary of the
Siegel disk of $f$.\label{fig_ui}}
\end{figure}

Figure \ref{fig_ui} shows the open set $U_1(\a)$ for an $\a$ of
bounded type.

\begin{proposition}
For all $\a\in \irrat_{\geq N}$ and all $j\geq 1$, the post-critical
set $\PC(P_\a)$ is contained in $U_j(\a)$.
\end{proposition}

\begin{proof}
We will show that for all $j\geq 1$, there is a point $z_0\in \C_j$
which is a precritical point of $P_\a$, and a sequence of positive
integers with $t_1<t_2<t_2<\ldots$ such that
\begin{itemize}
\item for all $m\geq 1$, $t_{m+1}-t_m< q_{j+1} +
\bigl(k_1+\lfloor R_3\rfloor -4)q_j$ and

\item $P_\a^{\circ t_m}(z_0) \in C_j$.
\end{itemize}
The proof follows immediately.

Denote by $\omega_{j+1}$ the critical point of $f_{j+1}$. According
to proposition \ref{prop_pcinfinite} the orbit of $\omega_{j+1}$
under iteration of $f_{j+1}$ is infinite. In particular, for all
$m\geq 0$, $f_{j+1}^{\circ m}(\omega_{j+1})$ is in the domain
$U_{f_{j+1}}$ of $f_{j+1}$. Remember that the map $\phi_j\eqdef
\Exp\circ \Phi_{f_j}:D_j\to U_{f_{j+1}}$ is surjective. So, for all
$m\geq 0$, we can find a point $w_m\in D_j$ such that
\[\phi_j(w_m) = f_{j+1}^{\circ m}(\omega_{j+1}).\]
Set
\[z_m\eqdef \Psi_j(w_m)\in C_j.\]
Then, $z_0$ is a precritical point of $P_\a$ and according to lemma
\ref{lemma_correspondanceorbit}, there is an increasing sequence
$(t_m)$ such that $z_m=P_\a^{\circ t_m}(z_0)$. It is therefore
enough to show that for all $m \geq 1$, $t_{m+1}-t_m< q_{j+1} +
\bigl(k_1+\lfloor R_3\rfloor -4)q_j$.

Note that for $m\geq 0$, $w_m\in D_j$, $w'_m\eqdef f_j^{\circ
k_1}(w_m) \in D'_j$. By definition of the renormalization $f_{j+1}$,
we have
\[\phi_j(w'_m) = f_{j+1}\bigl(\phi_j(w_m)\bigr) = f_{j+1}^{\circ (m+1)}(\omega_{j+1}) =
\phi_j(w_m).\] Thus, $\Phi_{f_j}(w_{m+1})-\Phi_{f_j}(w'_m)$ is a
positive integer $\ell_m$. Then,
\[w_{m+1} = f_j^{\circ \ell_m)}(w'_m).\]
We have
\[\re\bigl(\Phi_{f_j}(w'_m)\bigr) \geq 0\quad\text{and}\quad
\re\bigl(\Phi_{f_j}(w_{m+1})\bigr) < \frac{1}{\a_j}-R_3-5.\]
Remember ${\rm a}_{j+1}=\lfloor 1/\a_j\rfloor$. Thus,
\[\ell_m\leq
{\rm a}_{j+1}-\lfloor R_3 \rfloor -4.\]

Set $z'_m\eqdef \Psi_j(w'_m)$. According to propositions
\ref{prop_correspondanceorbit1} and \ref{prop_correspondanceorbit2},
we have
\[z'_m = P_\a^{\circ (k_1 q_j+q_{j-1})}(z_m)\quad\text{and}\quad
z_{m+1} = P_\a^{\circ \ell_m q_j}(z'_m).\] Thus,
\[t_{m+1}-t_m = k_1 q_j+q_{j-1} + \ell_m q_j\leq ({\rm a}_{j+1}+k_1+\lfloor R_3\rfloor -4)q_j+q_{j-1}.\]
The result now follows immediately from $q_{j+1} = {\rm a}_{j+1}q_j
+ q_{j-1}$.
\end{proof}

We will now assume that $\a\in {\cal S}_N$, i.e. $\a\in \irrat_{\geq
N}$ is a bounded type irrational number (the coefficients of the
continued fraction are bounded). We will use the additional
hypothesis that $\a$ has bounded type in order to obtain the
following result.

\begin{proposition}\label{prop_uindeltaeps}
For all $\a\in {\cal S}_N$, for all $\eps>0$, if $j$ is large
enough, the set $U_j(\a)$ is contained in the $\eps$-neighborhood of
the Siegel disk $\Delta_\a$.
\end{proposition}

\begin{proof}
Consider the renormalization tower associated to $f_0\eqdef P_\a$
and let us keep the notations we have introduced so far. Set
\[D''_j\eqdef
f_j^{\circ (a_{j+1}+\ell)}(D_j).\] Define
\[N_j\eqdef {\rm a}_{j+1} -\lfloor R_3\rfloor -1
<\frac{1}{\a_j}-R_3.\]
 Note that
\[D''_j=
\bigl\{z\in \C~;~N_j-3<\re\bigl(\Phi_{f_j}(z)\bigr)< N_j-1\text{ and
}\im(w)>-2\bigr\}.\] In particular, $D''_j\subset \Pet_{f_j}$. Set
\[C''_j\eqdef \Psi_j (D''_j).\]
According to propositions~\ref{prop_correspondanceorbit1} and
\ref{prop_correspondanceorbit2},
\[C''_j = P_\a^{\circ (q_{j+1}+\ell
q_j)}(C_j).\]

\begin{lemma}
There exists $M$ such that for all $j\geq 1$, the disk
$D\bigl(0,|v|e^{-2\pi M}\bigr)$ is
contained in the Siegel disk of $f_j$.
\end{lemma}

\begin{proof}
Let $B(\a_j)$ be the Brjuno sum defined by Yoccoz as
\[B(\a_j)\eqdef \sum_{k=0}^{+\infty} \a_j\cdots\a_{j+k-1}\log
\frac{1}{\a_{j+k}}.\]
Since $\a$ is of bounded type, there is a constant $B$ such that for
all $j\geq 1$, $B(\a_j)\leq B$.

The map $f_j$ has a univalent inverse branch
$g_j:D\bigl(0,|v|\bigr)\to \C$ fixing $0$ with derivative $e^{-2i\pi \a_j}$.
According to a theorem of Yoccoz \cite{y}, there is a constant $C$,
which does not depend on $j$, such that the Siegel disk of $g_j$
contains the disk centered at $0$ with radius
\[|v|e^{-2\pi (B(\a_j)+C)}\geq |v|e^{-2\pi(B+C)}.\]
The lemma is proved with $M\eqdef B+C$.
\end{proof}

Let us now show that for any $\eps>0$, for $j$ large enough, $C''_j$
is contained in the $\eps$-neighborhood of $\Delta_\a$. Denote by
${D''_j}^\sharp$ the set of points in $D''_j$ which are mapped by
$\phi_j= \Exp\circ \Phi_{f_j}$ in  $D\bigl(0,|v|e^{-2\pi
  M}\bigr)$ and set ${D''_j}^\flat\eqdef D''_j\setminus
{D''_j}^\sharp$. In addition, set
\[{C''_j}^\sharp\eqdef
\Psi_j\bigl({D''_j}^\sharp\bigr)\quad\text{and}\quad
{C''_j}^\flat\eqdef \Psi_j\bigl({D''_j}^\flat\bigr).\] Points in
$D\bigl(0,|v|e^{-2\pi M}\bigr)$ have an infinite orbit under
iteration of $f_{j+1}$. It follows that points in ${D''_j}^\sharp$
have an infinite orbit under iteration of $f_j$. Thus, the orbit of
points in ${C''_j}^\sharp$ remains in $U_j(\a)$, thus is bounded. As
a consequence, ${C''_j}^\sharp$ (which is open) is contained in the
Fatou set of $P_\a$, and since it contains $0$ in its boundary,
${C''_j}^\sharp$ is contained in the Siegel disk of $P_\a$.

So, in order to show that  $C''_j$ is contained in the
$\eps$-neighborhood of $\Delta_\a$, it is enough to show that
${C''_j}^\flat$ is contained in the $\eps$-neighborhood of
$\Delta_\a$. Note that ${D''_j}^\flat$ is the image of the rectangle
\[\bigl\{w\in \C~;~N_j-3<\re(w)<N_j-1\text{ and }-2<\im(w)\leq M\bigr\}\]
by the map $\Phi_{f_j}^{-1}$ which is univalent on the strip
\[\bigl\{w\in \C~;~ 0<\re(w)<1/\a_j-R_3\bigr\}.\]
Since
\[1<N_j-3<N_j < 1/\a_j-R_3,\] the modulus of the annulus $\Pet_{f_j}\setminus
\overline{{D''_j}^\flat}$ is bounded from below independently of
$j$.

It follows from Koebe's distortion lemma that there is a constant
$K$ such that \[{\rm diam}({C''_j}^\flat)\leq K\cdot d(z_j,z'_j)\]
where \[z_j\eqdef \Psi_j\circ
\Phi_{f_j}^{-1}(N_j-3)\quad\text{and}\quad z'_j\eqdef \Psi_j\circ
\Phi_{f_j}^{-1}(N_j-2).\] According to proposition
\ref{prop_correspondanceorbit1},
\[z_j =P_\a^{\circ (N_j-3)q_j}(\omega_\a)\quad \text{and}\quad
z'_j=P_\a^{\circ q_j}(z_j).\] The boundary of $P_\a$ is a Jordan
curve, and $P_\a:\partial \Delta_\a\to \partial \Delta_\a$ is
conjugate to the rotation of angle $\a$ on $\R/\Z$. It follows that
\[{\rm diam}({C''_j}^\flat)\leq K\cdot\max_{z\in \partial \Delta_\a} \bigl|P_\a^{\circ
q_j}(z)-z\bigr|.\]

Without loss of generality, we may assume that $M\geq 2$. If $z\in
U_j(\a)$, then there is a $k\leq q_{j+1}+\ell q_j$ such that
$P_\a^{\circ k}(z)\in C''_j$. Then,
\begin{itemize}
\item either $P_\a^{\circ k}(z) \in {C''_j}^\sharp$ in which case
$z\in \Delta_\a$,

\item or $P_\a^{\circ k}(z) \in {C''_j}^\flat$ in which case $z$
belongs to the connected component $O_j^{-k}$ of
$P_\a^{-k}({C''_j}^\flat)$ intersecting $\Delta_\a$.
\end{itemize}
In the second case, $O_j^{-k}$ contains two points $z_j^{-k}$ and
${z_j'}^{-k}$ which are in the boundary of $\Delta_\a$ and which are
respectively mapped to $z_j$ and $z'_j$ by $P_\a^k$. We have
${z'_j}^{-k} = P_\a^{\circ q_j}(z_j^{-k})$.

Note that since $\a$ is of bounded type, there is a constant $A$
such that \[\forall j\geq 1\qquad q_{j+1}+\ell q_j \leq  A\cdot
q_j.\] According to lemma \ref{lemma_pullbacksiegel} below, there is
a constant $K'$ such that for all $j\geq 1$ and all $k\leq
q_{j+1}+\ell q_j$
\[{\rm diam}(O_j^{-k})\leq K'\cdot \bigl|{z'_j}^{-k}-z_j^{-k}\bigr|
\leq K'\cdot\max_{z\in \partial \Delta_\a} \bigl|P_\a^{\circ
q_j}(z)-z\bigr|.\] So, we see that
\[\sup_{z\in U_j(\a)} d(z,\Delta_\a) \leq \max(K,K')\cdot \max_{z\in \partial \Delta_\a} \bigl|P_\a^{\circ
q_j}(z)-z\bigr|\underset{j\to +\infty}\longrightarrow 0.\] This
completes the proof of proposition \ref{prop_uindeltaeps}.
\end{proof}

Assume $\a\in \R\setminus \Q$ is of bounded type. If $z\in
\partial \Delta_\a$, we set
\[r_j(z) = \bigl|P_\a^{\circ q_j}(z)-z\bigr|.\]

\begin{lemma}\label{lemma_pullbacksiegel}
For all $\a\in \R\setminus \Q$ of bounded type, all $A\geq 1$ and
all $K\geq 1$, there exists a $K'$ such that the following holds. If
$j\geq 1$, if $k\leq A\cdot q_j$, if $z_0\in \partial \Delta_\a$, if
$z_k=P_\a^{\circ k}(z_0)$ and if $O$ is the connected component of
$P_\a^{-k}\bigl(D(z_k,K\cdot r_j(z_k))\bigr)$ containing $z_0$, then
\[{\rm diam}(O)\leq K'\cdot r_j(z_0).\]
\end{lemma}

\begin{proof}
The constants $M_1$, $M_2$ and $m$ which will be introduced in the
proof depend on $\alpha$, $A$ and $K$, but they do not depend on
$j$, $k$ or $z$.

Set
\[D\eqdef D\bigl(z_k,K\cdot r_j(z_k)\bigr)\quad\text{and}\quad
\widehat D\eqdef D\bigl(z_k,2K\cdot r_j(z_k)\bigr).\] Since
$\partial \Delta_\a$ is a quasicircle and since $P_\a:\partial
\Delta_\a\to \partial \Delta_\a$ is conjugate to the rotation of
angle $\a$ on $\R/\Z$, the number of critical values of $P_\a^{\circ
k}$ in $\widehat D$ is bounded by a constant $M_1$ which only
depends on $\a$, $A$ and $K$.

Let $O$ (respectively $\widehat O$) be the connected component of
$P_\a^{-k}(D)$ (respectively $P_\a^{-k}(\widehat D)$) containing
$z_0$. The degree of $P_\a^{\circ k}:\widehat O\to \widehat D$ is
bounded by $2^{M_1}$.

On the one hand, it easily follows from the Grötzsch inequality that
the modulus of the annulus $\widehat O\setminus \overline O$ is
bounded from below by $\log 2/(2\pi 2^M)$ (see for example \cite{st}
lemma 2.1).

On the other hand, it follows from Schwarz's lemma that the
hyperbolic distance in $\widehat O$ between $z_0$ and $P_\a^{\circ
q_j}(z_0)$ is greater than the hyperbolic distance in $\widehat D$
between $z_k$ and $P_\a^{\circ q_j}(z_j)$, i.e. a constant $m$ which
only depends on $\a$, $A$ and $K$.

Lemma \ref{lemma_pullbacksiegel} now follows easily from the Koebe
distortion lemma.
\end{proof}

Note that for each fixed $j$, the set $U_j(\a)$ depends continuously
on $\a$ as long as the first $j+1$ approximants remain unchanged.
Hence,  given $\a\in {\cal S}_N$ and $\delta>0$, if $\a'\in
\irrat_{\geq N}$ is sufficiently close to $\a$ (in particular, the
first $j$ entries of the continued fractions of $\a$ and $\a'$
coincide), then $\overline U_j(\a')$ is contained in the
$\delta$-neighborhood of $\overline U_j(\a)$. This completes the
proof of proposition \ref{theo_postcrit}.

\subsection{Lebesgue density near the boundary of a Siegel
disk\label{sec_mcmullen}}

\begin{definition}
If $\a$ is a Brjuno number and if $\delta>0$, we denote by $\Delta$
the Siegel disk of $P_\a$ and by $K(\delta)$ the set of points whose
orbit under iteration of $P_\a$ remains at distance less than
$\delta$ of $\Delta$.
\end{definition}

Our proof will be based on the following theorem of Curtis T.
McMullen \cite{mcm}.

\begin{theorem}[McMullen]\label{theo_mcm}
Assume $\a$ is a bounded type irrational and $\delta>0$. Then, every
point $z\in \partial \Delta$ is a Lebesgue density point of
$K(\delta)$.
\end{theorem}

\begin{figure}[htbp]
\centerline{\scalebox{.8}{\includegraphics{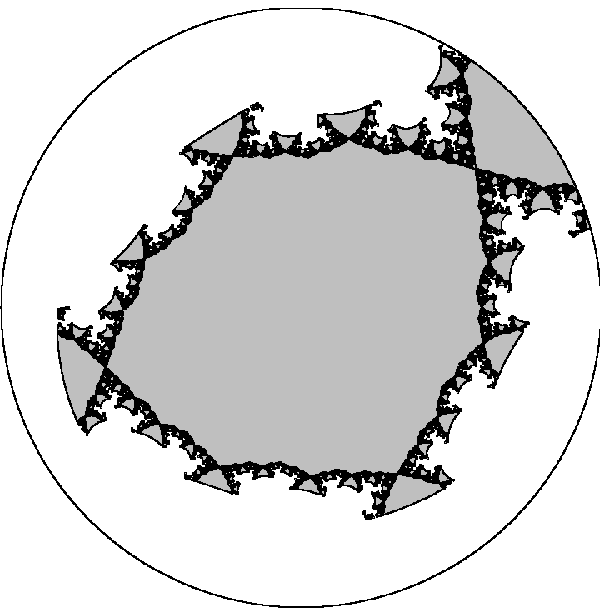}}\qquad
\scalebox{.8}{\includegraphics{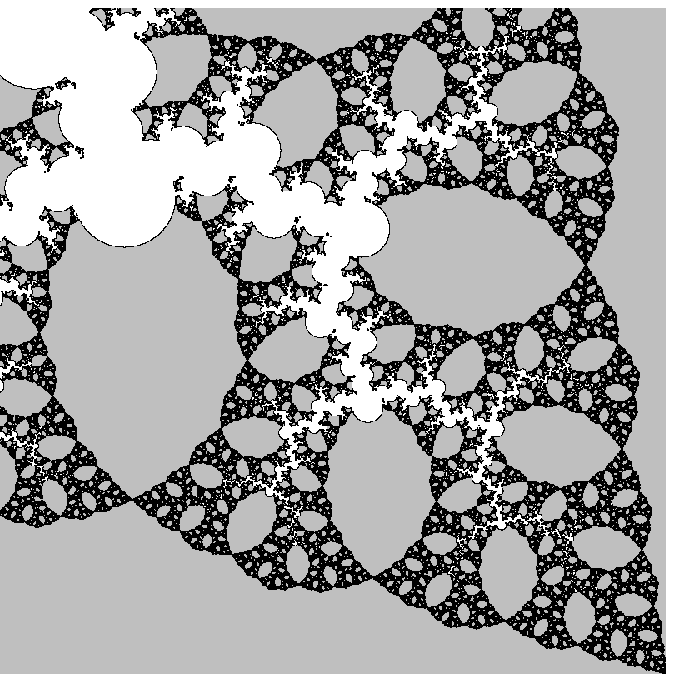}} } \caption{If
$\a=(\sqrt5-1)/2$, the critical point of $P_\a$ is a Lebesgue
density point of the set of points whose orbit remain in $D(0,1)$.
Left: the set of points whose orbit remains in $D(0,1)$. Right: a
zoom near the critical point.}
\end{figure}

\begin{corollary}\label{coro_mcm}
Assume $\a$ is a bounded type irrational and $\delta>0$. Then
\[d\eqdef d(z,\partial \Delta)\to 0~\text{ with }z\not\in \overline
\Delta \quad\Longrightarrow\quad  \dens_{D(z,d)}\bigl(\C\setminus
K(\delta)\bigr) \to 0.\]
\end{corollary}

\proof We proceed by contradiction. Assume we can find a sequence
$(z_j)$ such that
\begin{itemize}
\item $d_j\eqdef d(z_j,\partial \Delta)\to 0$ and

\item $\rho_j\eqdef \dens_{D(z_j,d_j)}\bigl(\C\setminus K(\delta)\bigr)\not\to
  0$.
\end{itemize}
Extracting a subsequence if necessary, we may assume that the
sequence $(z_j)$ converges to a point $z_0\in \partial \Delta$ and
that $\lim \rho_j= \rho>0$.

Set $\eta\eqdef \rho/5$ and for $i\geq 1$, set
\[X_i\eqdef \bigl\{w\in \partial \Delta\Tq (\forall r\leq 1/i)~ \dens_{D(w,r)}
\bigl( \C\setminus K(\delta)\bigr)\leq \eta\bigr\}.\] The sets $X_i$
are closed. By McMullen's theorem \ref{theo_mcm}, $\bigcup
X_i=\partial \Delta$. By Baire category, one of these sets $X_i$
contains an open subset $W$ of $\partial \Delta$. Then, for all
sequence of points $w_j\in W$ and all sequence of real number $r_j$
converging to $0$, we have
\begin{equation}\label{eq_Y}
\limsup_{j\to +\infty }\dens_{D(w_j,r_j)}\bigl(\C\setminus
K(\delta)\bigr)\leq \eta=\frac{\rho}{5}.
\end{equation}

We claim that there is a
map $g$ defined and univalent in a neighborhood $U$ of $z_0$, such
that
\begin{itemize}
\item $g(z_0)=w_0\in W$,

\item $g\bigl(K(\delta)\cap U\bigr) = K(\delta)\cap g(U)$ and

\item $g(\partial \Delta\cap U)= \partial \Delta\cap g(U)$.
\end{itemize}
Indeed, if $z_0$ is not precritical, we can find an integer $k\geq
0$ such that $f^{\circ k}(z_0)\in W$ and we let $g$ be the
restriction of $f^{\circ k}$ to a sufficiently small neighborhood of
$z_0$. If $z_0$ is precritical, we can find a point $w_0\in W$ and
an integer $k\geq 0$ such that $f^{\circ k}(w_0)=z_0$ and we let $g$
coincide the restriction of the branch of $f^{-k}$ sending $z_0$ to
$w_0$, to a sufficiently small neighborhood of $z_0$.

Let $z'_j\in \partial \Delta$ be such that $|z_j-z'_j|=d_j$. Then,
$z'_j\underset{j\to +\infty}\longrightarrow z_0$. Let $j$ be
sufficiently large so that $z'_j\in U$ and set $w_j:=g(z'_j)$. On
the one hand, $w_j\underset{j\to +\infty}\longrightarrow w_0$. Thus,
$w_j\in W$ for $j$ large enough. On the other hand,
\[\dens_{D(z'_j,2d_j)}
 \bigl(\C\setminus K(\delta)\bigr)\geq \frac{1}{4}\dens_{D(z_j,d_j)}
 \bigl(\C\setminus K(\delta)\bigr)\]
and so
\[\liminf_{j\to +\infty} \dens_{D(z'_j,2d_j)}
 \bigl(\C\setminus K(\delta)\bigr)\geq \frac{\rho}{4}.\]
Since $g$ is holomorphic at $z_0$,
\[\liminf_{j\to +\infty} \dens_{D(w_j,r_j)}
 \bigl(\C\setminus K(\delta)\bigr)\geq \frac{\rho}{4}\quad
 \text{with}\quad
r_j\eqdef \bigl|g'(w_0)\bigr|\cdot 2d_j\underset{j\to +\infty}
\longrightarrow 0.\] This contradicts (\ref{eq_Y}).
 \qed\par\medskip

\subsection{The proof\label{sec_theproof}}

We will now prove proposition \ref{prop_keycremer}. We let $N$ be
sufficiently large so that the conclusions of proposition
\ref{theo_postcrit} and corollary \ref{cor_postcrit} apply. Assume
$\alpha\in {\cal S}_N$ and choose a sequence $(A_n)$ such that
\[\sqrt[q_n]{A_n}\To_{n\to +\infty}+\infty\quad \text{and}\quad
\sqrt[q_n]{\log A_n}\To_{n\to +\infty}1.\] Set
\[\a_n\eqdef [{\rm a}_0,{\rm a}_1,\ldots, {\rm a}_n, A_n,
N,N,N,\ldots].\] Note that since $\a$ is of bounded type, the Julia
set $J_\alpha$ has zero Lebesgue measure (see \cite{pe}).
Proposition \ref{theo_disquesdigites} then easily implies that
\[\liminf \area(K_{\a_n}) \geq \frac{1}{2} \area(K_\a).\]
Everything relies on our ability to promote the coefficient $1/2$ to
a coefficient $1$.

Denote by $K$ (resp. $K_n$) the filled-in Julia set of $P_\a$ (resp.
$P_{\a_n}$) and by $\Delta$ (resp. $\Delta_n$) its Siegel disk. For
$\delta>0$, set
\begin{eqnarray*}
V(\delta) &\eqdef & \bigl\{ z\in \C \Tq d(z,\partial
\Delta)< \delta\bigr\},\\
K(\delta)&\eqdef& \bigl\{ z\in V(\delta) \Tq (\forall k\geq
0)~P_\a^{\circ k}(z)\in
V(\delta)\bigr\} \quad\text{and}\quad \\
K_n(\delta)&\eqdef& \bigl\{ z\in V(\delta) \Tq (\forall k\geq
0)~P_{\a_n}^{\circ k}(z)\in V(\delta)\bigr\}.
\end{eqnarray*}
Define $\rho_n:\left]0,+\infty\right[\to [0,1]$ by
\[\rho_n(\delta)\eqdef \dens_\Delta\bigl(\C\setminus
K_n(\delta)\bigr).\]

\begin{lemma}\label{lemma_densities}
For all $\delta>0$, there exist $\delta'>0$ (with $\delta'<\delta$)
and a sequence $(c_n>0)$ converging to $0$, such that
\[\rho_n(\delta)\leq \frac{3}{4} \rho_n(\delta')+c_n.\footnote{The
    coefficient $\frac{3}{4}$ could have been replaced by any
    $\lambda>\frac{1}{2}$}\]
\end{lemma}

This lemma enables us to complete the proof of proposition
\ref{prop_keycremer} as follows. We set
\[\rho(\delta)\eqdef \limsup_{n\to +\infty} \rho_n(\delta)\quad (\leq 1).\]
Then, for all $\delta>0$, there is a $\delta'>0$ such that
$\rho(\delta)\leq \frac{3}{4} \rho(\delta')$. Since $\rho$ is
bounded from above by $1$, this implies that $\rho$ identically
vanishes. In other words
\begin{equation}\label{eq_density1Kn}(\forall \delta>0)\quad
\dens_\Delta\bigl(K_n(\delta)\bigr)\underset{n\to
  +\infty}\longrightarrow 1.
\end{equation}
Since $K_n(\delta)\subset K_{\a_n}$, we deduce that $\ds
\dens_\Delta(K_{\a_n})\underset{n\to
  +\infty}\longrightarrow 1$.
We know that
\begin{itemize}
\item $P_{\a_n}$ converges
locally uniformly to $P_\a$,

\item the orbit under iteration of $P_\a$ of any point in $K_\alpha\setminus J_\alpha$
eventually lands in $\Delta$ and

\item $P_{\a_n}^{-1}(K_{\a_n})=K_{\a_n}$.
\end{itemize}
It follows that $\ds \dens_{K_\a\setminus
J_\a}(K_{\a_n})\underset{n\to +\infty}\longrightarrow 1$. Since the
Julia set $J_\a$ has Lebesgue measure zero, this implies that
$\liminf \area (K_{\a_n})\geq \area(K_\a)$. This completes the proof
of proposition \ref{prop_keycremer} up to Lemma
\ref{lemma_densities}.

\medskip
\noindent{\sc Proof of Lemma \ref{lemma_densities}.}
Let us sum up what we obtained in sections \ref{sec_pertsigeldisk},
\ref{sec_inoushishikura} and \ref{sec_mcmullen}.
\begin{enumerate}[(A)]
\item \label{prop_dempos1}
For all open set $U\subset \Delta$ and all $\delta>0$, $\ds
\liminf_{n\to +\infty} \dens_U\bigl(K_n(\delta)\bigr) \geq
\frac{1}{2}.$ This is an immediate consequence of proposition
\ref{theo_disquesdigites} in section \ref{sec_pertsigeldisk}.

\item \label{prop_dempos2} For all $\delta>0$, if $n$ is sufficiently
large, the post-critical set of $P_{\a_n}$ is contained in
$V(\delta)$. This is just a restatement of corollary
\ref{cor_postcrit} in section \ref{sec_inoushishikura}.

\item \label{prop_dempos3} For all $\eta>0$ and all $\delta>0$, there
exists $\delta'_0>0$ such that if $\delta'<\delta'_0$ and if $z\in
\overline {V(8\delta')}\setminus V(2\delta')$, then $\ds
\dens_{D(z,\delta')}\bigl(\C\setminus K(\delta)\bigr)<\eta.$ This is
an easy consequence of corollary \ref{coro_mcm} in section
\ref{sec_mcmullen}.
\end{enumerate}

\newcounter{comptestep}
\begin{list}{\textbf{Step \arabic{comptestep}.}}%
{\usecounter{comptestep}%
 \setlength{\leftmargin}{0cm}%
 \setlength{\rightmargin}{0cm}%
 \setlength{\itemindent}{0.2cm}%
 \setlength{\labelsep}{0.2cm}%
}

\item\label{step:kappa} By Koebe distortion theorem, there exists a
constant $\kappa$ such that for all map $\phi:D\eqdef D(a,r)\to \C$
which extends univalently to $D(a,3r/2)$, we have
\[\sup_D|\phi'| \leq \kappa \inf_D |\phi'|.\]
We choose $\eta>0$ such that
\[8\pi \kappa^2 \eta < \frac{1}{4}.\]

\item\label{step:deltaprime} Fix $\delta>0$. We claim that there exists
$\delta'>0$ such that:
\begin{enumerate}[(i)]
\item\label{step2.1} $9\delta'< \delta$ and $(2+3\kappa)\cdot
\delta'<\delta$,\footnote{Those requirements will be used in step
\ref{step:phiKndelta}.}

\item\label{step2.2} if $d(z,\Delta)< 2\delta'$, then
$d\bigl(P_\alpha(z),\Delta\bigr)<8\delta'$ and

\item\label{step2.3} if $z\in \overline {V(8\delta')}\setminus
 V(2\delta')$, then
$\dens_{D(z,\delta')}\bigl(\C\setminus K(\delta)\bigr)<\eta$.
\end{enumerate}
Indeed, it is well-known that for $\a\in \R$, $\bigl|P_\a'\bigr|<4$
on $K_\alpha$. As a consequence, if $\delta'>0$ is sufficiently
small, then $\bigl|P_\a'\bigr|<4$ on $V(2\delta')$. It follows that
(\ref{step2.2}) holds for $\delta'>0$ sufficiently small. Claim
(\ref{step2.3}) follows from the aforementioned point
(\ref{prop_dempos3}).

From now on, we assume that $\delta'$ is chosen so that the above
claims hold and we set
\[W\eqdef \overline {V(8\delta')}\setminus V(2\delta').\]

\item\label{step:Yl} Set
\[Y^\ell\eqdef \bigl\{z\in K(\delta)\Tq
P_\a^{\circ \ell}(z)\in
 \Delta\bigr\}.\]
The set of points in $K(\delta)$ whose orbits do not intersect
$\Delta$, is contained in the Julia set of $P_\alpha$. This set has
zero Lebesgue measure. Thus, $K(\delta)$ and $\bigcup Y^\ell$
coincide up to a set of zero Lebesgue measure. The sequence
$(Y^\ell)_{\ell\geq 0}$ is increasing. From now on, we assume that
$\ell$ is sufficiently large so that
\[\bigl(\forall w\in W\bigr)\quad
\dens_{D(w,\delta')}(\C\setminus Y^\ell)<\eta.\]

\item\label{step:phiYl} Assume $\phi$ is univalent on $D(w,3\delta'/2)$ with $w\in W$,
$r$ is the radius of the largest disk centered at $\phi(w)$ and
contained in $\phi\bigl(D(w,\delta')\bigr)$ and $Q$ is a square
contained in $\phi\bigl(D(w,\delta')\bigr)$ with side length at
least $r/\sqrt 8$. Set $D\eqdef D(w,\delta')$. Then, $r\geq \inf_D
|\phi'|\cdot \delta'$ and thus,
\[\area(Q)\geq \inf_D|\phi'|^2 \cdot
\frac{(\delta')^2}{8}.\] In addition, $\ds \sup_D|\phi'| \leq \kappa
\inf_D |\phi'|$ and so,
\[\dens_{Q} \bigl(\C\setminus \phi(Y^{\ell})\bigr) \leq
\frac{\area\bigl(\phi(D\setminus Y^\ell)\bigr)}{\area(Q)}\leq
\frac{\ds \sup_D|\phi'|^2 \cdot \pi(\delta')^2\cdot \eta}{\ds
\inf_{D}|\phi'|^2\cdot (\delta')^2/8} \leq 8\pi
\kappa^2\eta<\frac{1}{4}.\] As a consequence,
\[\dens_Q\bigl( \phi(Y^\ell)\bigr)> \frac{3}{4}.\]

\item\label{step:weaklimit} If $X\subset \C$ is a measurable set, we use the notation
$m|_X$ for the Lebesgue measure on $X$, extended by $0$ outside $X$.
If $U\subset \C$ is an open set, $(X_n)$ is a sequence of measurable
subsets of $\C$ and $\lambda\in [0,1]$, we say that
\[\liminf_{n\to +\infty} m|_{X_n} \geq \lambda \cdot m|_U\]
if for all non empty open set $U'$ relatively compact in $U$, we
have
\[\liminf_{n\to +\infty} \dens_{U'}(X_n)\geq \lambda.\footnote{Equivalently,
for all non empty open set $U'\subset \C$ with finite area, $\ds
\liminf_{n\to +\infty} \dens_{U'}(X_n)\geq \lambda\cdot
\dens_{U'}(U)$.}\] Assume $f:V\to U$ is a holomorphic map, nowhere
locally constant, and $(f_n:V_n\to \C)$ is a sequence of holomorphic
maps such that
\begin{itemize}
\item every compact subset of $V$ is eventually contained in $V_n$
and

\item the sequence $(f_n)$ converges uniformly to $f$ on every
compact subset of $V$.
\end{itemize}
Then,
\[\liminf_{n\to +\infty} m|_{X_n}\geq \lambda \cdot
m|_U\quad\Longrightarrow \quad \liminf_{n\to +\infty}
m|_{f_n^{-1}(X_n)}\geq \lambda \cdot m|_V.\]

\item\label{step:weaklimitinW} Set
\[Y_n^\ell\eqdef \bigl\{z\in V(\delta)\Tq (\forall j\leq \ell)~P_{\a_n}^{\circ j}(z)\in
V(\delta)\text{ and }P_{\a_n}^{\circ \ell}(z)\in \Delta\bigr\}.\] On
the one hand, if $z\in Y_n^\ell$ and $P_{\a_n}^{\ell}(z)\in
K_n(\delta)$, then $z\in K_n(\delta)$. On the other hand, every
compact subset of $Y^\ell$ is eventually contained in $Y_n^\ell$ and
the sequence $(P_{\a_n}^{\circ \ell})$ converges uniformly to
$P_\a^\ell$ on every compact subset of $Y^\ell$. By the
aforementioned point (\ref{prop_dempos1}), we have
\[\liminf_{n\to +\infty} m|_{K_n(\delta)} \geq
\frac{1}{2}m|_{\Delta}.\] So, according to step
\ref{step:weaklimit},
\[\liminf_{n\to +\infty} m|_{K_n(\delta)} \geq
\frac{1}{2}m|_{Y^\ell}.\]

\item\label{step:weaklimitinQn} Assume $\phi_n$ is univalent on $D(w_n,3\delta'/2)$ with
$w_n\in W$, $r_n$ is the radius of the largest disk centered at
$\phi_n(w_n)$ and contained in $\phi_n\bigl(D(w_n,\delta')\bigr)$
and $Q_n$ is a square contained in
$\phi_n\bigl(D(w_n,\delta')\bigr)$ with side length at least
$r_n/\sqrt 8$. Then,
\[\liminf_{n\to +\infty}
\dens_{Q_n}\bigl( \phi_n\bigl(K_n(\delta)\bigr)\bigr)\geq
\frac{3}{8}.\] Indeed, assume $\lambda$ is a limit value of the
sequence
\[\dens_{Q_n}\bigl( \phi_n\bigl(K_n(\delta)\bigr)\bigr).\]
Post-composing the maps $\phi_n$ with affine maps and extracting a
subsequence if necessary, we may assume that $(w_n)$ converges to
$w\in W$, $(\phi_n)$ converges locally uniformly to
$\phi:D(w,3\delta'/2)\to \C$, $r_n$ converges to the radius $r$ of
the largest disk centered at $\phi(w)$ and contained in
$\phi\bigl(D(w,\delta')\bigr)$ and $Q_n$ converges to a square $Q$
with side length at least $r/\sqrt 8$. According to steps
\ref{step:weaklimit} and \ref{step:weaklimitinW},
\[\liminf_{n\to +\infty} m|_{
\phi_n(K_n(\delta))}\geq \frac{1}{2}m|_{\phi(Y^\ell)}.\] According
to step \ref{step:phiYl}, it follows that
\[\lambda\geq \frac{1}{2} \dens_Q\bigl(\phi(Y^\ell)\bigr)\geq
\frac{3}{8}.\]

\item\label{step:choosingn} From now on, we assume that $n$ is sufficiently large, so that:
\begin{enumerate}[(i)]
\item $\Delta\setminus K_n(\delta)\subset X_n \subset \Delta\setminus
K_n(\delta')$ with
\[X_n\eqdef \bigl\{z\in \Delta \Tq (\exists k)~ P_{\a_n}^{\circ
k}(z)\in W\bigr\}\] (this is possible by step
\ref{step:deltaprime});

\item $s_n< \delta'$ with
\[s_n\eqdef \sup_{z\in \Delta}  d\bigl(z,K_n(\delta')\bigr)\]
(this is possible since $s_n\To_{n\to +\infty}0$ in order for the
aforementioned point (\ref{prop_dempos1}) to hold);

\item the post-critical set of $P_{\a_n}$ is contained
in $V(\delta'/2)$ (this is possible by the aforementioned point
(\ref{prop_dempos2}));

\item\label{stepn.4} if $\phi$ is univalent on
$D(w,3\delta'/2)$ with $w\in W$, if $r$ is the radius of the largest
disk centered at $\phi(w)$ and contained in
$\phi\bigl(D(w,\delta')\bigr)$ and if $Q$ is a square contained in
$\phi\bigl(D(w,\delta')\bigr)$ with side length at least $r/\sqrt
8$, then
\[\dens_Q\bigl( \phi\bigl(K_n(\delta)\bigr)\bigr)\geq \frac{1}{4}\]
(this is easily follows from step \ref{step:weaklimitinQn} by
contradiction).
\end{enumerate}

\item\label{step:phiKndelta} Assume $z_0\in X_n$. Then, we have
\[z_0\in X_n\overset{P_{\a_n}}\mapsto z_1\in V(2\delta')\overset{P_{\a_n}}\mapsto \cdots
\overset{P_{\a_n}}\mapsto z_{k-1}\in
V(2\delta')\overset{P_{\a_n}}\mapsto z_k\in W\] for some integer
$k>0$. Since the post-critical set of $P_{\a_n}$ is contained in
$V(\delta'/2)$, for $j\leq k$ there exists a univalent map
$\phi_j:D\eqdef D(z_k,\delta')\to \C$ such that
\begin{itemize}
\item $\phi_j$ is the inverse branch of $P_{\a_n}^{\circ k-j}$
which maps $z_k$ to $z_j$ and

\item $\phi_j$ extends univalently to $D(z_k,3\delta'/2)$.
\end{itemize}
In particular,
\[\sup_D|\phi_j'| \leq \kappa \inf_D |\phi_j'|.\]
Let $D(z_j,r_j)$ be the largest disk centered at $z_j$ and contained
in $\phi_j(D)$ and $D(z_j,R_j)$ be the smallest disk centered at
$z_j$ and containing $\phi_j(D)$.  Note that $D$ is contained in
$\C\setminus V(\delta')$ and so, for $j\leq k-1$, $D(z_j,r_j)\subset
\phi_j(D)\subset \C\setminus K_n(\delta')$. On the one hand,
$d(z_j,\Delta)<2\delta'$ and on the other hand, every point of
$\Delta$ is at distance at most $s_n$ from a point of
$K_n(\delta')$. It follows that
\[R_j\leq \kappa r_j\leq \kappa\cdot (s_n+2\delta').\]
If $w_0\in \phi_0(D)$ and $w_j\eqdef P_{\a_n}^{\circ j}(w_0)$, then
for $j\leq k-1$,
\[d(w_j,\Delta)\leq d(w_j,z_j)+d(z_j,\Delta) \leq \kappa\cdot
(s_n+2\delta') + 2\delta'<(2+3\kappa)\cdot \delta'<\delta\]
and for $j=k$,
\[d(w_k,\Delta)\leq d(w_k,z_k)+d(z_k,\Delta) \leq 9\delta'<\delta.\]
In other words, $w_0$, $w_1$, \ldots, $w_k$ all belong to
$V(\delta)$. As a consequence,
\[\phi_0\bigl(K_n(\delta)\bigr)\subset  K_n(\delta).\]

\item\label{step:douadic} Continuing with the notations of step \ref{step:phiKndelta},  we denote by
$Q_{z_0}$ the largest douadic square containing $z_0$ and contained
in $D(z_0,r_0)$. On the one hand, since $z_0\in \Delta$ and since
$\phi_0(D)\subset \C\setminus K_n(\delta')$, we have $r_0\leq s_n$,
and so
\[Q_{z_0}\subset D(z_0,r_0)\subset V(s_n)\setminus K_n(\delta').\]
On the other hand, $Q_{z_0}$ has an edge of length greater than
$r_0/2\sqrt 2$ and so, according to step \ref{step:choosingn} point
(\ref{stepn.4}),
\[\dens_{Q_{z_0}} \bigl(K_n(\delta)\bigr) > \frac{1}{4}.\] As a
consequence
\[\dens_{Q_{z_0}} \bigl(\C\setminus K_n(\delta)\bigr) <\frac{3}{4}.\]

Given two douadic squares $Q$ and $Q'$, either $Q\cap Q'=\emptyset$,
or $Q\subset Q'$ or $Q'\subset Q$. It follows that
\begin{eqnarray*} \area\bigl(\Delta\setminus K_n(\delta)\bigr)
&\leq & \frac{3}{4}\area\left(\bigcup_{z\in X_n} Q_z\right)
\\
&\leq& \frac{3}{4} \area\bigl(V(s_n)\setminus K_n(\delta')\bigr)
\\
&\leq& \frac{3}{4} \area\bigl(\Delta\setminus K_n(\delta')\bigr) +
\frac{3}{4}\area\bigl(V(s_n)\setminus \Delta\bigr)\\
&=& \frac{3}{4} \area\bigl(\Delta\setminus K_n(\delta')\bigr) +
c_n\cdot \area(\Delta)
\end{eqnarray*}
with
\[c_n\eqdef\frac{3}{4} \frac{\area\bigl(V(s_n)\setminus \Delta\bigr)}{ \area(\Delta)}.\]

\item\label{step:cn} Since $s_n\to 0$ and since the boundary of $\Delta$ has zero
Lebesgue measure,
\[\area\bigl(V(s_n)\setminus \Delta\bigr)\underset{n\to +\infty}
\longrightarrow 0.\] Thus,
\[\dens_{\Delta}\bigl(\C\setminus K_n(\delta)\bigr) < \frac{3}{4}
\dens_{\Delta}\bigl(\C\setminus K_n(\delta')\bigr) +
c_n\quad\text{with}\quad c_n\underset{n\to +\infty} \longrightarrow
0.\] This completes the proof of Lemma \ref{lemma_densities}.
\qed\par\medskip

\end{list}

\section{The linearizable case}

In order to find a quadratic polynomial with a linearizable fixed
point and a Julia set of positive area, we need to modify the
argument.

\begin{definition}
If $\a$ is a Brjuno number, we denote by $\Delta_\a$ the Siegel disk
of $P_\a$ and by $r_\a$ its conformal radius. For $\rho\leq r_\a$,
we denote by $\Delta_\a(\rho)$ the invariant sub-disk with conformal
radius $\rho$ and by $L_\a(\rho)$ the set of points in $K_\a$ whose
orbits do not intersect $\Delta_\a(\rho)$.
\end{definition}

\begin{figure}[htbp]
\centerline{ \scalebox{.45}{\includegraphics{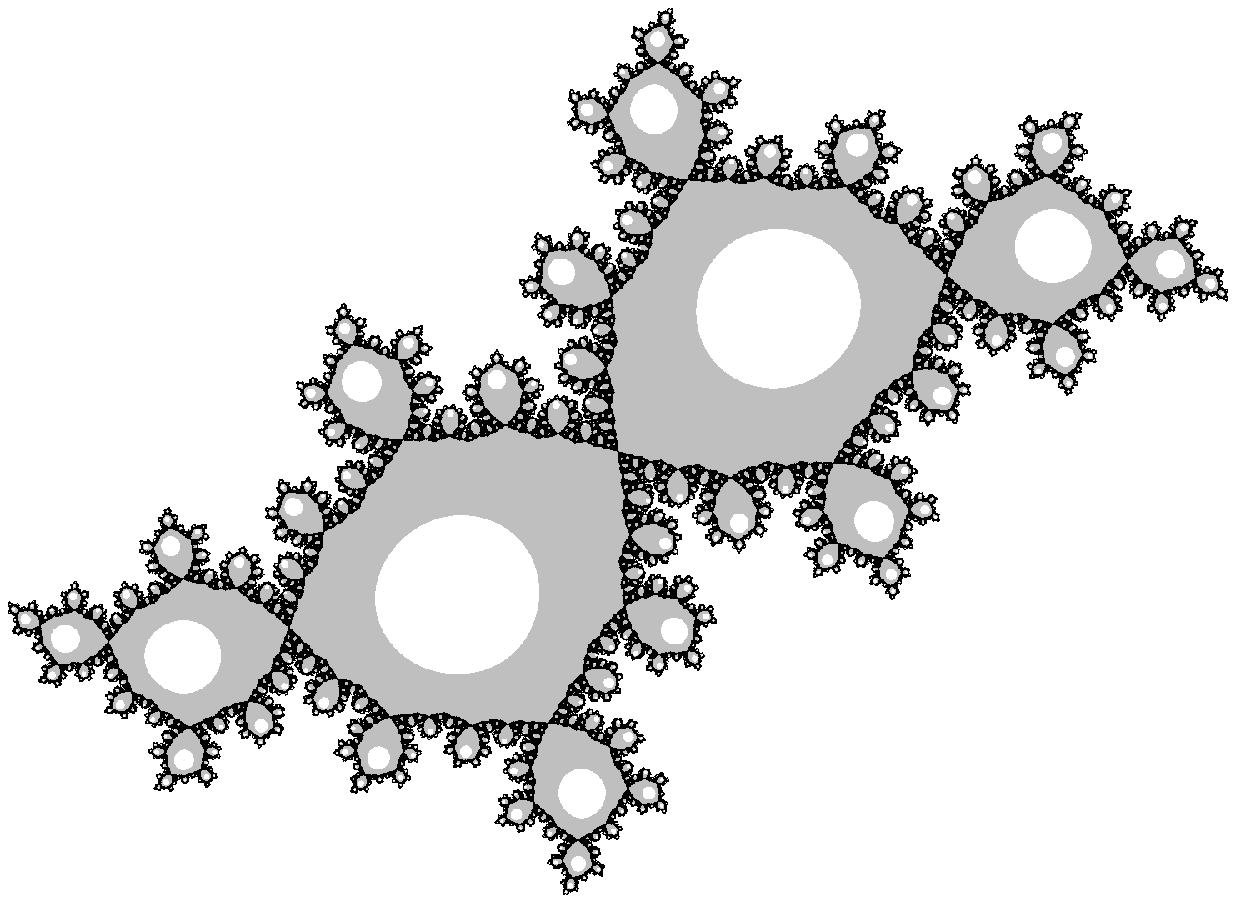}}
\scalebox{.45}{\includegraphics{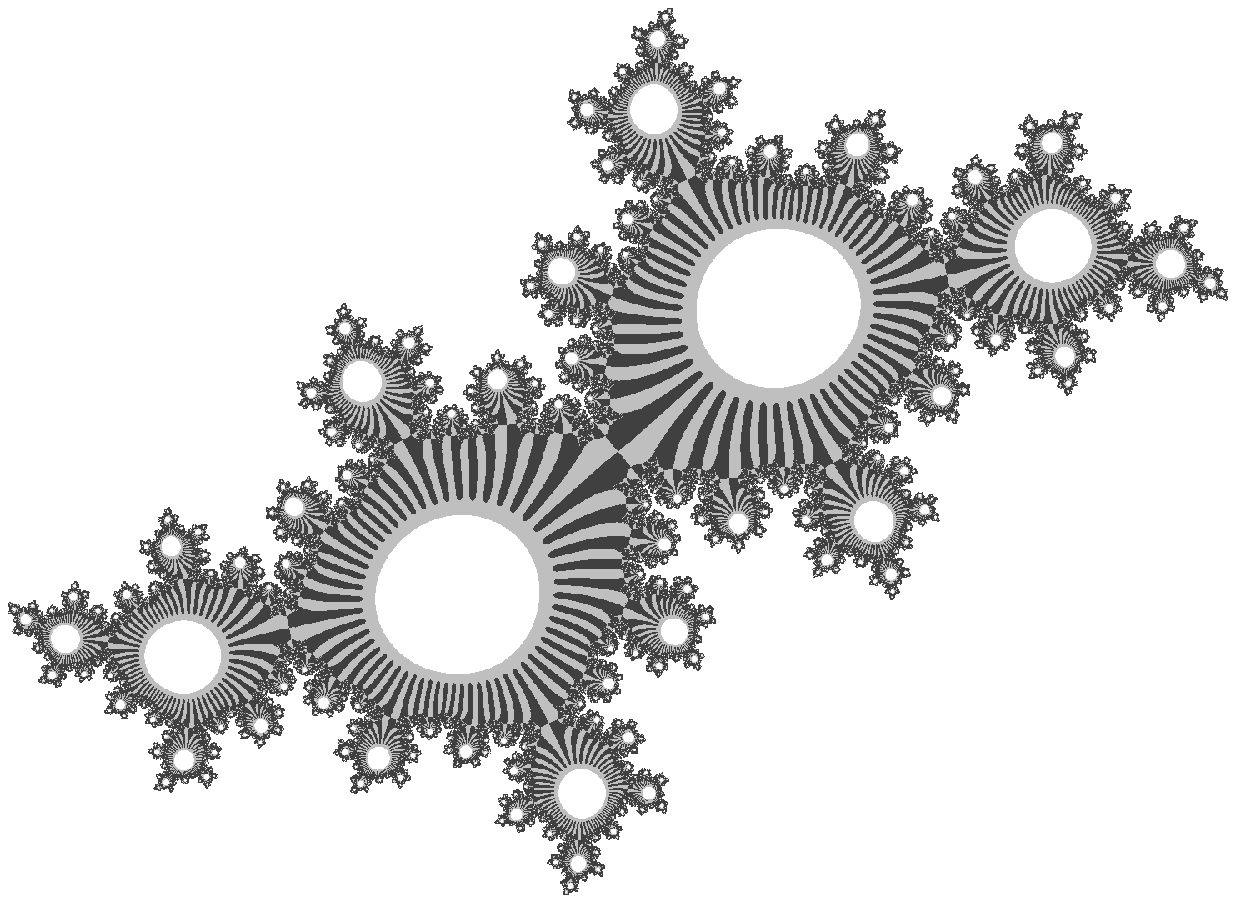}}} \caption{Two sets
$L_{\a}\rho)$ and $L_{\a'}(\rho)$, with $\a'$ a well-chosen
perturbation of $\a$ as in proposition
 \ref{prop_keylin0}. This proposition
asserts that if $\a$ and $\a'$ are chosen carefully enough, the loss
of measure from $L_{\a}(\rho)$ to $L_{\a'}(\rho)$ is small. We
colored white the basin of infinity, the invariant subdisks
$\Delta_\a(\rho)$ and $\Delta_{\a'}(\rho)$ and their preimages; we
colored light grey the remaining parts of the Siegel disks and their
preimages; we colored dark grey the pixels where the preimages are
too small to be drawn.}
\end{figure}

\begin{proposition}\label{prop_keylin0}
There exists a set ${\cal S}$ of bounded type irrationals such that
for all $\a\in {\cal S}$, all $\rho<\rho'<r_\a$ and all $\eps>0$,
there exists $\a'\in {\cal S}$ with
\begin{itemize}
\item $|\a'-\a|<\eps$,

\item $\max\bigl(\rho,(1-\eps)\rho'\bigr)<r_{\a'}<(1+\eps)\rho'$ and

\item $\area\bigl(L_{\a'}(\rho)\bigr)\geq (1-\eps)\area\bigl(L_{\a}(\rho)\bigr)$.
\end{itemize}
\end{proposition}

\proof We let $N$ be sufficiently large so that the conclusions of
proposition \ref{theo_postcrit} and corollary \ref{cor_postcrit}
apply. We will work with ${\cal S}= {\cal S}_N$. Assume $\alpha\in
{\cal S}_N$ and choose a sequence $(A_n)$ such that
\[\lim_{n\to +\infty} \sqrt[q_n]{A_n} =\frac{r_\a}{\rho'}.\]
Set
\[\a_n\eqdef [{\rm a}_0,{\rm a}_1,\ldots, {\rm a}_n, A_n,
N,N,N,\ldots].\] This guaranties that $\ds r_{\a_n}\To_{n\to
+\infty} \rho'$ (see \cite{abc}).

Let $\Delta$ be the Siegel disk of $P_\a$.  Let us use the notations
$V(\delta)$, $K(\delta)$ and $K_n(\delta)$ introduced in section
\ref{sec_theproof}. With an abuse of notations, set
$\Delta(\rho)\eqdef \Delta_\a(\rho)$ and $\Delta_n(\rho)\eqdef
\Delta_{\a_n}(\rho)$. Set
\[\Delta'(\rho)\eqdef
P_\a^{-1}\bigl(\Delta(\rho)\bigr)\setminus \Delta(\rho).\] Then,
$\Delta(\rho)$ and $\Delta'(\rho)$ are symmetric with respect to the
critical point of $P_\a$. The orbit under iteration of $P_\a$ of a
point $z\notin \Delta(\rho)$ lands in $\Delta(\rho)$ if and only if
it passes through $\Delta'(\rho)$. We have a similar property for
\[
\Delta'_n(\rho)\eqdef
P_{\a_n}^{-1}\bigl(\Delta_n(\rho)\bigr)\setminus \Delta_n(\rho).\]

We have proved -- see equation (\ref{eq_density1Kn})  -- that
\[(\forall \delta>0)\quad  \dens_{\Delta}
\bigl(K_n(\delta)\bigr) \underset{n\to +\infty}\longrightarrow 1.\]
The sequence of compact sets $\bigl(\overline \Delta_n(\rho)\bigr)$
converges to $\overline \Delta(\rho)$ for the Hausdorff topology on
compact subsets of $\C$, because $\lim r_{\a_n}>\rho$. It
immediately follows that for all $\delta>0$,
\[\dens_{\Delta\setminus \overline \Delta(\rho)} \bigl(K_n(\delta)\setminus \Delta_n(\rho)\bigr)
\underset{n\to +\infty}\longrightarrow
1.\]

Choose $\delta$ sufficiently small so that $V(\delta)$ does not
intersect $\overline \Delta'(\rho)$. Then, for $n$ large enough
$V(\delta)$ does not intersect $\overline \Delta'_n(\rho)$. In that
case, the orbit under iteration of $P_{\a_n}$ of a point in
$K_n(\delta)\setminus \Delta_n(\rho)$ cannot pass through
$\Delta'_n(\rho)$ and so,
\[K_n(\delta)\setminus \Delta_n(\rho)\subset L_{\a_n}(\rho).\]
Thus,
\[\dens_{\Delta\setminus \overline \Delta(\rho)} \bigl(L_{\a_n}(\rho)\bigr)
\underset{n\to +\infty}\longrightarrow 1.\] The points of
$L_\a(\rho)$ whose orbits do not intersect $\Delta\setminus
\overline \Delta(\rho)$ are contained in the union of the Julia set
$J_\a$ and the countably many preimages of $\partial \Delta(\rho)$.
Thus, they form a set of zero Lebesgue measure. It follows that
\[\area\bigl(L_{\a_n}(\rho)\bigr)\underset{n\to
+\infty}\longrightarrow \area\bigl(L_\a(\rho)\bigr).\]
 \qed\par\medskip

\noindent{\it Proof of theorem \ref{theo_arealin}.}~ We start with
$\a_0\in {\cal S}$ and set $\rho_0\eqdef r_{\a_0}$. We then choose
$\rho\in \left]0,\rho_0\right[$ and two sequences of real numbers
$\e_n$ in $(0,1)$ and $\rho_n$ in $(0,\rho_0)$ such that
$\prod(1-\e_n)>0$ and $\rho_n\searrow \rho>0$. We can construct
inductively a Cauchy sequence $(\a_n\in {\cal S})$ such that for all
$n\geq 1$,
\begin{itemize}
\item $r_{\a_n}\in (\rho_n,\rho_{n-1})$ and

\item
$\area\bigl(L_{\a_{n}}(\rho)\bigr) \geq (1-\e_n)
\area\bigl(L_{\a_{n-1}}(\rho)\bigr)$.
\end{itemize}
Let $\a$ be the limit of the sequence $(\a_n)$.  The conformal
radius of a fixed Siegel disk depends upper semi-continuously on the
polynomial (a limit of linearizations linearizes the limit). So,
$r_\a\geq \lim r_{\a_n}=\rho$. Also, by choosing $\a_n$ sufficiently
close to $\a_{n-1}$ at each step, we can guaranty that $r_\a\leq
\rho$, in which case $r_\a=\rho$.

In addition, the sequence of pointed domains
$\bigl(\Delta_{\a_n}(\rho),0\bigr)$ converges for the Carathéodory
topology to $(\Delta_\a,0)$. In particular, every compact subset of
$\Delta_\a$ is contained in $\Delta_{\a_n}(\rho)$ for $n$ large
enough. Similarly, every compact subset of $\C\setminus K_\a$ is
contained in $\C\setminus K_{\a_n}$ for $n$ large enough. It follows
that
\[\limsup L_{\a_n}(\rho) \eqdef \bigcap_m \overline{\bigcup_{n\geq m} L_{\a_n}(\rho)}
\subset L_\a(\rho).\] Since $r_\a=\rho$, $\Delta_\a(\rho)=\Delta_\a$
and $L_\a(\rho) = J_\a$. Thus, $\limsup L_{\a_n}(\rho)\subset J_\a$
and \[\area(J_\a)\geq \area(\limsup L_{\a_n}(\rho))\geq
\area\bigl(L_{\a_0}(\rho)\bigr)\cdot \prod (1-\e_n)>0.\]
\qed\par\medskip

\section{The infinitely renormalizable case}

In order to find an infinitely renormalizable quadratic polynomial
with a Julia set of positive area, we need a modification based on
S\o rensen's construction of an infinitely renormalizable quadratic
polynomial with a non-locally connected Julia set.

\begin{proposition}\label{prop_keyinfren0}
There exists a set ${\cal S}$ of bounded type irrationals such that
for all $\a\in {\cal S}$ and all $\eps>0$, there exists
$\a'\in\C\setminus \R$ with
\begin{itemize}
\item $|\a'-\a|<\eps$,

\item $P_{\a'}$ has a periodic Siegel disk with period $>1$ and rotation number in
${\cal S}$ and

\item $\area(K_{\a'})\geq (1-\eps)\area(K_\a)$.
\end{itemize}
\end{proposition}

\begin{figure}[htbp]
\centerline{ \scalebox{.5}{\includegraphics{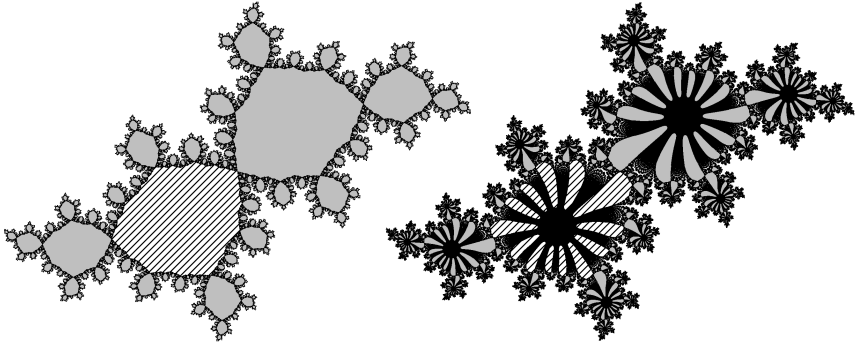}}}
\caption{Two filled-in Julia sets $K_{\a}$ and $K_{\a'}$, with $\a'$
a well-chosen perturbation of $\a$ as in proposition
 \ref{prop_keyinfren0}. This proposition
asserts that if $\a$ and $\a'$ are chosen carefully enough,
$P_{\a'}$ has a periodic Siegel disk and the loss of measure from
$K_{\a}$ to $K_{\a'}$ is small. Left: we hatched the fixed Siegel
disk. Right: we hatched the cycle of Siegel disks. }
\end{figure}

\proof We can choose ${\cal S}={\cal S}_N$ with $N$ large enough (in
order to be able to apply Inou and Shishikura techniques). The proof
essentially goes as in the Cremer case

Given $\a\in {\cal S}$, we let $p_k/q_k$ be its approximants, and we
consider the functions of explosion $\chi_k$ given by proposition
\ref{prop_chi1}. If $\a'$ belongs to the disk centered at $p_k/q_k$
with radius $1/q_k^3$, the set
\[{\cal C}_k(\a')\eqdef \chi_k\left\{ \sqrt[q_k]{\a_k-p_k/q_k}\right\}\]
is a cycle of $P_{\a'}$. Its multiplier is $e^{2i\pi \theta_k(\a')}$
with $\theta_k:D(p_k/q_k,1/q_k^3)\to \C$ a non-constant holomorphic
function which vanishes at $p_k/q_k$.

We consider a sequence $(\a_n)$ converging to $\a$ so that
\begin{itemize}
\item $\ds \limsup_{n\to +\infty}\sqrt[q_n]{\a_n-p_n/q_n}=+\infty$
and

\item $\theta_n(\a_n)\eqdef [A_n,N,N,N,\ldots]$ with
\[\lim_{n\to +\infty}\sqrt[q_n]{A_n}=+\infty\quad\text{and}\quad  \lim_{n\to
+\infty}\sqrt[q_n]{\log A_n}=1.\]
\end{itemize}

We control the shape of the cycle of Siegel disk as in the Cremer
case. For all $\rho<1$ and all $n$ sufficiently large, the cycle of
Siegel disks contains the $\chi_n\bigl(Y_n(\rho)\bigr)$ with
\[Y_n(\rho) \eqdef \left \{z\in \C ~;~\frac{z^{q_n}-\e_n}{z^{q_n}}
\in D(0,s_n)\right\}\quad\text{with}\quad s_n\eqdef
\frac{\rho^{q_n}-|\e_n|}{\rho^{q_n}}.\] For this purpose, we work in
the coordinate given by $\chi_n$ and compare the dynamics of the
conjugated map to the flow of a vector field.

We control the post-critical set of $P_{\a_n}$ via Inou-Shishikura's
techniques.

We then control the loss of area as in the Cremer case.
\qed\par\medskip

\begin{definition}
For $c\in \C$, we denote by $Q_c$ the quadratic polynomial
$Q_c:z\mapsto z^2+c$. With an abuse of notations, we denote by $K_c$
its filled-in Julia set and by $J_c$ its Julia set. We denote by $M$
the Mandelbrot set, i.e. the set of parameters $c$ for which $K_c$
is connected.
\end{definition}

The previous proposition can be restated as follows.

\begin{proposition}\label{prop_areainfren2}
Assume $P_c$ has a fixed Siegel disk with rotation number in ${\cal
S}$. Then, for all $\eps>0$, there exists $c'$ such that
\begin{itemize}
\item $|c'-c|<\eps$,

\item $P_{c'}$ has a periodic Siegel disk with period $>1$ and rotation number in ${\cal S}$ and

\item $\area(K_{c'})>(1-\eps)\area(K_c)$.
\end{itemize}
\end{proposition}

In fact, such a $c$ is on the boundary of the main cardioid of $M$
and the proof we proposed yields a $c'$ which is on the boundary of
a satellite component of the main cardioid of $M$.

Using the theory of quadratic-like maps introduced by Douady and
Hubbard \cite{dh2}, we can transfer this statement to perturbations
of quadratic polynomials having periodic Siegel disks. We will use
the notions of renormalization and tuning  (see for example
\cite{ha}).

If $0$ is periodic of period $p$ under iteration of $Q_{c_0}$, then
$c_0$ is the center of a hyperbolic component $\Omega$ of the
Mandelbrot set. This component $\Omega$ has a root: the parameter
$c_1\in \partial \Omega$ such that $Q_{c_1}$ has an indifferent
cycle with multiplier $1$. In addition, there exist
\begin{itemize}
\item a compact set
$M'\subset M$ such that $\partial M'\subset \partial M$,

\item a simply connected neighborhood $\Lambda$ of $M'\setminus
\{c_1\}$,

\item a continuous map $\chi:\Lambda\cup\{c_1\}\to \C$ and

\item two families of open sets $(U'_\lambda)_{\lambda\in \Lambda}$ and $(U_\lambda)_{\lambda\in \Lambda}$,
\end{itemize}
such that
\begin{itemize}
\item
$\ds \bigl(f_\lambda\eqdef Q_\lambda^{\circ p}:U'_\lambda\to
U_\lambda\bigr)_{\lambda\in \Lambda}$ is an analytic family of
quadratic-like maps

\item for all $\lambda\in \Lambda$, $f_\lambda$ is hybrid conjugate to
$Q_{\chi(\lambda)}$,

\item the Julia set of $f_\lambda$ is connected if and only if $\lambda\in
M'$ and

\item $\chi:M'\to M$ is a homeomorphism (sending $c_0$ to $0$
and $c_1$ to $1/4$).
\end{itemize}
We denote by $c_0\perp\cdot:M\to M'$ the homeomorphism
$(\chi|_{M'})^{-1}$. We say that $c_0\perp c$ is $c_0$ is tuned by
$c$ and that $\ds \bigl(f_\lambda\eqdef Q_\lambda^{\circ
p}:U'_\lambda\to U_\lambda\bigr)_{\lambda\in \Lambda}$ is a
Mandelbrot-like family centered at $c_0$.

\begin{proposition}\label{prop_areainfren3}
Assume $0$ is periodic under iteration of $Q_{c_0}$ and $c'\in M\to
c\in M$ with $\area(K_{c'})\to \area(K_c)$. Then
\[\area(K_{c_0\perp c'})\to \area(K_{c_0\perp c}).\]
\end{proposition}

\proof Let $p$ be the period of $0$ under iteration of $Q_{c_0}$ and
let $\bigl(f_\lambda\eqdef Q_\lambda^{\circ p}:U'_\lambda\to
U_\lambda\bigr)_{\lambda\in \Lambda}$ be a Mandelbrot-like family
centered at $c_0$.

Let $\phi_{c'}:U_{c_0\perp c'}\to \C$ be hybrid conjugacies. As
$c'\to c$, the modulus of the annulus $U_{c_0\perp c'}\setminus
\overline U'_{c_0\perp c'}$ is bounded from below. So, the
$\phi_{c'}$ can be chosen to have a uniformly bounded quasiconformal
dilatation. It follows that if $c'\in M\to c\in M$ with
$\area(K_{c'})\to \area(K_c)$, we have
\[\area\left(\phi_{c'}^{-1}(K_{c'})\right)\underset{c'\to c}\longrightarrow
\area\left(\phi_c^{-1}(K_{c})\right).\] It follows easily that
$\area(K_{c_0\perp c'})\to \area (K_{c_0\perp c})$ since almost
every point in $K_{c_0\perp c}$ has an orbit terminating in
$\phi_{c}^{-1}(K_{c})$. \qed\par\medskip

\noindent{\it Proof of theorem \ref{theo_areainfren}.}~ If $P_c$ has
a periodic Siegel disk then $c$ is on the boundary of a hyperbolic
component with center $c_0$. We denote by $\Omega_c$ this hyperbolic
component and we set $M_c\eqdef c_0\perp M$.

We will denote by $S$ the image of ${\cal S}$ by the map $\a\mapsto
e^{2i\pi\a}/2-e^{4i\pi\a}/4$. Then, $c\in S$ if and only if $P_c$
has a fixed Siegel disk with rotation number in ${\cal S}$.
Moreover, $P_c$ has a periodic Siegel disk with rotation number in
${\cal S}$ if and only if $c=c_0\perp s$ with $c_0$ the center of
the hyperbolic component containing $c$ in its boundary and $s\in
S$.

It follows from propositions \ref{prop_areainfren2} and
\ref{prop_areainfren3} that if $Q_c$ has a periodic Siegel disk with
rotation number in ${\cal S}$, then for all $\e>0$, we can find
$c'\in M_c\setminus \overline\Omega_c$ such that
\begin{itemize}
\item $|c'-c|<\eps$,

\item $P_{c'}$ has a periodic Siegel disk with rotation number in ${\cal S}$ and

\item $\area(K_{c'})>(1-\eps)\area(K_c)$.
\end{itemize}

Let us choose a parameter $c_0\in S$ and a sequence of real number
$\e_n$ in $(0,1)$ such that $\prod(1-\e_n)>0$. We can construct
inductively a sequence $(c_n)$ such that
\begin{itemize}
\item $(c_n)$ is a Cauchy sequence that converges to a parameter $c$,

\item $Q_{c_n}$ has a periodic Siegel disk with rotation number in
${\cal S}$,

\item for $n\geq 1$, $c_n\in M_{c_{n-1}}\setminus \overline
\Omega_{c_{n-1}}$ and

\item $\area(K_{c_n})>(1-\eps_n)\area(K_{c_{n-1}})$.
\end{itemize}

Then, $P_c$ is infinitely renormalizable (it is in the intersection
of the nested copies $M_{c_n}$). Thus, $J_c=K_c=\lim K_{c_n}$.
Finally,
\[\area(J_c)=\area(K_c) \geq \area(K_{c_0})\cdot\prod (1-\e_n) >0.\]
\qed\par\medskip

\appendix

\section{Parabolic implosion and perturbed petals\label{appendix_fatou}}

The notations used in this appendix are those of section
\ref{subsubperturbedfatou}. We postponed the proof of the following
lemma to this appendix.

\begin{lemma}\label{lemma_perturbedpetal}
If $R>0$ and $K>0$ are sufficiently large, then for $n$ large
enough:
\begin{enumerate}
\item\label{item_1perturbedpetal}
$\Phi^n(\Omega^n)$ contains the vertical strip
\[U^n\eqdef \bigl\{w\in \C~;~R<\re(w)<1/\a_n-R\bigr\},\]

\item\label{item_2perturbedpetal} $\tau_n$ is injective on $\Pet^n\eqdef
(\Phi^n)^{-1}(U^n)$ and

\item\label{item_4perturbedpetal} there is a branch of argument defined on $\tau_n(\Pet^n)$ such that
\[\sup_{z\in \tau_n(\Pet^n)} \arg(z) - \inf_{z\in \tau_n(\Pet^n)} \arg(z) <K.\]
\end{enumerate}
\end{lemma}

\begin{proof}
As in \cite{sh}, the argument consists in comparing the Fatou
coordinate $\Phi^n$ to the Fatou coordinate $\Psi^n$ of the time one
map of the vector field $\zeta_n$ defined on ${\cal D}_n$ by
\[\zeta_n= \zeta_n(w)\frac{\partial}{\partial w}\eqdef
\bigl(F_n(w)-w\bigr)\frac{\partial}{\partial w}.\] In other words,
set $\ds w_n\eqdef \frac{1}{2\alpha_n}$ and let $\Psi^n:\Omega_n\to
\C$ be defined by
\[\Psi^n(w) = \Phi^n(w_n) +
\int_{w_n}^w \frac{du}{F_n(u)-u}.\]

\medskip
\noindent{\bf Claim 1}. Increasing $R_1$ if necessary, there is a
constant $C>0$ such that for all $n$ sufficiently large
\[\sup_{w\in \Omega^n} \bigl|\Phi^n(w)-\Psi^n(w)\bigr|< C.\]

\medskip
\noindent{\em Proof of Claim 1}. According to Prop. 2.6.2 in
\cite{sh}, there are constants $R$ and $C$ such that for all
sufficiently large $n$ and for all $w\in \Omega^n$ with
$d(w,\partial \Omega^n)\geq R$, we have
\[\bigl|(\Phi^n)'(w) - (\Psi^n)'(w)\bigr|\leq
C\left(\frac{1}{d(w,\partial\Omega^n)^2} +
\bigl|F_n'(w)-1\bigr|\right).\] We will first show that we can get
rid of $\bigl|F_n'(w)-1\bigr|$. Set
\[G_n(w)\eqdef F'_n(w)-1\quad\text{and}\quad
S_n(w)\eqdef \left(\frac{\pi\a_n}{\sin(\pi\a_nw)}\right)^2.\] Those
functions are $1/\a_n$ periodic. On the one hand, as $n\to +\infty$,
\begin{itemize}
\item the functions $G_n$ are uniformly bounded by $1/4$ on $\partial
\Omega^n$ and

\item the sequence $(S_n)$ converges uniformly to $w\mapsto 1/w^2$ on $\partial \Omega^n$, and thus,
the functions $S_n$ are uniformly bounded away from $0$ on $\partial
\Omega^n$.
\end{itemize}
As a consequence, the functions $G_n/S_n$ are uniformly bounded on
$\partial \Omega^n$. On the other hand, as $\im(w)\to \pm \infty$,
$G_n(w)\to 0$. Thus, in $\C/\frac{1}{\a_n}\Z$, $G_n$ has removable
singularities at $\pm i\infty$ and vanishes at those points. Since
in $\C/\frac{1}{\a_n}\Z$, $S_n$ has simple zeros at $\pm i\infty$,
the function $G_n/S_n$ has removable singularities at $\pm i\infty $
in $\C/\frac{1}{\a_n}\Z$. It follows that from the maximum modulus
principle that there is a constant $C_1$ such that for all
sufficiently large $n$ and all $w\in \Omega^n$, we have
\[\bigr|F'_n(w)-1\bigr|\leq
C_1\left|\frac{\pi\a_n}{\sin(\pi\a_nw)}\right|^2.\] Note that there
is a constant $C_2>0$ such that
\[\forall w\in \C,
\quad  d(w,\Z)\leq C_2\bigl|\sin(\pi w)\bigr|.\] Indeed, the
quotient $\ds\frac{d(w,\Z)}{\bigl|\sin(\pi w)\bigr|}$ extends
continuously to $(\C/\Z)\cup \{\pm i\infty\}$ which is compact. It
follows that for all $w\in \Omega^n$,
\[\left|\frac{\pi\a_n}{\sin(\pi \a_n w)}\right|^2\leq \frac{C_2^2\pi^2|\a_n|^2}{d(\a_nw,\Z)^2}\leq
\frac{C_2^2\pi^2}{d(w,\partial \Omega^n)^2}.\] Thus, there is a
constant $C'$ such that for all sufficiently large $n$ and for all
$w\in \Omega^n$ with $d(w,\partial \Omega^n)\geq R$, we have
\[\bigl|(\Phi^n)'(w) -(\Psi^n)'(w)\bigr|\leq \frac{C'}{d(w,\partial
\Omega^n)^2}.\] Taking $R\geq 1$ and replacing $R_1$ by $R_1+\sqrt
2R$, this can be rewritten as: there is a constant $C$ such that for
all sufficiently large $n$ and for all $w\in \Omega^n$
\[\bigl|(\Phi^n)'(w) - (\Psi^n)'(w)\bigr|\leq\frac{C'}{\bigl(1+d(w,\partial
\Omega^n)\bigr)^2}.\]

Let us now assume $n$ is sufficiently large, so that
\[X_n\eqdef \frac{1}{2\a_n}-R_1>0.\]
Then, $\ds w_n\eqdef \frac{1}{2\a_n}$ belongs to $\Omega^n.$ Fix
$w\eqdef w_n+x+iy\in \Omega^n$. Note that
\[|x|< X_n+|y|\quad\text{and}\quad
d(w,\partial \Omega^n)> \sqrt2\bigl(X_n+|y|-|x|\bigr).\] It follows
that
\begin{align*}
\bigl|\Phi^n(w)-\Psi^n(w)\bigr|&\leq
\int_{[w_n,w_n+iy]\cup [w_n+iy,w]}\frac{C'|du|}{\bigl(1+d(u,\partial\Omega^n)\bigr)^2}\\
&\leq \int_0^{+\infty} \frac{C'ds}{\bigl(1+\sqrt 2(X_n+s)\bigr)^2} +
\int_0^{X_n+|y|} \frac{C'dt}{\bigl(1+\sqrt 2(X_n+|y|-t)\bigr)^2}\\
&\leq 2C'.\end{align*} This completes the proof of Claim 1. \qed
\medskip

\medskip
\noindent{\bf Claim 2}. The map $\Psi^n$ is univalent on $\Omega^n$,
$\Psi^n(\Omega^n)$ contains the vertical strip
\[V^n\eqdef \bigl\{w\in \C~;~\re\bigl(\Psi^n(R_1)\bigr)<\re(w)<
\re\bigl(\Psi^n(1/\a_n-R_1)\bigr)\bigr\}\] and $\tau_n$ is injective
on ${\cal Q}^n\eqdef (\Psi^n)^{-1}(V^n)$.

\medskip
\noindent{\em Proof of Claim 2}. Note that $\Psi^n$ is a
straightening map for the vector field $\zeta_n$:
\[(\Psi^n)_*\zeta_n =
\ds\frac{\partial}{\partial w}.\] Since $F_n(w)-w\in D(1,1/4)$ on
$\Omega^n$, the trajectories of the vector field $\zeta_n$ are
curves which enter $\Omega^n$ through its left boundary and exit
$\Omega^n$ through the right boundary. In particular, no trajectory
is periodic. Since two distinct trajectories cannot intersect, the
map $\Psi^n$ is injective.

Observe that for $w\in \partial \Omega^n$,
\[\arg\bigl((\Psi^n)'(w)\bigr) = -\arg\bigl(F_n(w)-w\bigr)\in
\bigl]-\arcsin(1/4),\arcsin(1/4)\bigr[\subset\bigl]-\pi/12,\pi/12\bigr[.\]
Integrating $(\Psi^n)'(w)$ along $\partial \Omega^n$, we conclude
that
\[\frac{2\pi}{3}<\arg\bigl(\Psi^n(w)-\Psi^n(R_1)\bigr)<\frac{4\pi}{3}\]
on the left boundary of $\Omega^n$ and that
\[-\frac{\pi}{3}<\arg\bigl(\Psi^n(w)-\Psi^n(1/\a_n-R_1)\bigr)<\frac{\pi}{3}\]
on the right boundary of $\Omega^n$. This proves that
$\Psi^n(\Omega^n)$ contains the vertical strip $V^n$.

Assume by contradiction that $\tau_n$ is not injective on $V^n$.
Then, there is an integer $k\in \Z\setminus \{0\}$ and a point $w\in
V^n$ such that $w+k/\a_n$ is in $V^n$. Note that $V^n$ is a union of
trajectories for the rotated vector field $i\zeta_n$. As $w$ runs
along those trajectories, the imaginary part of $w$ increases from
$-i\infty$ to $+i\infty$. In particular, every trajectory intersects
$\R$. Since for all $w\in {\cal D}_n$, we have $i\zeta_n(w) =
i\zeta_n(w+1/\a_n)$, the trajectory for $i\zeta_n$ passing through
$w+k/\a_n$ is obtained from the trajectory passing through $w$ by
translation by $k/\a_n$. This is not possible since the intersection
of those trajectories with $\R$ is contained in $\Omega^n\cap \R =
\left]R_1,1/\a_n-R_1\right[$. This completes the proof of Claim 2.
\qed
\medskip

Let us now come to the proof of parts (\ref{item_1perturbedpetal})
and (\ref{item_2perturbedpetal}) of lemma
\ref{lemma_perturbedpetal}. Assume $n$ is sufficiently large, so
that
\[\sup_{w\in \Omega^n}\bigl|\Phi^n(w)-\Psi^n(w)\bigr|\leq C.\]
Then, $\Phi^n({\cal Q}^n)$ contains the vertical strip
\[\bigl\{w\in \C~;~\re\bigl(\Psi^n(R_1)\bigr)+C<\re(w)<
\re\bigl(\Psi^n(1/\a_n-R_1)\bigr)-C\bigr\}.\] Note that
\[\Psi^n(R_1) =\Phi^n(R_1) +{\cal O}(1) = {\cal O}(1)\]
and
\[\Psi^n(1/\a_n-R_1) = \Phi^n(1/\a_n-R_1) +{\cal O}(1) =
1/\a_n+{\cal O}(1).\] Thus, if $R$ is large enough and if $n$ is
sufficiently large, then $\Phi^n({\cal Q}^n)$ contains the vertical
strip
\[U^n\eqdef \bigl\{w\in \C~;~R<\re(w)<
1/\a_n-R\bigr\}.\] Since $\tau_n$ is injective on ${\cal Q}^n$, this
proves parts (\ref{item_1perturbedpetal}) and
(\ref{item_2perturbedpetal}) of lemma \ref{lemma_perturbedpetal}.

\medskip
Let us now come to the proof of part (\ref{item_4perturbedpetal}) of
lemma \ref{lemma_perturbedpetal}. Note that $\tau_n$ sends the
segment $\left]0,1/\a_n\right[$ to the perpendicular bisector of the
segment $[0,\sigma_n]$. It sends the lower half-plane $\H^-\eqdef
\bigl\{w\in \C~;~\im(w)<0\bigr\}$ in the half-plane $\bigl\{z\in
\C~;~|z|<|z-\sigma_n|\bigr\}$. It is a universal covering from the
upper half-plane
\[\H^+\eqdef \bigl\{w\in \C~;~\im(w)>0\bigr\}\] to the punctured
half-plane $\bigl\{z\in \C~;~0<|z|<|z-\sigma_n|\bigr\}$, with
covering transformation group generated by the translation
$T_n:w\mapsto w+1/\a_n$. It sends the lines
\[L_k\eqdef \left\{w\in
\C~;~\re(w)=\frac{2k+1}{2\a_n}\right\},\quad k\in \Z\] to the
segment $\left]0,\sigma_n\right[$.

We must show that there is a constant $M$ such that for $n$ large
enough, $\Pet^n\cap \H^+$ is contained in the vertical strip
\[\left\{w\in
\C~;~-\frac{M}{\a_n}<\re(w)<\frac{M}{\a_n}\right\}.\] For all $w\in
\Pet^n$, we have \[R\leq \re\bigl(\Phi^n(w)\bigr)\leq
\frac{1}{\a_n}-R.\] It is therefore enough to show that
\[\sup_{w\in \Omega^n\cap \H^+}\bigl|\Phi^n(w)-w\bigr|
={\cal O}\left( \frac{1}{\a_n}\right)\] or equivalently that
\[\sup_{w\in \Omega^n\cap \H^+}\bigl|\Psi^n(w)-w\bigr| ={\cal
O}\left( \frac{1}{\a_n}\right).\] Note that
$\ds\frac{1}{F_n(w)-w}-1$ is periodic of period $1/\a_n$, bounded by
$1/3$ in $\Omega^n$ and tends to $0$ as $\im(w)$ tends to $+\infty$.
It follows from the maximum modulus principle that
\[\left|\frac{1}{F_n(w)-w}-1\right|<\frac{1}{3}\cdot \left(\inf_{w\in
\partial (\Omega^n\cap \H^+)}|e^{2i\pi \a_n w}|\right)\cdot |e^{2i\pi \a_n w}|
\leq C e^{-2\pi \a_n \im(w)}\] for some constant $C$ which does not
depend on $n$. If $w\eqdef R+x+iy\in \Omega^n\cap \H^+$, then
$|x|<y+ 1/\a_n$. So
\begin{align*}
\sup_{w\in \Omega^n\cap \H^+} \bigl|\Psi^n(w)-w\bigr| &\leq
\bigl|\Psi^n(R)-R\bigr| \\
&\quad +\sup_{y>0\atop |x|<y+1/\a_n }\left(\int_0^{y} Ce^{-2\pi \a_n
t}dt + \int_0^{|x|} Ce^{-2\pi \a_n
y}dt\right) \\
&= \bigl|\Psi^n(R)-R\bigr|+\frac{C}{\a_n}
\left(\frac{1}{2\pi}+1\right) ={\cal O}\left(\frac{1}{\a_n}\right).
\end{align*}
This completes the proof of part (\ref{item_4perturbedpetal}) of
lemma \ref{lemma_perturbedpetal}.
\end{proof}

\end{document}